\input amstex
\documentstyle{amsppt}
\magnification1200
\tolerance=10000
\overfullrule=0pt
\def\n#1{\Bbb #1}
\def\p{\Bbb C_{\infty}}

\def\cb{\Cal B}
\def\ci{\Cal I}

\def\cl{\hbox{ cl }}

\def\Gal{\hbox{Gal }}

\def\Hom{\hbox{Hom}}

\def\End{\hbox{End}}

\def\Ker{\hbox{Ker }}
\def\Lie{\hbox{Lie}}

\def\Pic{\hbox{Pic}}

\def\diag{\hbox{ diag }}

\def\s{\sigma}

\def\e11{I_{11}}

\def\ga{\goth A}

\def\Im{\hbox{Im }}
\def\vf{\varphi}

\def\ve{\varepsilon}

\def\vf{\varphi}
\def\cl{\Cal L}
\def\de{\delta}

\def\vi{v_\infty}
\def\ga{\gamma}

\def\la{\lambda}

\def\be{\beta}

\def\th{\theta}
\def\ze{\zeta}
\def\al{\alpha}
\def\r{\n R_\infty}
\def\si{\sigma}
\def\oo{\Omega}
\def\om{\omega}
\def\g{\goth }
\def\vk{\varkappa}

\NoRunningHeads
\topmatter
\title Introduction to Anderson t-motives: a survey
\endtitle
\author A. Grishkov, D. Logachev\footnotemark \footnotetext{E-mail: logachev94{\@}gmail.com\phantom{*********************************************************}}
\endauthor
\address
First author: Instituto de Matem\'atica e estatistica, 
Universidade de S\~ao Paulo. Rua de Mat\~ao 1010, CEP 05508-090, S\~ao Paulo, Brasil, 
\medskip
and Sobolev Institute of Mathematics, Omsk, Russia.
\medskip
Second author: Departamento de Matem\'atica, Universidade Federal do Amazonas, Manaus, Brasil.
\endaddress
\NoRunningHeads
\keywords Anderson t-motives; Lattices of Anderson t-motives  \endkeywords
\subjclass 11G09 \endsubjclass
\abstract This is a survey on Anderson t-motives --- high-dimensional generalizations of Drinfeld modules. They are the function field analogs of abelian varieties with multiplication by an imaginary quadratic field. We describe their lattices, their groups $H^1$ and $H_1$, their tensor products, the duality functor and the duality theorem, endomorphisms of Drinfeld modules in finite characteristic, and their L-functions of a certain type. Further on, we introduce the notion of affine equations, $T$-divisible $\Bbb F_q[[T]]$-modules, holonomic sequences in the function field case, analogs of Siegel matrices as elements of flag varieties, and some other notions (to be continued). Many examples of explicit calculations are given, some elementary research problems are stated. Some results (Sections 16; 19) are apparently new. 
\endabstract
\endtopmatter
\document
\rightline{Preliminary version}
\medskip
{\bf Acknowledgements}
\medskip
This survey grew out of the discussion with students of a course given by the second author in POSTECH, Pohang, Korea, in Spring 2024. The second author is grateful to YoungJu Choie for her kind invitation to POSTECH and for the possibility to give this course. 
\medskip
The first author was supported by FAPESP (grant 2018/23690-6), 
CNPq (grants 307593/2023-1 and the universal grant), and in accordance
 with the state task of the IM SB RAS, project FWNF-2022-003.
\medskip
The second author is grateful to Vladimir Drinfeld for numerous explanations on the subject. The authors are grateful to Urs Hartl and to anonymous reviewers of their papers for indicating simple proofs of some theorems. 
\medskip

{\bf  Contents}
\medskip
\settabs 21 \columns
\+0. Introduction &&&&& &&&&& &&&&& &&&&& 3\cr
\medskip
\+1. Abelian varieties in characteristic 0&&&&& &&&&& &&&&& &&&&&8\cr
\medskip
\+&1.5. Tate modules&&&&& &&&&& &&&&& &&&&8\cr
\medskip
\+&1.7. Reduction&&&&& &&&&& &&&&& &&&&9\cr
\medskip
\+&1.8. Abelian varieties with MIQF&&&&& &&&&& &&&&& &&&&10\cr
\medskip
\+2. Abelian varieties in characteristic $p$&&&&& &&&&& &&&&& &&&&&11\cr
\medskip
\+3. Initial rings and fields&&&&& &&&&& &&&&& &&&&&11\cr
\medskip
\+4. Lattices in functional fields&&&&& &&&&& &&&&& &&&&&12\cr
\medskip
\+5. Anderson t-motives&&&&& &&&&& &&&&& &&&&&15\cr
\medskip
\+&5.14. Change of basis&&&&& &&&&& &&&&& &&&&17\cr
\medskip
\+&5.15. Pure t-motives&&&&& &&&&& &&&&& &&&&19\cr
\medskip
\+&5.17. Tensor products&&&&& &&&&& &&&&& &&&&21\cr
\medskip
\+&5.18. Exterior powers of Drinfeld modules&&&&& &&&&& &&&&& &&&&22\cr
\medskip
\+&5.20. Analogy with rings of differential operators&&&&& &&&&& &&&&& &&&&23\cr
\medskip
\+6. Action of $\n F_q[T]$ on $\p^n$ associated to Anderson t-motives&&&&& &&&&& &&&&& &&&&&23\cr
\medskip
\+&6.11. Torsion points and Tate modules&&&&& &&&&& &&&&& &&&&26\cr
\medskip
\+7. Drinfeld modules over finite fields; reductions&&&&& &&&&& &&&&& &&&&&27\cr
\medskip
\+&7.9. Valuations over $\g p$&&&&& &&&&& &&&&& &&&&29\cr
\medskip
\+&7.10. Valuations over $\infty$&&&&& &&&&& &&&&& &&&&30\cr
\medskip
\+&7.13. Weil numbers&&&&& &&&&& &&&&& &&&&31\cr
\medskip
\+&7.16. Case of Anderson t-motives&&&&& &&&&& &&&&& &&&&31\cr
\medskip
\+&7.17. Reductions&&&&& &&&&& &&&&& &&&&32\cr
\medskip
\+8. Lattices of Drinfeld modules&&&&& &&&&& &&&&& &&&&&34\cr
\medskip
\+9. Lie(M)&&&&& &&&&& &&&&& &&&&&35\cr
\medskip
\+10. Lattice associated to some Anderson t-motives, and a principal exact \cr
\+&sequence&&&&& &&&&& &&&&& &&&&37\cr
\medskip
\+11. $\th$-shift&&&&& &&&&& &&&&& &&&&&42\cr
\medskip
\+12. Dualities and pairings&&&&& &&&&& &&&&& &&&&&45\cr
\medskip
\+&12.1. Dual of a t-motive&&&&& &&&&& &&&&& &&&&45\cr
\medskip
\+&12.3. Definitions, isomorphisms and pairings for $H_1(M)$, $H_1(M')$, &&&&& &&&&& &&&&& &&&&\cr
\+&$H^1(M)$, $H^1(M')$&&&&& &&&&& &&&&& &&&&46\cr
\medskip
\+&12.7. List of formulas for pairings and isomorphisms&&&&& &&&&& &&&&& &&&&52\cr
\medskip
\+&12.8. Connections between $\Psi(M), \ \Phi(M), \ \Psi(M'), \ \Phi(M')$&&&&& &&&&& &&&&& &&&&53\cr
\medskip
\+&12.9. Duality of lattices&&&&& &&&&& &&&&& &&&&53\cr
\medskip
\+&12.11. Theorem of Anderson, complete form&&&&& &&&&& &&&&& &&&&55\cr
\medskip
\+13. The lattice functor&&&&& &&&&& &&&&& &&&&&56\cr
\medskip
\+&13.7. Surjectivity of the lattice functor&&&&& &&&&& &&&&& &&&&57\cr
\medskip
\+14. Case $N\ne0$&&&&& &&&&& &&&&& &&&&&58\cr
\medskip
\+15. L-functions&&&&& &&&&& &&&&& &&&&&60\cr
\medskip
\+&15.5. Example: Explicit calculations&&&&& &&&&& &&&&& &&&&64\cr
\medskip
\+&15.6. Global case&&&&& &&&&& &&&&& &&&&66\cr
\medskip
\+&15.7. $v_\infty$ of eigenvalues of pure t-motives&&&&& &&&&& &&&&& &&&&68\cr
\medskip
\+16. Affine equations&&&&& &&&&& &&&&& &&&&&69\cr
\medskip
\+&16.11. Small solutions&&&&& &&&&& &&&&& &&&&72\cr
\medskip
\+&16.12. Holonomic sequences&&&&& &&&&& &&&&& &&&&72\cr
\medskip
\+17. The Drinfeld "upper half plane"&&&&& &&&&& &&&&& &&&&&73\cr
\medskip
\+&17.8. Sets $D(n,x)$ and the Bruhat-Tits building $\Cal T$ of $\oo$&&&&& &&&&& &&&&& &&&&75\cr
\medskip
\+18. Relations between $\vf$-modules and Galois modules&&&&& &&&&& &&&&& &&&&&75\cr
\medskip
\+19. Flag varieties related to lattices of t-motives&&&&& &&&&& &&&&& &&&&&80\cr
\medskip
\+&19.1. Definitions related to Grassmannians&&&&& &&&&& &&&&& &&&&80\cr
\medskip
\+&19.2. Map $\de$ from the set of bases $\cb$ to the Grassmannian $Gr(r-n,r)$ &&&&& &&&&& &&&&& &&&&81\cr
\medskip
\+&19.3. Action of $GL_r(\n F_q[\th])$ on $\cb$ and on $Gr(r-n,r)$&&&&& &&&&& &&&&& &&&&83\cr
\medskip
\+&19.4. t-modular forms&&&&& &&&&& &&&&& &&&&84\cr
\medskip
\+&19.5. Application to the duality&&&&& &&&&& &&&&& &&&&85\cr
\medskip
\+& 19.6. Case of $N\ne0$&&&&& &&&&& &&&&& &&&&85\cr
\medskip
\+&&19.6.1. Properties of the "classical" flag varieties&&&& &&&&& &&&&& &&&&85\cr
\medskip
\+&&19.6.3. Definition of numbers $k_i$ --- invariants of $N$&&&& &&&&& &&&&& &&&&87\cr
\medskip
\+&&19.6.4. t-flag variety&&&& &&&&& &&&&& &&&&88\cr
\medskip
\+&&19.6.5. $tFV(N)$ is an analog of $Gr(r-n,r)$ for the case $N\ne0$&&&& &&&&& &&&&& &&&&91\cr
\medskip
\+&&19.6.6. Description of $\de(l_*)$ via its Siegel set &&&& &&&&& &&&&& &&&&93\cr
\medskip
\+&& 19.6.7. Case of an arbitrary Schubert cell &&&& &&&&& &&&&& &&&&94\cr
\medskip
\+&& 19.6.8. Duality and tensor product &&&& &&&&& &&&&& &&&&96\cr
\medskip
\+99. Generalizations of Anderson t-motives&&&&& &&&&& &&&&& &&&&&97\cr
\medskip
\+References&&&&& &&&&& &&&&& &&&&&98\cr
\medskip
\medskip
{\bf 0. Introduction.}
\medskip
{\bf 0.1.} Drinfeld in [D76] defined new objects --- elliptic modules, now called Drinfeld modules. Anderson in [A86] defined their high-dimensional generalizations --- Anderson t-motives. Their theory forms --- in some sense --- a parallel world to the classical theory of abelian varieties (see Section 6, Table 6.14), like the non-Euclidean geometry (Bolyai, Gauss, Lobachevsky) forms a parallel world to the Euclidean geometry. Namely, Anderson t-motives can be considered as analogs of abelian varieties with MIQF (multiplication by imaginary quadratic fields), see 1.8. This analogy is a source of numerous research papers: a scientist considers a theorem of the theory of abelian varieties and proves its analog for the theory of Anderson t-motives. This is not a routine activity, because the analogy is not complete. As an example of absence of the complete analogy, we can mention that we do not know whether there exists, or not, a 1 -- 1 correspondence between the pure uniformizable t-motives and lattices (for the Drinfeld modules it does exist), the Siegel matrices of lattices are elements not of affine spaces like in the case of abelian varieties, but of Grassmannians and, more generally, of flag varieties (see Section 19) etc. 
%\medskip

Also, there is an important analogy of the theory of Anderson t-motives and the theory of differential operators, see [Mu78] for a source. See Sections 5.20 and 16.12 for details. 
%\medskip

A detailed introduction to the subject can be found in [G96] and [Pa23]. Nevertheless, all these four sources [D76], [A86], [G96], [Pa23] can be too difficult for the beginners. One of the reasons of this fact is the subject itself. Namely, important objects of the theory are elements $T$ and $\th$, see below. Their roles are close one to another, there is a map $\iota$ sending $T$ to $\th$, until now some authors confuse them, identifying them by this map. This happened quite often at earlier stages of development of the theory, and is not satisfactory. In fact, it is necessary to distinguish carefully $T$ and $\th$, like we distinguish the nearby notions of mass and weight in physics. 
\medskip
The purpose of the present survey is to give an elementary introduction to the subject. On the one hand, we indicate the most important ideas omitting technical
details. On the other hand, we give lacking details of proofs of [G96], [A86]. For example, we give a detailed coordinate description of the objects of the principal exact sequence (10.4), namely of $E(M)$ --- the t-module associated to an Anderson t-motive $M$ of dimension $n$ (see below for the definition of $n$), and its Lie group $Lie(M)$. It turns out that although $E(M)$ is isomorphic to $\p^n$ (here $\p$ is a finite characteristic analog of $\n C$, see Section 3), but this isomorphism is not canonical for $n>1$: it depends on a choice of a basis of $M$ over $\p\{\tau\}$, see 5.2.2. Particularly, there is no well-defined $c\cdot x$, where $c\in\p$, $x\in E(M)$. On the other hand, there exists a canonical isomorphism $Lie(M)\to \p^n$ (maybe this is wrong: we have only a canonical $\p$-structure on $Lie(M)$, but not an isomorphism). 
\medskip
These facts are not emphasized in both [A86], [G96], although a non-canonical coordinade description of $E(M)$ (i.e. an isomorphism $E(M)\to \p^n$) is used essentially in the proofs of theorems. This can cause confusion for beginners. The present survey contains a detailed exposition of this subject (Sections 6 and 9). 
\medskip
Further on, we give a detailed proof of coincidence of two definitions of local L-factors of t-motives, which come, on the one hand, from a consideration of a t-motive $M$ as a $\p[T]$-module, and on the other hand, from its consideration as a $\p\{\tau\}$-module (see below for the definitions). We give explicit calculations of these local L-factors for small $r$ by two ways (Section 15.5, cases A and B). Concordence of these two calculations is not a trivial fact. 
\medskip
Also, we give explicit examples and formulas. Hence, this paper can be used as a reference book. For example, we give explicit formulas for duals of the Drinfeld modules, for pairings and isomorphisms between $H^1(M),\ H_1(M),\ H^1(M'),\ H_1(M')$. 
\medskip
We give a simple explicit proof of the duality theorem for the lattice functor, i.e. that the functors of lattice and of duality for t-motives commute, for the case $N=0$. The most important formula in this proof is the relation between a scattering matrix $\Psi$ of a t-motive $M$, and a Siegel matrix of its lattice, see (12.10.5) (apparently this formula fails in [A86]). 
\medskip
{\bf 0.2.} From a formal point of view, Anderson t-motives are simple objects: they are modules (having some specific properties (5.2.1 -- 5.2.3)) over a ring of non-commutative polynomials in two variables over a field. But in order to understand that they are really analogs of abelian varieties, it is necessary to plunge into the theory.
\medskip
The generalizations of the theory of Anderson t-motives are very deep and complicated. Nevertheless, until now there exist elementary problems that are not solved yet. Another purpose of this survey is to indicate some "down-to-earth" problems that can be a research subject for the beginners. For example, calculation of $h^1$, $h_1$ of some Anderson t-motives is a long, but easy problem that will definitely give a result (continuation of [GL21], [EGL]). See Section 13.8. Other elementary research subject is described in Section 15.6.6.
\medskip
Let us describe the analogy between Anderson t-motives and abelian varieties in more details. First of all, while abelian varieties depend on one discrete parameter --- their dimension $g$, Anderson t-motives depend on two parameters --- dimension $n$ and rank $r$. It turns out that Anderson t-motives are analogs not of general abelian varieties, but of abelian varieties with MIQF of signature $(n,r-n)$, see 1.8 below and [GL09]. Surprisingly, this fact is not emphasized in most survey papers on the subject.
\medskip
Continuing the analogy, we have: many objects attached to abelian varieties, for example Tate modules, Galois action on them, lattices, modular curves, L-functions etc.,
also can be attached to Anderson t-motives. 
\medskip
Nevertheless, as we mentioned above, this analogy is far to be complete. Other examples of absence of analogy: there is no functional equation for L-functions of Anderson t-motives; notion of their algebraic rank is not known yet; not all Anderson t-motives have a lattice of "correct" dimension, etc. 
\medskip
{\bf 0.3.} The paper is organized as follows. In Section 1 we give briefly some properties of abelian varieties over number fields, in order to show the analogy. Particularly, in Section 1.8 we give a definition of abelian varieties with MIQF. In Section 2 we show how these properties are modified for abelian varieties over global functional fields. In Section 3 we define $\n F_q[\th], \ \n F_q(\th), \ \n R_\infty, \ \p$ --- analogs of $\n Z, \ \n Q, \ \n R, \ \n C$ for the finite characteristic case, and we consider an explicit calculation in $\p$. In Section 4 we define $\n F_q[\th]$-lattices in vector spaces over $\p$, and their Siegel matrices. We prove an important theorem 4.2 that there are only finitely many elements of a lattice in any bounded domain. 
\medskip
Section 5 gives a definition and properties of Anderson t-motives, considered as modules over Anderson ring, having some specific properties. We give explicit formulas for them and for Drinfeld modules as their particular cases. We give a definition of pure t-motives, we give its explicit matrix form, we give an explicit proof that Drinfeld modules are pure. Further, we give some examples of equations on pure and non-pure t-motives with $n=2$ and small $r$. Also, we define the tensor product of t-motives. Finally, we show in (5.20) an analogy between the Anderson ring and a ring of differential operators. 
\medskip
In Section 6 we consider t-modules $E(M)$ associated to Anderson t-motives $M$ (historically, it was the first definition). We consider identifications of $E(M)$ with $\p^n$ --- there is no canonical identification unless $n=1$. Also, we define the Tate modules of $M$.
\medskip
In Section 7 we give --- without proofs, we refer to [G96], [Pa23], --- the description of endomorphism rings and fields of Drinfeld modules over finite fields. We give a diagram of involved endomorphism fields, we give an explicit calculation of characteristic polynomial of the Frobenius endomorphism for some simple cases. By the way, to find explicit formulas for coefficients of the Frobenius endomorphism for some more complicated cases would be a good exercise for a student. Further, we consider the behaviour of valuations in these field extensions and the values of valuations of the Frobenius endomorphism. Finally, we give briefly some properties of reductions of Anderson t-motives. We give an explicit calculational proof of a theorem that the kernels of the reduction map for a t-motive $M$ and its dual $M'$ are in duality, for a particular case when $M$ is a Drinfeld module and the reduction is in $\th$, see (7.17). 
\medskip
In Section 8 we define lattices of Drinfeld modules, and we give a detailed explicit calculation of the lattice map for the Carlitz module. In Section 9 we define $Lie(M)$, its coordinate description, the action of operator $T$ on it, and the map $exp: \ Lie(M)\to E(M)$. We show that $Lie(M)$ is canonically isomorphic to $\p^n$. 
\medskip
Section 10 gives three equivalent descriptions of the principal exact sequence (10.4) associated to an Anderson t-motive: 
$$0\to L(M)\to Lie(M)\overset{exp}\to{\to} E(M)$$
where $L(M)$ (denoted also by $H_1(M)$ --- the homology group of $M$) is a lattice of $M$. 
\medskip
An important Theorem 10.5 (see also Theorem 12.11.2 for its generalization) affirms, first, that three descriptions of (10.4) are equivalent, and second, that 
$$h_1(M):=\dim (L(M))=r\iff exp \hbox{ is surjective. }$$

Also, we give some formulas that relate elements of terms of this exact sequence. In particular, we define $\partial_\la(f)$, where $\la\in Lie(M)$, $f\in M$ - it enters in the statement of the below Theorem 11.9. 
\medskip
Section 11 is a purely technical one. It describes the $\th$-shift --- a tool used in the statement and the proof of Theorem 11.9. It is the most important among the theorems relating elements of the terms of the principal exact sequence. This theorem will be used later for the proof of the duality theorem for lattices. 
\medskip
In Section 12, first, we define duality for Anderson t-motives, as a particular case of the Hom functor (case $N=0$). We give explicit formulas for the dual of a Drinfeld module. Further, we define $H^1(M)$ --- the cohomology group of $M$, and we give an example of calculation of $H^1$ of the Carlitz module (a Drinfeld module having $r=n=1$). We give formulas for pairings and isomorphisms between $H_1$ and $H^1$ of $M$ and its dual (sections 12.6, 12.7). We define the notion of the dual lattice (12.9), and we give a proof of the main theorem: the lattice of the dual t-motive is the dual of the lattice of the initial t-motive (case $N=0$). Finally, we formulate the principal theorem on uniformizable t-motives, i.e. t-motives having $h_1=h^1=r$ (Theorem 12.11.2). 
\medskip
In Section 13 we describe some properties of the lattice functor of t-motives: cases of duals of Drinfeld modules, local surjectivity of the lattice map, its injectivity for mixed t-motives ([HJ]), its non-injectivity for non-mixed ones ([GL24]). We state a research problem of calculation of $h_1, \ h^1$ for a large class of t-motives.  
\medskip
Section 14 extends the results of Sections 10, 13 to the case of t-motives having $N\ne0$. The main technical tool used for this case is a Hodge-Pink structure ([P], [HJ]). In Remark 14.3 we give a reference to [GL18] where the set of $N$-lattices is interpreted in terms of a generalized flag variety. 
\medskip
In Section 15 we give a definition of a L-function of an Anderson t-motive (only one type of two existing types of L-functions). We show that its two definitions are equivalent, giving in details the proof of [G96]. We give its explicit calculation for a Drinfeld module of rank 2 over the field $\n F_q$ (a local factor), Section 15.5, and for the Carlitz module over $\n F_2$ (the global product), Example 15.6.3. 
\medskip
In Section 16 we consider the theory of affine equations and their systems --- an important technical tool to find $H^1(M), \ H_1(M),\ T_\g L(M)$. We introduce the notion of $T$-divisible submodules of $\p[[T]]$ and their small elements. Finally, we introduce the notion of holonomic sequences in $\p$. We conjecture that their properties are analogous to the properties of classical holonomic sequences. 
\medskip
In Section 17 we start to describe the Drinfeld upper half plane $\oo$. At the moment we give a detailed description of the sets $D(n,x)$ that cover $\oo$, and we show how to use them in order to define the Bruhat-Tits building of $\oo$ (To be continued). 
\medskip
In Section 18 we give an explicit calculational proof of the fact that for any prime ideal $\pi$ of $\n F_q[T]$ we have: an analog of $H_1(M)$ for $\pi$ is isomorphic to the $\pi$-Tate module $T_\pi(M)$, see [G95], p. 186. 
\medskip
In Section 19 we give the function field case analog of the notion of the Siegel matrix for an abelian variety with MIQF. Recall (see the formula (1.8.1) and the above lines) that we have a map from the set of $O_K$-bases of the lattices of abelian varieties with MIQF, to the set of complex $(n\times (r-n))$-matrices satisfying (1.8.1). The function field case analog of this construction is a map $\de$ from the set of bases of lattices (up to an equivalence) to the Grassmannian $Gr(n,r)$ for the case $N=0$ (Sections 19.2 -- 19.3), or to a flag variety (denoted by tFV) for the case $N\ne0$, Sections 19.6.6 (case of the maximal Schubert cell) and 19.6.7 (general case). The flag variety tFV which is the target of the map $\de$ for the case $N\ne0$ is (apparently) of a new type, i.e. it was never considered earlier. 
\medskip
These results can be applied to the duality theory (Section 19.5) and to the definition of analogs of the modular forms, called t-modular forms (Section 19.4).\footnotemark \footnotetext{At the moment we did not check that non-0 t-modular forms really exist; there are some doubts.}
\medskip
Finally, in Section 99 we consider some generalizations and modifications of Anderson t-motives.
\medskip
Many important subjects are not considered in this short survey. We do not consider problems related to parametrizations of Drinfeld modules and their generalizations, hence we do not consider analogs of Eichler-Shimura theorem on reductions of Hecke correspondences on modular curves for the function field case, and all further theory leading to proofs of Langlands conjectures. This will be made in future.
\medskip
%\newpage
{\bf 1. Abelian varieties in characteristic 0. }
\medskip
Here we describe briefly main objects related to abelian varieties over number fields.
\medskip
An abelian variety $A$ over $\n C$ is a projective variety having a group structure. Let dim $A=g$. Analytically, it is $\n C^g/L$, where $L=\n Z^{2g}$ is a lattice in $\n C^g$ satisfying the
Riemann condition:
\medskip
$\exists \ H=B+i\Omega$ --- a non-degenerated hermitian form on $\n C^g$ (here $B$, $\Omega$ are its real and imaginary parts) such that
\medskip
{\bf 1.1.1.} $H$ is a positively defined;
\medskip
{\bf 1.1.2.} $\Omega|_L\in \n Z$, i.e. for $u, \ v \in L\subset \n Z^{2g}$ we have $\Omega(u,v)\in \n Z$.
\medskip
There exists a basis $\{e\}=\{e_1,\dots, e_g, e_{g+1},\dots, e_{2g}\}$ of $L$ over $\n Z$ such that the matrix of $\Omega$ in this basis is $\left(\matrix 0&D\\ -D&0\endmatrix \right)$ (entries are $g\times g$-blocks) where $D=\diag(d_1,d_2,\dots,d_g)$ is a diagonal $g\times g$-matrix with integer positive entries satisfying $d_1\ | \ d_2 \ | \ ... \ |\ d_g$. If all $d_i$ are 1 then $A$ (more exactly, a pair $\{A, \ H\}$ ) is called a principally polarized variety; for simplicity, we shall consider only them.

There exists a matrix $S\in M_{g\times g}(\n C)$ such that
\medskip
$\left(\matrix e_{g+1}\\ \dots \\ e_{2g}\endmatrix \right)= S\left(\matrix e_{1}\\ \dots \\ e_{g}\endmatrix \right)$ (equality in $\n C^g$).
\medskip
$S$ is called a Siegel matrix of $A$ (and of the basis $\{e\}$ ).
\medskip
Conditions 1.1.1, 1.1.2 are equivalent to (case of principally polarized varieties):
\medskip
{\bf 1.2.} $S$ is symmetric, and $Im(S)$ is positively defined.
\medskip
The set of Siegel matrices is denoted by $\g H_g$ (the Siegel upper half plane).
\medskip
{\bf 1.3.} The symplectic group $Sp_{2g}(\n Z)$ acts on $\g H_g$. Action of $Sp_{2g}(\n Z)$ corresponds to a change of basis of $L$ over $\n Z$.
\medskip
Two lattices are called equivalent (notation: $L_1\sim L_2$) if there exists a $\n C$-linear map $\vf: \n C^g \to \n C^g$ such that $\vf(L_1)=L_2$. Equivalent lattices have the same Siegel matrices (in appropriate bases).
\medskip
{\bf Theorem 1.4.} There exists a 1 -- 1 equivalence between the set of abelian varieties (up to isomorphism) and the set of $L$ satisfying the Riemann condition, up to equivalence.
\medskip
{\bf Corollary.} The set of principally polarized abelian varieties (up to isomorphism) is isomorphic to $\g H_g/Sp_{2g}(\n Z)$.
\medskip
Particularly, the dimension of the moduli space of abelian varieties of dimension $g$ is $\binom{g+1}2=\frac{(g+1)g}2$, because $S$ is symmetric.
\medskip
{\bf 1.5. Tate modules.} Let $A_n$ be the group of $n$-torsion points of an abelian variety $A$. We have $A_n=(\n Z/n)^{2g}$. Let $A$ be defined over $\n Q$. In this case Gal$(\bar \n Q/\n Q)$ --- the Galois group of $\n Q$ --- acts on $A_n$.
\medskip
Let $l$ be a prime. The Tate module of $A$ is defined as follows: $T_l(A):=\underset{\underset{n\to\infty}\to{\longleftarrow}}\to{\hbox{lim}} A_{l^n}$. We have $T_l(A)=\n Z_l^{2g}$. Gal($\n Q$) acts on $T_l(A)$.
\medskip
Let $p$ be a prime. We tell that $A$ has a good reduction at $p$ if there exists a system of equations defining $A$ such that all coefficients of these equations belong to $\n Q$, are $p$-integer and after reduction of these coefficients modulo $p$, we get a system of equations defining an abelian variety over $\n F_p$ (this is only a rough definition; for an exact definition we should distinguish forms of $A$ over $\n Q$ and to prove that the reduction does not depend on a choice of a system of equations defining $A$).
\medskip
If $A$ has a good reduction at $p$ then the Frobenius automorphism Fr($p$) acts on $T_l(A)$ (for simplicity, we consider the case $p\ne l$). It is defined up to a conjugation, and its characteristic polynomial $C_{A,p,l}[U]$ is uniquely defined ($U$ is an independent variable).
\medskip
A priori, $C_{A,p,l}[U]\in \n Z_l[U]$ and hence depends on $l$. Really, we have
\medskip
{\bf Theorem 1.5.1.} $C_{A,p,l}[U]\in \n Z[U]$ and does not depend on $l$.
\medskip
{\bf Idea of the proof.} Instead of $\tilde A:=$ the reduction of $A$ at $p$, we can consider as an initial object any abelian variety $B$ of dimension $g$ defined over $\n F_p$. For a prime $l\ne p$ the Tate module $T_l(B)$, the action of Frobenius on $T_l(B)$ and its characteristic polynomial $C_{B,l}[U]$ are defined. 
\medskip
But since $B$ defined over $\n F_p$, the Frobenius map is an algebraic homomorphism: $fr: B\to B$. Really, if $t\in B$ and $t$ has coordinates $(x_0, \dots, x_N)$ in an embedding $B\hookrightarrow P^N(\bar \n F_p)$, then $fr(t)$ has coordinates $(x_0^p, \dots, x_N^p)$. They are polynomials, hence the Frobenius map is an algebraic homomorphism. 
\medskip
The endomorphism ring $\End(B)$ is $\n Z^k$ for some $k$. For almost all $B$ we have $k=2g$. In this case we have $\End(B)=\n Z[fr]$ and the minimal polynomial of $fr\in \End(B)$ is exactly $C_{B,l}[U]$. $\square$
\medskip
%\newpage
{\bf Theorem 1.6. (A. Weil)}: The roots $ \alpha_{p,1}, \dots, \alpha_{p,2g}\in \n C$ of $C(U)$ have properties:
\medskip
1.6.1. $\alpha_i\alpha_{g+i}=p$;
\medskip
1.6.2. $|\alpha_i|=p^{1/2}$.
\medskip
Particularly, for an elliptic curve the valuation $v_p$ splits in $v_{p,1}$ and $v_{p,2}$ in $\n Q(\al_*)/\n Q$, and we have 
\medskip
1.6.3. $v_{p,1}(\al_1)=1, \ v_{p,1}(\al_2)=0, \ v_{p,2}(\al_1)=0, \ v_{p,2}(\al_2)=1$.
\medskip
{\bf 1.7. Reduction.} We need the following definition. Let $A$ be an abelian variety of dimension $g$ over $\bar \n F_p$, and let $A_p$ be the group of its closed points of order $p$. We have $A_p=(\n Z/p)^{g_0}$ where $0\le g_0\le g$. $A$ is called ordinary if $g_0=g$, i.e. $A_p=(\n Z/p)^g$. An equivalent definition: let $\Cal A_p$ be the group scheme of $p$-torsion of $A$, i.e. the kernel of multiplication by $p$. $A$ is ordinary iff $\Cal A_p=(\mu_p)^g\oplus (\n Z/p)^g$ as a group scheme, where $\mu_p$ is the group scheme Spec $\bar \n F_p[x]/(x^p-1)$.
\medskip
Here we consider only two properties concerning reduction of abelian varieties.
We consider an abelian variety $A$ of dimension $g$ defined over $\n Q$ such that it has a good reduction $\tilde A$ at a prime $p$, and this $\tilde A$ is ordinary. We have a reduction map on points of order $p$: $A_p \to \tilde A_p$.
\medskip
{\bf Theorem 1.7.1.} It is surjective, hence its kernel has dimension $g$ over $\n F_p$. This kernel is
an isotropic subspace of $A_p$ with respect to a skew form on it (coming from the above $\Omega$).
\medskip
{\bf Theorem 1.7.2.} For any $k=0,\dots, g$ the dimension of the moduli space of abelian varieties of dimension $g$ over $\bar \n F_p$ whose $g_0$ is $\le k$ is $\binom{k+1}2=\frac{(k+1)k}2$.
\medskip
Particularly, "almost all" abelian varieties over $\bar \n F_p$ are ordinary; there are only finitely many abelian varieties having $g_0=0$.
\medskip
{\bf 1.8. Abelian varieties with MIQF.}
\medskip
Since Anderson t-motives are analogs of abelian varieties with MIQF, we recall briefly their definition and properties. Let $K$ be an imaginary quadratic field, $A$ an abelian variety of dimension $r$ such that there exists an inclusion $K\hookrightarrow \End^0(A):=\End(A)\otimes \n Q$ (we fix this inclusion).
\medskip
Rings of endomorphisms of abelian varieties are described for example in [Sh]. Our case of $A$ with MIQF
is the type IV  of [Sh], Proposition 1, p. 153. For this case (for generic $A$), in notations of [Sh], we have $F=\n Q, \ F_0 = K $,
i.e. the degrees $g$, $q$ of [Sh], 2.1, p. 155, are both 1. Further, $m$ of [Sh], (7), p. 156 is $r$.
\medskip
A description of the action of $K$ on $A$ is the following. Let $V=\n C^r$ and $A=V/L$ as above. The field $K$ acts on $V$. The space $V$ is the sum of two subspaces: $V=V^+\oplus V^-$ where the action of $K$ on $V^+$ is the "direct" one (i.e. simply the multiplication by the corresponding element of $K$), while the action of $K$ on $V^-$ is the "conjugate" one: for $x\in K, \ v\in V^-$ we have $x(v)=\bar x\cdot v$ (here $x(v)$ is the action of $x$ on $v$, bar means the complex conjugation, and $\bar x\cdot v$ is the multiplication in the vector space $V$).
\medskip
Let $n:=\dim V^+$, hence $\dim V^-=r-n$.
The pair $(n, r-n)$ is called the signature of $A$, it is $(r_1, s_1)$ of [Sh], (8), p. 156.
\medskip
We do not give a definition of a Siegel matrix for the case of these $A$ (see [Sh] or, more generally, [De] for a general definition). We only indicate a general principle. Let $\End(A)=O_K$ (the simplest possible case; it can happen that $\End(A)$ is strictly less than $O_K$), and let the above lattice $L$ be a free $O_K$-module with respect to the action of $K$ on $V$. Let $l_*=\{l_1, \dots, l_r\}$ be a $O_K$-basis of $L$. We can associate to $l_*$ a $n \times (r-n)$-complex matrix $S$ satisfying (see [Sh], 2.6, p. 162)
$$ \hbox{$I_n-S\bar S^t$ is positively defined} \eqno{(1.8.1)}$$
It is called a Siegel matrix of $A$. The set of these matrices is denoted in [Sh] by $\Cal H^3_{n,r-n}$. We see that it belongs to the affine space of $n \times (r-n)$-matrices over $\n C$ (to compare with (19.2.1): the analog of the Siegel matrix for the function field case belongs to the Grassmannian $Gr(r-n,r)$; the above affine space is its maximal Schubert cell).
\medskip
{\bf 1.8.2.} Particularly, the dimension of the moduli space of abelian varieties with MIQF of signature $(n,r-n)$ is $n(r-n)$, because the condition (1.8.1) does not impose algebraic relations on entries of $S$.
\medskip
We can associate a reductive group over $\n Q$ to any type of abelian varieties with a fixed endomorphism ring (see [De] for a much more general situation). For abelian varieties with MIQF of signature $(n, r-n)$ this group is $GU(n,r-n)$.
\medskip
We do not give here a definition of ordinariness of reduction of abelian varieties with MIQF over $\overline{\n F_q}$. We indicate only that if $A$ is ordinary as a variety with MIQF then it is not ordinary as a variety obtained from $A$ by forgetting the MIQF-structure (unless $n=r-n$). Further, like for the general abelian varieties, "almost all" $A$ with MIQF are ordinary as a variety with MIQF. In this case the kernel of the reduction map of
Theorem 1.7.1 is $n$ (if $n\ge r-n$).
\medskip
{\bf 1.8.3.} Finally, we indicate an amusing construction of a lattice for an abelian variety with MIQF. This lattice is an analog of a lattice of an Anderson t-motive. Namely, we shall see (Section 10) that if $M$ is an uniformizable Anderson t-motive of rank $r$ and dimension $n$ then we can associate it a lattice in $\p^n$ of dimension $r$ over $\n F_q[\th]$ --- the functional analog of $\n Z$. Therefore, we can expect that if $A$ is an abelian variety with MIQF of signature $(n, r-n)$ then we can associate it a "lattice" in $\n C^n$. This is really so! This "lattice" is an $O_K$-module of rank $r$ in $\n C^n$ having some properties. See [GL09], Theorem 2.6 for the exact statement; there is a 1 -- 1 correspondence between abelian varieties with MIQF and such "lattices".
\medskip
This is a rare example of "an analogy to the opposite direction": we consider a construction in the theory of Anderson t-motives, and we find its analog in the theory of abelian varieties with MIQF.
\medskip
{\bf 2. Abelian varieties in characteristic $p$.}
\medskip
Let $q$ be a power of a prime $p$, $\n F_q$ a finite field of order $q$.
Let $\th$ be an abstract transcendent element. The analog of
$\n Z \subset \n Q$ is $\n F_q[\theta]   \subset \n F_q(\theta)$. Let $\overline{\n F_q(\theta)}$ be an algebraic closure of $\n F_q(\theta)$.

Let $A$ be an abelian variety over $\overline{\n F_q(\theta)}$ of dimension $g$. There is no analog of the above formula $A=\n C^g/L$, but we have
\medskip
$A_n=(\n Z/n)^{2g}$ if $(n,p)=1$ (as earlier),
\medskip
and hence $T_l(A)=\n Z_l^{2g}$ ($l\ne p$),
\medskip
If $A$ is defined over $\n F_q(\theta)$, then Gal($\n F_q(\theta)$) acts on $T_l(A)$, and we have an analog of the Weil Theorem.
\medskip
\medskip
Hence, we have the following table for abelian varieties:
\settabs 6 \columns
\medskip
%\newpage
\+ &&& {\bf Table 2.1} \cr
\medskip
\+ &&& Tate module && Galois group\cr
\medskip
\medskip
\+Abelian varieties over $\n Q$ &&& $\n Z_l^{2g}$ && Gal($\n Q$)\cr
\medskip
\+Abelian varieties over $\n F_q(\theta)$ &&& $\n Z_l^{2g}$ && Gal($\n F_q(\theta)$)\cr
\medskip
\medskip
{\bf 3. Initial rings and fields}
\medskip
For the characteristic 0 we have inclusions of rings and fields:
$\n Z \subset \n Q  \subset \n R \subset \n C$.
\medskip
Let us consider their analogs in characteristic $p$. As it was mentioned above, the analog of
$\n Z \subset \n Q$ is $\n F_q[\theta]   \subset \n F_q(\theta)$.
\medskip
Valuations of $\n Q$ are $v_p$ where $p$ is prime, and $v_\infty$ the Archimedean valuation.
\medskip
Valuations of $\n F_q(\theta)$ have a similar description. Let $P\in \n F_q[\theta] $ be an irreducible polynomial. It defines a valuation $v_P$ on $\n F_q(\theta)$. The analog of $v_\infty$ on $\n F_q(\theta)$ (it has the same notation $v_\infty$) is the valuation "minus degree": for $S\in \n F_q(\theta) $ we have $v_\infty(S):=-$degree $(S)$. Equivalently, this is the order of zero of a function at infinity; also $v_\infty$ can be defined as the only valuation satisfying $v_\infty(\th)=-1$. Unlike $v_\infty$ for the number field case, $v_\infty$ for the function field case is a non-archimedean valuation.

Valuation $v_\infty$ defines a topology in $\n F_q(\theta)$. Later we shall consider only this valuation and its topology.
We have in it: $\theta^{-n} \to 0$ for $n\to +\infty$.
\medskip
We have $\n R_\infty=\n F_q((\theta^{-1}))$ is the completion of $\n F_q(\theta)$. Like in the number field case, $\n F_q[\theta]$ is discrete in $\n R_\infty$, and the quotient $\n R_\infty/\n F_q[\theta]=\th^{-1}\n F_q[[\theta^{-1}]]$ is compact.

Finally, $\p:=\widehat {\overline {\n R_\infty}}$ is the completion of the algebraic closure of $\n R_\infty$. It is
complete by definition, and algebraically closed (see [G96], Proposition 2.1). $\n R_\infty \subset \p$ is a characteristic $p$ analog
of $\n R \subset \n C$.
\medskip
We denote by $\tau$ the Frobenius automorphism of $\p$: $\tau(x):=x^q$.
\medskip
{\bf Remark 3.1.} Let us consider a field formed by convergent series generated
by rational powers of $\th^{-1}$, with coefficients in $\bar \n F_p$. More exactly, let
$\al_1 < \al_2 < \al_3 < ... $ be a sequence of rational numbers tending to $+\infty$, and let
$c_i \in \bar \n F_p$ be coefficients. The series $\sum_{i=0}^\infty c_i \th^{-\al_i}$ form a field denoted by ${\p}_s$. We have $\n R_\infty\subset {\p}_s\subset \p$. The reader can think that ${\p}_s=\p$, but this is wrong. A well-known example of $\g r\in \p-{\p}_{s}$ (see [GL21], Remark 4.3) is a root to the equation $$x^2+x+\th^2=0\eqno{(3.2)}$$ (here $q=2$). Really, formally we have $$\g r=\th+\th^{\frac12}+\th^{\frac14}+\th^{\frac18}+...\eqno{(3.3)}$$ but this series $\not\in {\p}_{s}$, because $-\frac1{2^n}$ does not tend to $+ \infty$ as $n\to\infty$. The $n$-th approximation $\g r_{n,i}$ to $\g r$ (here $i=1,2$: there are two roots to (3.2), it is separable) is given by the formula $$\g r_{n,i}=\th+\th^{\frac12}+\th^{\frac14}+\th^{\frac18}+...+\th^{\frac1{2^n}}+\delta_{in}$$ We have: $\delta_{in}$ is a root to $$y^2+y+\th^{\frac1{2^n}}=0$$ and hence both $\delta_{1n}$, $\delta_{2n}$ have $v_\infty(\delta_{1n})=v_\infty(\delta_{2n})=-\frac1{2^{n+1}}$. This shows once again that the series (3.3) does not converge to $\g r$.
\medskip
This phenomenon plays an important role in explicit calculations, see for example [GL21], proof of 4.1, 4.2. It shows that not all equations can be solved by a method of consecutive approximations.
\medskip
{\bf 4. Lattices in functional fields}
\medskip
The dimension of $\n C$ over $\n R$ is 2, hence a $2g$-dimensional lattice $L$ over $\n Z$ in $\n C^g$ is complete: the quotient $\n C^g/L$ is compact, and a basis of $L$ over $\n Z$ is also a basis of $\n C^g$ over $\n R$.

Unlike the number field case, the dimension of $\p$ over $\n R_\infty$ is infinite (and moreover of cardinality continuum), hence all lattices in $\p^n$ are "incomplete".
\medskip
{\bf Definition 4.1.} Let $L=\n F_q[\theta]^r$ and
$L \subset \p^n$. $L$ is called a lattice if:
\medskip
1. $L$ generates all $\p^n$ over $\p$;
\medskip
2. $L$ generates a space of dimension $r$ over $\n R_\infty$ (i.e. elements of a basis of $L$ over $\n F_q[\theta]$ are linearly independent over $\n R_\infty$).
\medskip
We see that the pair $(L, \p^n)$ has two discrete parameters: $r$ and $n$. Their analogs in the number field case are $2g$, resp. $g$.
\medskip
Later we shall need a generalization of the notion of lattice, it is called a $N$-lattice where $N$ is a nilpotent operator on $\p^n$. See 10.5.1 and the below formulas. 
\medskip
{\bf Remark.} It is meaningful to consider incomplete lattices in the number field case as well. For example, we have the exponential map
$$0\to L \to \n C \overset{exp}\to{\to}\n C^*\to0$$
where $L$ is a 1-dimensional lattice in $\n C$.
\medskip
{\bf Definition 4.1a.} Two lattices $L_1\subset \p^n$, $L_2\subset \p^n$ are called (strongly) isomorphic if there exists an isomorphism $\vf:\p^n\to\p^n$ (i.e. $\vf\in GL_n(\p)$) such that $\vf(L_1)=L_2$. 
\medskip
{\bf Remark 4.1b.} The group $GL_n(\p\{\tau\})$ (see 5.1 below for the definition of $\p\{\tau\}$) acts on $\p^n$ considered as a $\n F_q$-module. Hence, if $L$ is a lattice and $\vf \in GL_n(\p\{\tau\})$, then $\vf(L)$ is not a lattice. 
\medskip
{\bf Question 4.1c:} Is it possible (and meaningful) to define a notion of weak equivalence of lattices $L_1$, $L_2$ based on a fact that "$L_1$, $L_2$ are related via $\vf \in GL_n(\p\{\tau\})$"? Apparently not. 
\medskip
Let $(a_1,\dots,a_n)\in \p^n$, where $a_i\in \p$. We define $v_\infty(a_1,\dots,a_n):=\min(v_\infty(a_i))$. It satisfies to all properties of valuation. 
\medskip
For any $k\in \n Z$ we denote $U_k:=\{x\in \p^n\ |\  v_\infty(x)\ge k\}$. 
\medskip
{\bf Theorem 4.2.} For any $k\in \n Z$ (here $k$ can be greatly negative) we have: $L\cap U_k$ is finite. 
\medskip
 {\bf Proof.} Let $x_1, \dots, x_r$ be a basis of $L$ over $\n F_q[\th]$. 
\medskip
{\bf Lemma 4.3.} $\exists \ c$ (it can be big) such that $\forall \ k $, $\forall \ \al_1,\dots, \al_r\in \n R_\infty$  we have: 
$$v_\infty(\sum_{i=1}^r \al_i x_i)\ge k \ \implies \forall \ i \ \ \  v_\infty(\al_i)\ge k-c.$$

{\bf Proof} of 4.3. We shall prove it by induction by $r$. For $r=1$ it is obvious. Let it be true for $m$, i.e. 
\medskip
$\exists \ c=c_m$ such that $\forall \ k $, $\forall \ \al_1,\dots, \al_m\in \n R_\infty$  we have: 

$$v_\infty(\sum_{i=1}^m \al_i x_i)\ge k \ \implies v_\infty(\al_i)\ge k-c_m, \ \forall \ i \eqno{(4.4)}$$

and we prove that the same is true for $m+1$. We need a 
\medskip
{\bf Sublemma 4.5.} $\exists \ d=d_m$ such that $\forall \ \al_1,\dots, \al_m\in \n R_\infty$  we have: 

$$v_\infty(x_{m+1}-\sum_{i=1}^m \al_i x_i)\le d\eqno{(4.6)}$$
(i.e. $x_{m+1}$ is not very close to the $\n R_\infty$-linear envelope of $x_1, \dots, x_m$). 
\medskip
{\bf Proof} of 4.5. Let, conversely, $\forall j\in \n Z$ $\exists \ \al_{j1},\dots, \al_{jm}\in \n R_\infty$  such that 
$$v_\infty(x_{m+1}-\sum_{i=1}^m \al_{ji} x_i)\ge j\eqno{(4.7)}$$
Hence (here  $j_1< j_2$) we have 
$$v_\infty(\sum_{i=1}^m (\al_{j_1i}-\al_{j_2i}) x_i)\ge j_1$$
and (4.4) implies $v(\al_{j_1i}-\al_{j_2i})\ge j_1-c_m$. Hence, $\forall \ i$ the sequence $\al_{ji}$ is a fundamental sequence (as $j\to \infty$). Since $\al_{ji}\in \n R_\infty$, there exists $\al_i:=\lim(\al_{ji})$. Because of (4.7), we have that $\forall \ j$ 
$$v_\infty(x_{m+1}-\sum_{i=1}^m \al_i x_i)\ge j$$
hence $x_{m+1}=\sum_{i=1}^m \al_i x_i$ --- a contradiction to the linear independence of $x_*$. 4.5 is proved. 
\medskip
Let us continue the proof of 4.3. Let 
$$v_\infty((\sum_{i=1}^m \al_i x_i)+\al_{m+1}x_{m+1})\ge k $$
If $\al_{m+1}=0$ then all is good. If not then 
$$v_\infty((\sum_{i=1}^m \frac{\al_i }{\al_{m+1}}x_i)+x_{m+1})\ge k - v(\al_{m+1})$$
According (4.6), 
$$d_m\ge v_\infty((\sum_{i=1}^m \frac{\al_i }{\al_{m+1}}x_i)+x_{m+1})\ge k - v_\infty(\al_{m+1})$$
i.e. $$v_\infty(\al_{m+1})\ge k-d_m.$$
Since all coefficients are equivalent, 4.3 is proved. 
\medskip
The theorem follows immediately, because $\forall \ k$ there are only finitely many 
\medskip
$a\in (\n F_q[\th])^r$ such that $v_\infty(a)\ge k-c$. $\square$
\medskip
{\bf 4.8.} The definition of a Siegel matrix for the function field case is also the same as the one for the number field case. Let $e_*:=\{e_1,\dots, e_n, e_{n+1},\dots, e_{r}\}^t$ (here and below $t$ means transposition) be a basis of $L$ over $\n F_q[\theta]$ such that
$e_1,\dots, e_n$ is a basis of $\p^n$ over $\p$. Hence, there exists a matrix $S=S_{ij}\in M_{(r-n)\times n}(\p)$ such that
\medskip
$\left(\matrix e_{n+1}\\ \dots \\ e_{r}\endmatrix \right)= S\left(\matrix e_{1}\\ \dots \\ e_{n}\endmatrix \right)$, i.e. $\forall \ j=1,\dots, r-n$ we have $e_{n+j}=\sum_{i=1}^n S_{ji}e_i$.
\medskip
It is called a Siegel matrix of $L$ in the basis $\{e_*\}$.
\medskip
Siegel matrices are used to check whether lattices are isomorphic,or not. Namely, lattices $L\subset \p^n$, $L'\subset \p^n$ are isomorphic iff there exist bases $e_*$, $e'_*$ of $L$, $L'$ respectively such that the Siegel matrix of $L$ in $e_*$ is equal to the Siegel matrix of $L'$ in $e'_*$.
\medskip
Let us consider a function field case analog of (1.3). Let $e'_*:=\{e'_1,\dots,  e'_{r}\}^t$ be another basis of $L$ over $\n F_q[\theta]$ and $\g g\in GL_r(\n F_q[\theta])$ be the matrix of the change of basis from $e_*$ to $e'_*$ (we use the agreement $e'_*= \g g \cdot e_* $, where $e_*$, $e'_*$ are considered as column matrices of size $r \times 1$).
\medskip
Unlike the number field case, it can happen that $e'_1,\dots,  e'_{n}$ is not a $\p$-basis of $\p^n$. See Section 19 where it is written that we should consider the Siegel matrix as an element of a Grassmannian $Gr(r-n,r)$. 

%In this case, a Siegel matrix for $e'_*$ does not exist. If $e'_1,\dots,  e'_{n}$ is a $\p$-basis of $\p^n$ then $S'$ --- the Siegel matrix of $L$ in the basis $\{e_*\}$ --- is related with $S$ by the same formula as in the number field case. Namely, we consider a block form of $\g g$: $\g g=\left(\matrix  \ga_{11}& \ga_{12}\\ \ga_{21} & \ga_{22} \endmatrix\right)$ where blocks $\ga_{11},\ \ga_{12},\ \ga_{21},\ \ga_{22}$ are of sizes respectively $n\times n$, $n\times (r-n)$, $(r-n)\times n$, $(r-n)\times (r-n)$. We have: $$S'=(\ga_{21}+\ga_{22}S)(\ga_{11}+\ga_{12}S)^{-1}\eqno{(4.9)}$$
\medskip
Condition $\ga_{11}+\ga_{12}S\in GL_{n}(\p)$ is equivalent to the condition that $e'_1,\dots,  e'_{n}$ form a $\p$-basis of $\p^n$.
\medskip
We see that the difference with the number field case is the following. First, there is no notion of polarization, i.e. there is no conditions like symmetry imposed on $S$ --- instead of the symplectic group $GSp$ we have the linear group $GL$. Second, we should consider Grassmannians instead of affine spaces. %formula 4.9 defines not the action of $GL_r(\n F_q[\theta])$ on the set of Siegel matrices, but only an "almost action": if $\ga_{11}+\ga_{12}S\not\in GL_{n}(\p)$ then the action is not defined.
\medskip
We see that really we have an analogy with the case of abelian varieties with MIQF. We have the same size of Siegel matrices and the same reductive group acting on them.
\medskip
\medskip
{\bf 5. Anderson t-motives}
\medskip
Let $\p[T,\tau]$ be a ring of non-commutative polynomials in two variables $T$, $\tau$ with
the following relations (here $a\in \p$):
$$aT=Ta; \ \ \tau T=T\tau; \ \ \tau a=a^q\tau\hbox{ (and hence $\tau^k a=a^{q^k}\tau^k$ )}\eqno{(5.1)}$$
It is called the Anderson ring. $\p[T,\tau]$ has subrings $\p[T]$, $\p\{\tau\}$.
\medskip
{\bf Definition 5.2.} Anderson t-motive $M$ is a left $\p[T,\tau]$-module satisfying conditions:
\medskip
5.2.1. $M$ as a $\p[T]$-module is free of finite dimension (denoted by  $r$);
\medskip
5.2.2. $M$ as a $\p\{\tau\}$-module is free of finite dimension (denoted by  $n$);
\medskip
5.2.3. The action of $T-\th$ on $M/\tau M$ is nilpotent.
\medskip
Homomorphisms of Anderson t-motives are module homomorphisms.

Numbers $r$, resp. $n$ are called the rank (resp. dimension) of $M$. We shall use by default these notations $r$, $n$. The number $\frac{n}{r}$ is called the weight of $M$, notation $w(M)$. It will be used in 5.17 and in 15.7. 
\medskip
For generalizations of the notion of Anderson t-motive see Section 16. In most generalizations $M$ is considered as a $\p[T]$-module with a skew $\tau$-action, i.e. for $m\in M$ we consider $\tau(m)$ instead of $\tau m$; skew action means that $\tau(am)=a^q\tau(m)$ (here $a\in\p$), according (5.1).
\medskip
{\bf Definition 5.3.} Drinfeld module\footnotemark \footnotetext{See Remark 6.1 why similar objects have different names --- modules and motives.} is an Anderson t-motive of dimension $n=1$.
\medskip
{\bf Example 5.4.} Let $M$ be a Drinfeld module, and $\{e\}=e_1$ the only element of a basis of $M$ over $\p\{\tau\}$.
This means that any $m\in M$ can be uniquely written
\medskip
$$m=(c_0+c_1\tau+c_2\tau^2+\dots +c_{k-1}\tau^{k-1}+c_k\tau^k)e\eqno{(5.5)}$$
\medskip
where $c_i\in \p$.
\medskip
To define a left $\p[T,\tau]$-module structure on $M$, it is sufficient to define the element $Te\in M$.
(5.5) implies that there exist $a_0, \dots, a_r\in\p$, $a_r\ne 0$, such that

$$Te=(a_0+a_1\tau+a_2\tau^2+\dots +a_{r-1}\tau^{r-1}+a_r\tau^r)e\eqno{(5.6)}$$

Hence, elements $a_0, \dots, a_r$ define $M$ uniquely.
Condition (5.2.3) implies that $$a_0=\theta\eqno{(5.7)}$$

{\bf Exercise 5.8.} $r$ is the rank of $M$, elements

$$e, \tau e,  \tau^2 e,\dots,  \tau^{r-1} e$$

form a basis of $M$ over $\p[T]$.
\medskip

We have an analog of (5.6) for Anderson t-motives. Let $\{e\}=\{e_1, \dots, e_n\}^t$ be a basis of $M$ over $\p\{\tau\}$ considered as a matrix column. Instead of $a_i$ of (5.6), we have $n\times n$ matrices $A_i$ with entries in $\p$. (5.6) becomes a matrix equality
$$T\{e\}=(A_0+A_1\tau+A_2\tau^2+\dots +A_{k-1}\tau^{k-1}+A_k\tau^k)\{e\}$$
$$=A_0\{e\}+A_1\tau\{e\}+A_2\tau^2\{e\}+\dots +A_{k-1}\tau^{k-1}\{e\}+A_k\tau^k\{e\}\eqno{(5.9)}$$
Condition (5.2.3) implies that $N:=A_0-\th I_n$ is a nilpotent matrix.
Anderson t-motives having $N=0 \ \ ( \iff A_0=\th I_n)$ are more simple objects than
Anderson t-motives having $N\ne 0$. Clearly $A_i$ and hence $N$ depend on a choice of basis $\{e\}$. 
\medskip
An analog of Exercise 5.8 for this case is
\medskip
{\bf Exercise 5.10.} If $\det A_k\ne0$ then elements $\tau^i e_j$, $i=0,\dots,k-1$, $j=1,\dots,n$ form a basis of $M$ over $\p[T]$.
\medskip
Hence if $\det A_k\ne0$ then $r=kn$. We see that for interesting examples of Anderson t-motives we should have
$\det A_k=0$.
\medskip
{\bf 5.11.} We can consider an Anderson t-motive $M$ not as a $\p\{\tau\}$-module with $T$-action, but as a $\p[T]$-module with $\tau$-action. Moreover, this type of consideration is used more frequently in applications. Let $\{f\}=\{f_1, \dots,
f_r\}^t$ be a basis of $M$ over $\p[T]$ considered as a matrix column. There exists a matrix
$Q\in M_{r\times r}\p[T]$ defining the multiplication by $\tau$, namely
$$\tau \{f\}=Q\{f\}\eqno{(5.11.1)}$$

{\bf Lemma 5.11.2.} $\det Q=c(T-\th)^n$ where $c\in \p^*$.
\medskip
{\bf Proof.} Condition (5.2.3) implies that a power of $T-\th$
annihilates $M/\tau M$. Hence, $M/\tau M$ is isomorphic as a $\p[T]$-module to the direct sum of finite cyclic $(T-\th)$-power torsion $\p[T]$-modules. Its dimension over $\p$ is $n$. On the other hand, (5.11.1) implies that as $\p[T]$-modules, we have
$$M/\tau M\cong M_{1\times r}(\p[T]) \ / \ M_{1\times r}(\p[T])\cdot Q.$$
This implies the lemma. $\square$
\medskip

{\bf Exercise 5.12.} For a Drinfeld module $M$ defined by (5.6) such that $a_r=1$, the matrix $Q$ is the following:
$$Q=\left(\matrix 0&1&0&0&\dots&0&0
\\ 0&0&1&0&\dots&0&0
\\ \dots&\dots&\dots&\dots&\dots&\dots&\dots
\\ 0&0&0&0&\dots&0&1
\\ T-\th & - a_1&-a_2&-a_3&\dots&-a_{r-2}&-a_{r-1}\endmatrix\right)\eqno{(5.12.1)}$$
For an Anderson t-motive given by (5.9) such that $A_k=I_n$, the matrix $Q$ is analogous, all entries are
$n\times n$-blocks.
\medskip
{\bf Example 5.13.} Let us consider a particular case of (5.9), namely an Anderson t-motive (denoted by $M(A)$ ) given by the formula (here $A\in M_{n\times n}(\p)$ )
$$T\{e\}=\th\{e\}+A\tau\{e\}+\tau^2\{e\}\eqno{(5.13.1)}$$
i.e. we have $A_0=\th I_n$, $N=0$, $A_1=A$, $A_2=I_n$. It has $r=2n$, it is a close analog of abelian varieties. These Anderson t-motives $M(A)$ are studied in [GL17], [GL21], [EGL].
\medskip
{\bf 5.14. Change of basis.} First, we consider the change of basis for a Drinfeld module $M$. Let $\{e\}=e_1$ be from 5.4 --- the only element of a basis of $M$ over $\p\{\tau\}$. Let $c\in \p^*$ and $e'=c^{-1}e$ be the only element of another basis of $M$ over $\p\{\tau\}$. (5.6) becomes (because $\tau^kc=c^{q^k}\tau^k$)
$$cTe'=c\th e'+c^q a_1\tau e'+c^{q^2}a_2\tau^2 e' +\dots +c^{q^{r-1}}a_{r-1}\tau^{r-1}e'+
c^{q^r}a_r\tau^re'\eqno{(5.14.1)},$$
or
$$Te'=\th e'+c^{q-1} a_1\tau e'+c^{q^2-1}a_2\tau^2 e' +\dots +c^{q^{r-1}-1}a_{r-1}\tau^{r-1}e'+
c^{q^r-1}a_r\tau^re'\eqno{(5.14.2)}$$

Hence, coefficients $\{\th, \ c^{q-1} a_1, \ c^{q^2-1}a_2, \dots, c^{q^{r-1}-1}a_{r-1}, \
c^{q^r-1}a_r\}$ define the same Drinfeld module as coefficients $\{\th, \ a_1, \dots, a_r\}$ (both are coefficients of (5.6)). Particularly, any Drinfeld module can be defined by an equation (5.6) with $a_r=1$.
\medskip
{\bf Corollary 5.14.3.} Two Drinfeld modules of rank two defined by the equations
$$Te=\th e+ a_1\tau e+\tau^2 e$$
$$Te=\th e+ a_2\tau e+\tau^2 e$$
are isomorphic iff $a_2=\be a_1$ where $\be^{q+1}=1$. Really, we can choose the above $c$ as $c=\be^{1/(q-1)}$; we have $c^{q^2-1}=1$.
\medskip
{\bf Corollary 5.14.4.} The moduli space of Drinfeld modules of rank $r$ has dimension $r-1$.
\medskip
Really, a Drinfeld module of rank $r$ is defined by $r$ parameters $a_1,\dots, a_r$, and their equivalence is defined by one above parameter $c$.
\medskip
For Anderson t-motives the notion of isomorphism is more complicated, because a matrix of change of basis from $e_*$ to $e'_*$ belongs to $GL_n(\p\{\tau\})$. For $n=1$ it is simply $\p^*$, but for $n>1$ \ \ $GL_n(\p\{\tau\})$ is a "doubly non-abelian" group (the first "non-abelianity" comes from $GL_n$, the second one from non-abelianity of $\p\{\tau\}$).
\medskip
If a matrix $C$ of change of basis has constant coefficients, i.e. $C\in GL_n(\p)\subset GL_n(\p\{\tau\})$ then (5.9) becomes (a calculation similar to (5.14.2), case $A_0=\th I_n$; here $C^{(m)}$ is a matrix obtained by elevation of all entries of $C$ to the $q^m$-th degree):

$$T\{e'\}=\th I_n\{e'\}+C^{-1}A_1C^{(1)}\tau\{e'\}+C^{-1}A_2C^{(2)}\tau^2\{e'\}+\dots $$ $$ +C^{-1}A_{k-1}C^{(k-1)}\tau^{k-1}\{e'\}+C^{-1}A_kC^{(k)}\tau^k\{e'\}\eqno{(5.14.5)}$$

According a theorem of Lang that $H^1(\n F_{q^k}, GL_n(\p))=1$ (the frobenius acts not only on $\overline{\n F_q}$ but also on $\p$ by the formula $fr(x)=x^q$, hence this cohomology group is well-defined) we get that if $\det A_k\ne0$ then $\exists \ C\in GL_n(\p)$ such that $C^{-1}A_kC^{(k)}=I_n$, i.e. in this case any Anderson t-motive of this type can be defined by an equation (5.9) with $A_k=I_n$.
\medskip
It can happen that $C\in GL_n(\p\{\tau\})-GL_n(\p)$. This phenomenon occurs even for $M$ defined by (5.13.1). More exactly, there exist two matrices $A$, resp. $B$ such that the Anderson t-motives $M(A)$, resp. $M(B)$ defined by (5.13.1) are isomorphic, and the matrix $C$ of a change of basis cannot be chosen in $GL_n(\p)$, but only in $GL_n(\p\{\tau\})-GL_n(\p)$. See [GL17] for examples.
\medskip
{\bf 5.14.6.} Moreover, according the knowledge of the authors, we have no algorithm to check whether two Anderson t-motives $M_1$, $M_2$ defined by (5.9), case $A_k=A'_k=I_n$, are isomorphic, or not. For the case of $C\in GL_n(\p)$ this algorithm clearly exists, because in this case we have $C=C^{(k)}$, i.e. $C\in GL_n(\n F_{q^k})$; there exists only finitely many such $C$. But really $C\in GL_n(\p\{\tau\})$, so we have an infinite problem. A natural way to solve such problems is finding an invariant of $M$. For mixed t-motives (see 5.17.3) the lattice (see below) is an invariant, because the lattice functor is injective for mixed t-motives ([HJ]), but for non-mixed ones it is not injective, see Section 13.
\medskip
{\bf 5.15. Pure t-motives.} There exists an important class of t-motives called pure. See [A86], p. 467, middle of the page, or [G96], 5.5.2 for the definition. We give here its matrix form (see [GL07], Lemma 10.3.1). We define: 
\medskip
$S=T^{-1}$; 
\medskip
For any $A\in M_{r\times r}(\p[T])$ let

$$A^{[k]}:=A^{(k-1)}\cdot A^{(k-2)}\cdot \dots\cdot A^{(1)}\cdot A\eqno{(5.15.1)}$$

{\bf Matrix form of a definition of pure.} Let $M$ be a t-motive and $Q$ from (5.11.1). $M$ is pure iff there exists $C\in
GL_r(\p((S)) \ )$ such that
for some $\al$, $\be>0$
$$S^\al C^{(\be)}Q^{[\be]}C^{-1}\in GL_r(\p[[S]])\eqno{(5.15.2)}$$
i.e. iff $S^\al q C^{(\be)}Q^{[\be]}C^{-1}$ is $S$-integer and its inicial coefficient is invertible.
\medskip
It is possible to show that in this case $\exists \ k$ such that $\al=kn, \ \be=kr$ (see [A86], Lemma 1.10.1, and [G96], Proposition 5.5.6). This follows from the equality $\det Q=c(T-\th)^n$ (Lemma 5.11.2).
\medskip
{\bf 5.15.3. Examples.} 
\medskip
(a) Drinfeld modules are pure. Really, we can take $C=1_r$. For $k\le r$ we have 
$$Q^{[k]}=\left(\matrix 0_{(r-k)\times k}&I_{r-k}
\\ B_{21}(k)&B_{22}(k)\endmatrix\right)$$
(block structure) where $B_{22}(k)\in M_{k\times (r-k)}(\p)$ and $$B_{21}(k)=\left(\matrix S^{-1}*_1&*_2&*_2&\dots &*_2\\  \\
*_2S^{-1}*_1&S^{-1}*_1&*_2&\dots &*_2\\  \\
*_2S^{-1}*_1&*_2S^{-1}*_1&S^{-1}*_1&\dots &*_2\\ &\dots\\ *_2S^{-1}*_1&*_2S^{-1}*_1&*_2S^{-1}*_1&\dots&S^{-1}*_1\endmatrix\right)\in M_{k\times k}(S^{-1}\p[S])$$ where $*_1\in 1+S\p[S]$, $*_2\in\p$. This formula is proved immediately by induction, starting from (5.12.1) for $Q$. We get $$SQ^{[r]}=SB_{21}(r)=\left(\matrix *_1&*_2S&*_2S&\dots &*_2S\\  \\
*_2*_1&*_1&*_2S&\dots &*_2S\\  \\
*_2*_1&*_2*_1&*_1&\dots &*_2S\\ &\dots\\ *_2*_1&*_2*_1&*_2*_1&\dots&*_1\endmatrix\right)$$ Obviously $SQ^{[r]}\in GL_r(\p[[S]])$. 
\medskip
(b) All Anderson t-motives defined by (5.9) such that $\det A_k\ne0$ are pure. The proof is the same; instead of elements of the above matrices we consider $n\times n$ blocks.  
\medskip
(c) Let $M$ be given by an equation 
$$T\left(\matrix e_1 \\ e_2\endmatrix \right)= \left(\matrix \th &0 \\ 0&\th \endmatrix \right) \left(\matrix e_1 \\ e_2\endmatrix \right) +  \left(\matrix a_{11} & 0 \\ a_{21} & 1 \endmatrix \right) \tau \left(\matrix e_1 \\ e_2 \endmatrix \right)
 + \left(\matrix 1&0\\ 0&0 \endmatrix \right) \tau^2 \left(\matrix e_1 \\ e_2 \endmatrix \right)$$
It has $n=2, \ r=3$. It is not pure. Proof: the eigenvalues of its characteristic polynomial do not satisfy Theorem 15.7.1.
\medskip
(d) Let $M$ be given by an equation 
$$T\left(\matrix e_1 \\ e_2\endmatrix \right)= \left(\matrix \th &0 \\ 0&\th \endmatrix \right) \left(\matrix e_1 \\ e_2\endmatrix \right) +  \left(\matrix 0&a_{12}  \\ 1&a_{22}  \endmatrix \right) \tau \left(\matrix e_1 \\ e_2 \endmatrix \right)
 + \left(\matrix 0&1\\ 0&0 \endmatrix \right) \tau^2 \left(\matrix e_1 \\ e_2 \endmatrix \right)$$
It has $n=2, \ r=3$. It is pure. It is dual to a Drinfeld module of rank 3 (see Section 13 for duality). 
\medskip
{\bf Remark 5.15.4.} $M$ of Example (c) is reducible: it has a Drinfeld submodule of rank 2 generated by $e_1$. This explains why it is not pure. Nevertheless, its analog of rank 5 given by the equation 
$$T\left(\matrix e_1 \\ e_2\endmatrix \right)= \left(\matrix \th &0 \\ 0&\th \endmatrix \right) \left(\matrix e_1 \\ e_2\endmatrix \right) + A_1\tau \left(\matrix e_1 \\ e_2\endmatrix \right) + \left(\matrix a_{211} & 0 \\ a_{221} & 1 \endmatrix \right) \tau^2 \left(\matrix e_1 \\ e_2 \endmatrix \right)
 + \left(\matrix 1&0\\ 0&0 \endmatrix \right) \tau^3 \left(\matrix e_1 \\ e_2 \endmatrix \right)$$
 is irreducible and not pure. 
\medskip
More explicit examples of pure and not pure t-motives having explicit equations can be found in [GL07], Section 11. 
\medskip
The terminology "pure" comes from the theorem 15.7.1, see below. 
\medskip
Further, we have
\medskip
{\bf Theorem 5.15.5.} ([H],  Theorem 3.2). The dimension of the moduli space of pure t-motives of
dimension $n$ and rank $r$ is $n(r-n)$.
\medskip
We see that this dimension coincides with the dimension of the moduli space of abelian varities with MIQF, of dimension $r$ and signature $(n,r-n)$, and with the set of Siegel matrices of lattices of dimension $r$ in $\p^n$, see Section 4.
\medskip
{\bf Example 5.16.} Case $r=1$.
The only Drinfeld module of rank 1 satisfies
$$Te=\th e +\tau e\eqno{(5.16.1)}$$
It is called the Carlitz module. It is denoted by $\g C$.
\medskip
{\bf Example 5.16.2.} Forms of Drinfeld module of rank 1.
\medskip
Let $P\in \n F_q(\th)^*$. We consider a Drinfeld module of rank 1 (denoted by $\g C_P$) defined by the equation
$$Te=\th e +P\tau  e$$

It is defined over $\n F_q(\th)$ (we do not give here the exact definition of "defined over" and (see below) "isomorphic over", it is clear). The above arguments show that $\g C$ is isomorphic to $\g C_P$ over $\p$. But they are not isomorphic over $\n F_q(\th)$ unless $P^{\frac1{q-1}}\in\n F_q(\th)^*$. So, the Carlitz module has forms (or twists) over $\n F_q(\th)$ (if $q\ne2$). The situation is completely analogous to the number field case, where, for example, an elliptic curve $E$ defined by the equation $y^2=x^3+ax+b$, $a,\ b \in \n Q$, has twists $E_d$ defined by the equations $dy^2=x^3+ax+b$, where $d\in \n Z$ is squarefree. $E_d$ is isomorphic to $E$ over $\n Q(\sqrt{d})$ (and hence over $\n C$) but not over $\n Q$.

Properties of $\g C_P$ can differ from the properties of $\g C$. See, for example, [GL16] for their $L$-functions.

Clearly forms exist for other Anderson t-motives, not only for $\g C$.
\medskip
{\bf 5.17. Tensor products.} Let $M_1$, $M_2$ be Anderson t-motives. We can consider their tensor product $M_1\otimes_{\p[T]}M_2$ over
$\p[T]$. It is clear that it is free of dimension $r_1r_2$ over $\p[T]$. Let us define the $\tau$-action on $M_1\otimes_{\p[T]}M_2$ by the formula $$\tau(m_1\otimes m_2):= \tau(m_1) \otimes \tau(m_2)\eqno{(5.17.1)}$$

{\bf Theorem 5.17.2.} Let $M_1$, $M_2$ be pure t-motives. Then (5.2.2), (5.2.3) hold for $M_1\otimes_{\p[T]}M_2$ with this $\tau$-action, hence it is an Anderson t-motive. It is also pure. Formula for $n$ of $M_1\otimes M_2$ ([G96], 5.7.2, (3)): $w(M_1\otimes M_2)=w(M_1)+w(M_2)$, i.e. 
$$n=n_1r_2+n_2r_1.$$

{\bf 5.17.3.} The same is true for those t-motives which can be written as appropriate extensions of pure t-motives. They are called mixed in [HJ], Definition 3.5b. Namely, $M$ is mixed if there is an increazing filtration 
$$0=M_0\subset M_1\subset M_2\subset \dots \subset M_k=M$$
of $M$ such that all quotients $M_i/M_{i-1}$ are pure of increasing weights $w(M_i/M_{i-1})$. 
\medskip
Example (a version of [HJ], Example 3.9): Let $\g C_2$ be the Drinfeld module of rank 2 such that its $a_1$ from (5.6) is 0. Let $M$ be a t-motive that enters in an exact sequence 
$$0\to \g C\to M\to \g C_2\to0\eqno{(5.17.4)}$$
Arguments similar to [HJ], Example 3.9 show that $M$ is not mixed (if $M\ne \g C\oplus\g C_2$), because $w(\g C)=1>w(\g C_2)=\frac12$. 
\medskip
From another side, let $M$ be a t-motive that enters in an exact sequence 
$$0\to \g C_2\to M\to \g C\to0\eqno{(5.17.5)}$$
These $M$ are mixed. 
\medskip
See Section 13 for an application of this notion. 
\medskip
{\bf Remark. }The nilpotent operator $N$ of $M_1\otimes M_2$ is not 0 even if $N$ of both $M_1$, $M_2$ are 0.
\medskip
Analogously, we can consider multiple tensor products, symmetric and external tensor powers.
\medskip
{\bf Example 5.17.6.} Let us consider $\g C^{\otimes n}$ --- the $n$-th tensor power of the Carlitz module (see [AT90] for details). For the $i$-th factor of $\g C^{\otimes n}$ we denote $e$ from (5.16.1) by $e_i$, hence the only element of a basis of $\g C^{\otimes n}$ over $\p[T]$ is $e_1\otimes e_2\otimes ... \otimes e_n$. We denote $e_1\otimes e_2\otimes ... \otimes e_n$ by $\g e_1$ and for $i=0,\dots, n-1$ we denote $(T-\th)^i\g e_1$ by $\g e_{i+1}$. According (5.16.1) and (5.17.1), we have $\tau(\g e_1)=(T-\th)^n\g e_1$.
Finally, we denote the matrix column $(\g e_1, \g e_2, \dots, \g e_n)^t$ by $\g e_*$. 
\medskip
Let $\ve_{ij}$ denotes the matrix whose $(i,j)$-th element is 1 and all other elements are 0.
\medskip
{\bf Exercise.} Elements $\g e_1, \g e_2, \dots, \g e_n$ form a basis of $\g C^{\otimes n}$ over $\p\{\tau\}$. The multiplication by $T$ in this basis is defined by the formula $$T\cdot \g e_*=(\th I_n+N)\g e_* + \ve_{n1}\tau(\g e_*)$$ where
$N=\left(\matrix 0&1&0&0&\dots&0&0
\\ 0&0&1&0&\dots&0&0
\\ &&&\dots
\\ 0&0&0&0&\dots&0&1 \\ 0&0&0&0&\dots&0&0 \endmatrix\right)$.
\medskip
{\bf Example 5.17.7.} Let $M_1$, $M_2$ be Drinfeld modules of ranks $r_1$, $r_2$. A description of $M_1\otimes M_2$ is the following. Let $e_1$, $e_2$ be $e$ from (5.6) for $M_1$, $M_2$ respectively. We have: elements
\medskip
$e_1\otimes e_2$, $e_1\otimes \tau e_2$, $e_1\otimes \tau^2 e_2$, ... , $e_1\otimes \tau^{r_2-1} e_2$,
\medskip
$\tau e_1\otimes e_2$, $\tau^2 e_1\otimes e_2$, ... , $\tau^{r_1-1} e_1\otimes e_2$, \ \ $(T-\th) (e_1\otimes e_2)$
\medskip
are a $\p\{\tau\}$-basis of $M_1\otimes M_2$. Equality $T( e_1\otimes e_2)= \th( e_1\otimes e_2)+(T-\th) (e_1\otimes e_2)$ shows that for this basis we have $N=\ve_{1,r_1+r_2}$. Details are left as an exercise for the reader.
\medskip
{\bf 5.18. Exterior powers of Drinfeld modules.}
\medskip
There exists the only case when a tensor power of Anderson t-motives having $N=0$ is also a t-motive having $N=0$. These are exterior powers of Drinfeld modules (and their duals, see Section 13). Namely, if $M$ is a Drinfeld module of rank $r$ then its $k$-th exterior power $\la^k(M)$ is an Anderson t-motive of rank $\binom{r}{k}$, dimension $\binom{r-1}{k-1}$ having $N=0$. For a proof see for example [GL09]; the corresponding fact for lattices (see Sections 8, \ 10) is an exercise.
\medskip
According the analogy between Anderson t-motives and abelian varieties with MIQF, we have the following proposition. Let $A$ be an abelian variety with MIQF of dimension $r$ and signature $(1,r-1)$. There exists an abelian variety with MIQF (denoted by $\la^k(A)$) of dimension $\binom{r}{k}$ and signature $(\ \ \binom{r-1}{k-1},\binom{r-1}{k}\ \ )$ such that the lattice of $\la^k(A)$ mentioned in (1.8.3) is obtained from the lattice of $A$ using the same construction which gives the lattice of $\la^k(M)$ starting from the lattice of $M$.\footnotemark \footnotetext{Prof. Claire Voisin told to the second author that experts in abelian varieties knew the construction of $\la^k(A)$ for an abelian variety with MIQF $A$; apparently, it is not published. But she did not know that this construction (which looks quite artificial for the case of abelian varieties with MIQF) comes from a very natural construction for Drinfeld modules.}
\medskip
{\bf 5.19.} It is natural to consider not only the tensor product of t-motives, but also their Hom over $\p[T]$. But here we have a difficulty. The standard action of $\tau$ on Hom group is the following. Let $\vf\in \Hom_{\p[T]}(M_1,M_2)$. By definition, we have $$\tau(\vf)(m_1):=\tau(\vf(\tau^{-1}(m_1)))\eqno{(5.19.1)}$$ But $\tau^{-1}$ is not necessarily defined.
\medskip
A way to get rid from this difficulty is to use a version of a definition of t-motives. See, for example, a definition of [HJ]; for this definition Hom always exists.
\medskip
Another way is to consider not all t-motives. For example, let $M$ be a t-motive having $N=0$. It turns out that for many $M$ its dual $$M':= \Hom_{\p[T]}(M,\g C)\eqno{(5.19.2)}$$ exists as a t-motive according definition 5.2; the action of $\tau$ is defined by (5.19.1)). See Section 12 for moredetails; here we indicate only that the dimension of $M'$ is $r-n$.
\medskip
If $N$ satisfies $N^\g m=0$ then for many such $M$ there exists $${M'}^\g m:= \Hom(M,\g C^\g m)\eqno{(5.19.3)}$$ which is called the $\g m$-th dual of $M$. See the end of Section 14 for its properties.
\medskip
{\bf 5.20. Analogy with a ring of differential operators.}
\medskip
Relations (5.1) can be described as follows. Let $R:=\p[T]$, and for $P=\sum a_iT^i\in R$ let $P^{(1)}:=\sum a_i^qT^i$. We have $\p[T,\tau]=R[\tau]$, and (5.1) is
$$\forall\ P\in R\ \ \ \tau P=P^{(1)}\tau\eqno{(5.20.1)}$$

Let now $Fun$ be a ring of infinitely derivable functions (for example, of functions from $\n R$ to $\n R$), and $D$ the operator of derivation on $Fun$. Let $Fun[D]$ be the ring of differential operators on $Fun$. Its elements are sums
$$F:=f_0+f_1D+f_2D^2+...+f_rD^r.\eqno{(5.20.2)}$$
The action of this operator $F$ on $g\in Fun$ is given y the formula 
$$F(g)=(f_0+f_1D+f_2D^2+...+f_rD^r)(g)=f_0\cdot g + f_1\cdot g' + f_2\cdot g'' + ... + f_r\cdot g^{(r)}.$$
The analog of (5.20.1) is (the multiplication is in $Fun[D]$)
$$D \cdot f = f' +f\cdot D\eqno{(5.20.2)}$$
It is checked immediately: $D \cdot f(g)=(fg)'$ and $(f'+f\cdot D)(g)=f'\ g + f\ g'$. 
\medskip
The analogy between (5.20.1) and (5.20.2), and hence between $R[\tau]$ and $Fun[D]$, is a source of some results. The initial reference is [Mu78]; see also [G95], Remark 2.5. Also, this analogy is used for the definition of holonomic sequences in $\p$, see Section 16. 

\medskip

{\bf 6. Action of $\n F_q[T]$ on $\p^n$ associated to an Anderson t-motive.}
\medskip
First, we consider a case of a Drinfeld module $M$ defined by (5.6). Associated is an action of
$\n F_q[T]$ on $\p$. The actions of the elements $a\in \n F_q\subset \n F_q[T]$, resp. $T\in \n F_q[T]$ are given by the formulas (here $x\in \p$):

$$a(x)=ax, \hbox{ resp. } T(x)=\theta x + a_1 x^q +a_2 x^{q^2}+\dots + a_r x^{q^r}\eqno{(6.0)}$$

The definition of the action of other elements of $\n F_q[T]$ is defined by the formulas:
\medskip
$(P_1+P_2)(x)=P_1(x)+P_2(x)$, $P_1P_2(x)=P_1(P_2(x))$. For example, $T^2(x)=T(T(x))$.
\medskip
It is easy to see that we have a formula: if $P\in \n F_q[T]$ and
\medskip
$Pe=(b_0+b_1\tau+b_2\tau^2+\dots +b_{k-1}\tau^{k-1}+b_k\tau^k)e$ in $M$ (here $b_i\in\p$),
\medskip
then $P(x)=b_0x + b_1 x^q +b_2 x^{q^2}+\dots + b_k x^{q^k}$.
\medskip
{\bf Example.} Let $\g C$ be the Carlitz module. Its action:
\medskip
$T(x)=\th x+x^q$, $T^2(x)=T(\th x+x^q)=\th^2 x+(\th+\th^q)x^q+x^{q^2}$,
$(T^2+T)(x)=(\th^2 +\th)x+(\th^q+\th+1)x^q+x^{q^2}$ etc.
\medskip

This action defines a structure of a $\n F_q[T]$-module on $\p$, because the formulas
\medskip
$P(x_1+x_2)=P(x_1)+P(x_2)$, $P(ax)=aP(x)$ where $a\in \n F_q$, trivially hold (here $P\in \n F_q[T]$). This module is denoted by $E(M)$.
\medskip
{\bf Remark 6.1.} Initially, Drinfeld gave in [D76] exactly this definition, this explains the word "module" in the terminology "Drinfeld module". This definition is equivalent to the definitions 5.2, 5.3 above, but definitions 5.2, 5.3 are more comprehensible for beginners.
\medskip
{\bf Example 6.2.} Let us indicate how to understand the definition of Drinfeld module given in [G96], Definition 4.4.2. {\bf A}, $\Cal F$ and $\iota$ of [G96] are respectively $\n F_q[T]$, $\p$ and an inclusion $\n F_q[T] \to \p$ sending $T$ to $\th$. A map $\phi:$ {\bf A} $\to \Cal F\{\tau\}$ of [G96] corresponding to a Drinfeld module defined by (5.6) is defined by
$$\phi(T)=\th+a_1\tau+a_2\tau^2+\dots +a_{r-1}\tau^{r-1}+a_r\tau^r\in \Cal F\{\tau\}$$
Condition [G96], 4.4.2, (1) is (5.7) of the present paper, and [G96], 4.4.2, (2) is a condition of non-triviality.
\medskip
Really, Goss considers a slightly more general case. First, his {\bf A} can be a ring bigger than $\n F_q[T]$.
\medskip
Second, $\iota$ can have a non-trivial kernel $\wp$ (a finite characteristic case, see [G96], 4.4.1). These Drinfeld modules in finite characteristic are useful for a description of reductions.
\medskip
Let us define $E(M)$ for an Anderson t-motive $M$  (see [G96], proof of 5.6.3). Let us consider $\p$ as a $\p\{\tau\}$-module defined by the formula $\tau(x)=x^q$ where $x\in\p$. 
\medskip
{\bf Definition 6.3.} $E(M)=\Hom_{\p\{\tau\}}(M,\p)$.
\medskip
The functor $M \mapsto E(M)$ is contravariant. According a general formalism of homological algebra, $E(M)$ is a module over $\n F_q$ --- the center of $\p\{\tau\}$.  Particularly, there is no canonical structure of $\p$-module on $E(M)$.
\medskip
Further, since $M$ is a $\p[T]$-module, $E(M)$ is a module over $\n F_q[T]$. Namely, for $\vf: M\to \p$ we have $T\vf(m):=\vf(Tm)$. 
\medskip
{\bf 6.3.1.} Since later (Sections 9 and further) the symbol $T$ will have different meanings, we shall use notation $T_E$ to indicate the action of $T$ on $E(M)$. 
\medskip
{\bf 6.4.} Nevertheless, it is convenient to have a coordinate description of $E(M)$ and an explicit description of the action of $T$ on $E(M)$. We fix the following notation: $\{e\}$ will mean a basis of $M$ over $\p\{\tau\}$, written as a vector column: $\{e\}=(e_1, \dots, e_n)^t$. 
\medskip
We use notations of (5.9) for the action of $T$ on $\{e\}$. 
\medskip
For $x\in E(M)$ we denote $x_i:=x(e_i)$ ($x$ is considered as a map from $M$ to $\p$). The column vector $\{x\}=\{x_1, \dots, x_n\}^t\in \p^n$ is denoted by $\g i_{\{e\}}(x)$, i.e. there is an isomorphism $\g i_{\{e\}}: E(M)\to\p^n$ (coordinates on $E(M)$). The action of $T$ on $E(M)$ in these coordinates becomes the following action on $\p^n$ via $\g i_{\{e\}}$ (here $\{x\}=\g i_{\{e\}}(x)$ ):
$$\g i_{\{e\}}(T(x))=A_0\{x\}+A_1\{x\}^{(1)}+A_2\{x\}^{(2)}+\dots +A_{k-1}\{x\}^{(k-1)}+A_k\{x\}^{(k)}\eqno{(6.5)}$$
Here, as above for matrices, $\{x\}^{(i)}:=\{x_1^{q^i}, \dots, x_n^{q^i}\}^t$, and $A_i$ are from (5.9). This formula is checked immediately. 
\medskip
We denote $\underline{A}:=\sum _{i=0}^kA_i\tau^i\in M_{n\times n}(\p\{\tau\})$. Further, for $\{x\}=(x_1,\dots, x_n)^t\in \p^n$ we denote 
$$\underline{A}_{\p}(\{x\}):=A_0\{x\}+A_1\{x\}^{(1)}+...+A_k\{x\}^{(k)}\eqno{(6.6)}$$
(and analogously $\underline{W}_{\p}(\{x\})$ for $W$, see 6.8 below). Hence, in this notation (6.5) becomes (here we neglect the isomorphisms $\g i_{\{e\}}$): $$T(\{x\}):=\underline{A}_{\p}(\{x\}).\eqno{(6.6.1)}$$ 
\medskip
More exactly, if do not neglect $\g i_{\{e\}}$, (6.5) is a commutative diagram: 
$$\matrix E(M)&\overset{\g i_{\{e\}}}\to{\to}&\p^n\\ \\ T_E\downarrow&& \downarrow \underline{A}_{\p}\\ \\ E(M)&\overset{\g i_{\{e\}}}\to{\to}&\p^n\endmatrix\eqno{(6.7)}$$

{\bf 6.8.} Let us consider the change of basis. Let $\{e'\}=(e'_1, \dots, e'_n)^t$ be another basis of $M$ over $\p\{\tau\}$ and $\underline{W}=\sum_{i=0}^kW_i\tau^i\in GL_n(\p\{\tau\})$ (here $W_i\in GL_n(\p)$) the matrix of change of basis: 
$$\{e'\}=W_0\{e\}+W_1\tau(\{e\})+...+W_k\tau^k(\{e\})$$
$\g i_{\{e\}}$, $\g i_{\{e'\}}$ enter in the following commutative diagram: 
$$\matrix E(M)&\overset{\g i_{\{e\}}}\to{\to}&\p^n\\ \\ id\downarrow&& \downarrow \underline{W}_{\p}\\ \\ E(M)&\overset{\g i_{\{e'\}}}\to{\to}&\p^n\endmatrix\eqno{(6.9)}$$

{\bf Remark.} For $n>1$ there is no canonical structure of $\p$-module on $E(M)$. Really, the only way to try to define this structure is the following. Let $x\in E(M)$ and $c\in \p$. We define $cx$ via $\g i_{\{e\}}$, namely $cx:=cx_{\{e\}}=\g i_{\{e\}}^{-1}(c\cdot \g i_{\{e\}}(x))$. But if $\exists \ i>0$ such that $W_i\ne0$ then $cx_{\{e\}}\ne cx_{\{e'\}}$.
\medskip
For $n=1$ (Drinfeld modules) we have $\underline{W}\in GL_1(\p\{\tau\})=GL_1(\p)$, hence for $n=1$ there exists a canonical structure of $\p$-module on $E(M)$. 
\medskip
{\bf 6.10.} Also, we can give an invariant definition of $E(M)$ in coordinates. Let $\g B$ be the set of $\p\{\tau\}$-bases of $M$ (a principal homogeneous space over $GL_n(\p\{\tau\})$). $E(M)$ is the quotient set of $\p^n \times \g B$ by the equivalence relation coming from (6.9): a pair $(\{x\}, \{e\})$ is equivalent to a pair $(\{x'\}, \{e'\})$ iff $\{x'\}=\underline{W}_{\p}(\{x\})$. 
\medskip
{\bf 6.11. Torsion points and Tate modules.} Let $M$ be an Anderson t-motive and $P\in \n F_q[T]$. We define $M_P$ --- the set of $P$-torsion points of $M$ --- as follows, see [G96], Proposition 5.6.3. It is the following subset of
$E(M)$ (not of $M$ itself!)
\medskip
$M_P: =\{ x\in E(M)\ \ |\ \ P(x)=0\}$.
\medskip
Choosing $\{e\}$ and identifying $E(M)$ with $\p^n$ via $\g i_{\{e\}}$, we can consider $M_P$ as a subset of $\p^n$. 
\medskip
{\bf 6.12. Example:} Let $P=T$ be the simplest irreducible polynomial, and $M$ a Drinfeld module defined by (5.6). Then $M_T\subset \p$, it is a set of the roots
of the following polynomial:
\medskip
$\theta x + a_1 x^q +a_2 x^{q^2}+\dots + a_r x^{q^r}$
\medskip
This is a $\n F_q$-vector space of dimension $r$ --- a phenomenon that never occurs in characteristic 0 !
\medskip
For any Anderson t-motive $M$ we have: $M_T$ is an abelian group, moreover, a $\n F_q$-module of dimension $r$. Analogously,
for any $P\in \n F_q[T]$ we have: $M_P$ is a free $\n F_q[T]/P$-module of dimension $r$, see [G96], Corollary 5.6.4.
\medskip
{\bf 6.13. Definition.} Let $\g L$ be a monic irreducible polynomial in $\n F_q[T]$, i.e. a finite place of $\n F_q(T)$. The $\g L$-Tate module of $M$ is:
$$T_\goth L(M):=\underset{\underset{k\to\infty}\to{\longleftarrow}}\to{\hbox{lim}} M_{\goth L^k}$$
\medskip
Namely, $x\in T_\goth L(M)$ is a sequence $(\dots, x_3, x_2, x_1, x_0=0)$, $x_i\in \p^n$ such that
$\goth L(x_i)=x_{i-1}$.
\medskip
We have: $T_\goth L(M)=(\n F_q[\theta]_\goth L)^r$ (here $\n F_q[\theta]_\goth L$ is the $\g L$-adic completion of $\n F_q[\theta]$). 
\medskip
Let $M$ be defined over a field $K\supset \n F_q(\th)$ (or, more exactly, we consider $M$ over $K$ with a fixed $K$-structure). The absolute Galois group Gal($K$) acts on $T_\goth L(M)$ for any $\g L$.
An analog of the Weil theorem holds for it. See Section 15.7, and [G96], Section 4.10 for details and proofs.
\medskip
Hence, we have the following table for Tate modules and Galois groups of our objects. Problem: define an object giving the fourth line!
\medskip
\newpage
{\bf \centerline{Table 6.14}}
\settabs 6 \columns
\medskip
\medskip
\+ &&& Tate module && Galois group\cr
\medskip
\+Abelian varieties over $\n Q$ &&& $\n Z_l^{2g}$ && Gal($\n Q$)\cr
\medskip
\+Abelian varieties over $\n F_q(\theta)$ &&& $\n Z_l^{2g}$ && Gal($\n F_q(\theta)$)\cr
\medskip
\+Anderson t-motives over $\n F_q(\theta)$ &&& $\n F_q[\theta]_\goth L^{r}$ && Gal($\n F_q(\theta)$)\cr
\medskip
\+? ? ? ? ? ? ? ? ? ? ? ?  &&& $\n F_q[\theta]_\goth L^{r}$ && Gal($\n Q$)\cr
\medskip
\medskip
{\bf 7. Drinfeld modules over finite fields; reductions}
\medskip
We give here only a descriptions of involved objects, some simple proofs and examples. For detailed proofs see [G96], Chapter 4; [Pa23], Chapter 4. We give here a table of notations of the present paper, of [G96] and of [Pa23]. 
\medskip
{\bf \centerline{Table of correspondence of notations}}
\medskip
\settabs 18 \columns
\+Present paper &&&& $M$&$q$&$r$&$\g p$&$d$&$s$&$q^s$&$t$&$\tau^s$&$m$&$e_{\infty1}$&$f_{\infty1}$&$e_{\infty2}$&$f_{\infty2}$&$ $&$ $&$ $&$ $\cr
\medskip
\+[G96]&&&&$\phi$&$r$&$d$&$\wp$&$ $&$s$&$q$&$t$&$F$&$[E$:{\bf k}]\cr
\medskip
\+[Pa23]&&&&$\phi$&$q$&$r$&$\g p$&$d$&$n$ &$q^n$&$\frac{r}{m_1}$&$\pi$ &$m_1$&$e_{1}$&$f_{1}$&$e_{2}$&$f_{2}$&\cr
\medskip
\medskip
\settabs 18 \columns
\+Present paper &&&& $s/t$&$e_{\g p1}$&$f_{\g p1}$&$e_{\g p2}$&$f_{\g p2}$&$M(\n F_q(T))$&&&$\n F_q(\tau^s)$&&&$M(\n F_q(T))(\tau^s)$\cr
\medskip
\+[G96]&& &&$s/t$&& &&&{\bf k}& &&$\n F_r(F)$&& &$E=$ {\bf k}$(F)$\cr
\medskip
\+[Pa23]&&&&$m_2$&$e_{2}$&$f_{2}$&$e_{1}$&$f_{1}$&$F=\n F_q(T)$&&&$K=\n F_q(\pi)$&&&$\tilde F=F(\pi)$\cr
\medskip
\medskip
Let $\g p$ be an irreducible polynomial in $\n F_q[T]$ of degree $d$, $\g p=\sum_{i=0}^dc_iT^i$, and $s$ a mulptiple of $d$. Hence, we have a map 
$$\iota:\ \n F_q[T]\to\n F_q[T]/\g p=\n F_{q^d}\hookrightarrow\n F_{q^s}$$
We denote $\th_0:=\iota(T)\in \n F_{q^s}$. 
A Drinfeld module of rank $r$ of characteristic $\g p$ over $\n F_{q^s}$ is a map $M: \n F_q[T]\to\n F_{q^s}\{\tau\}$ defined by $$M(T)=\th_0+a_1\tau+...+a_r\tau^r\eqno{(7.0)}$$
where $a_i\in \n F_{q^s}, \ a_r\ne0$. 

Let us consider $M(\g p)\in \n F_{q^s}$, $M(\g p)=\sum_{i=0}^{dr}\be_i\tau^i$. Its free term $\be_0=\sum_{i=0}^dc_i\th_0^i$ is 0. Let $\bar h$ be the minimal number such that $\be_{\bar h}\ne0$. 
\medskip
{\bf Proposition 7.1.} $\bar h$ is a multiple of $d$.
\medskip
{\bf Proof.} The group of $\g p$-torsion points of $E(M)$ (considered as a set) is a vector space over $\n F_q[T]/\g p=\n F_{q^d}$. Its order is $q^{dr-\bar h}$. $\square$
\medskip
{\bf Definition 7.2.} The number $h:=\bar h/d$ is called the height of $M$. If $h=1$ then $M$ is called ordinary, if $h=r$ then $M$ is called supersingular (obviously $1\le h  \le r$). 
\medskip
{\bf Example 7.3.} Let $r=1$ (the Carlitz module), $d=s=3$, $\g p=T^3+c_2T^2+c_1T+c_0$. We have 
$$M(\g p)=(\th_0+\tau)^3+c_2(\th_0+\tau)^2+c_1(\th_0+\tau)+c_0=$$ $$(\th_0^3+c_2\th_0^2+c_1\th_0+c_0)+(\th_0^{}+\th_0^{2}+\th_0^{q+1}+\th_0^{2q}+(\th_0+\th_0^{q})c_2+c_1)\tau+(\th_0+\th_0^{q}+\th_0^{q^2}+c_2)\tau^2+\tau^3$$ Since $c_{3-i}=(-1)^i\s_{i}(\th_0, \ \th_0^{q}, \ \th_0^{q^2})$ (symmetric polynomials), we get: $M(\g p)=\tau^3$. 
\medskip
We shall describe some endomorphisms rings related to $M$. They are extensions of $\n F_q[T]$. It is simpler to consider their tensor products by $\n F_q(T)$, i.e. by analogy with abelian varieties, we consider $\End^0$ instead of End itself. 
\medskip
All below rings are subrings of $\n F_{q^s}\{\tau\}$, and their tensor products by $\n F_q(T)$ are (skew) fields which are subfields of $\n F_{q^s}(\tau)$ --- its skew field of fractions. 
\medskip
For any $a\in \n F_q[T]$ there is an element $M(a)\in \n F_{q^s}\{\tau\}$, defined according (7.0). Hence, we have a field $M(\n F_q(T))$ which is a subfield of $\End^0(M)$. Further, we have the Frobenius map $\tau^s$, it is in the center of $\n F_{q^s}(\tau)$. It gives us an inclusion $\n F_{q}(\tau^s)\subset \End^0(M)$. Hence, we have their compositum $M(\n F_q(T))(\tau^s)\subset \End^0(M)$ and a diagram (7.4) (all maps in (7.4) are inclusions, numbers over arrows are dimensions, see below for their definitions and related formulas):

$$\matrix M(\n F_q(T))\\ & \overset{m}\to{\searrow}\\  &&M(\n F_q(T))(\tau^s)&\overset{t^2}\to{\hookrightarrow}&\End^0(M)&\overset{s/t}\to{\hookrightarrow}&\n F_{q^s}(\tau) 
\\ & \overset{s/t}\to{\nearrow} \\ \n F_{q}(\tau^s) \endmatrix\eqno{(7.4)}$$

We have: 

{\bf 7.4.1.} $\End^0(M)$ is a division algebra, and $M(\n F_q(T))(\tau^s)$ is its center;
\medskip
{\bf 7.5.} $r=mt$. 
\medskip
{\bf Example 7.6.} Let $r=d=s=2$, $m=2, \ t=1$ (the ordinary case, see below). Let us find $\Cal F$ --- the characteristic polynomial of Frobenius, i.e. the equation satisfied by $\tau^2$ over $M(\n F_q(T))$. Let as above $\g p=T^2+c_1T+c_0$, and a Drinfeld module $M$ of rank 2 is defined by $$M(T)=\th_0+a_1\tau+\tau^2$$
where $a_1\in \n F_{q^2}$. We have (here and below $N$, $Tr$ are the norm, trace from $\n F_{q^2}$ and (below) $\n F_{q^3}$ to $\n F_q$):
$$\tau^4+M(-2T-c_1-N(a_1))\ \tau^2+M(\g p)=0\eqno{(7.6.1)}$$
Really,

$M(\g p)=(\th_0+a_1\tau+\tau^2)^2+c_1(\th_0+a_1\tau+\tau^2)+c_0=(2\th_0+N(a_1)+c_1)\tau^2+2a_1\tau^3+\tau^4$, 
\medskip
$(-2T-c_1-N(a_1))\tau^2=(-2\th_0-N(a_1)-c_1)\tau^2-2a_1\tau^3-2\tau^4$, hence the result. 
\medskip
{\bf Example 7.7.} Let now $r=3, \ d=s=2$. The only possibility for $m, \ t$ is $m=3, \ t=1$. An analog of (7.6.1) has the form $$\tau^6+M(C_2)\ \tau^4+M(C_1)\ \tau^2+M(\g p)=0\eqno{(7.7.1)}$$
where $C_2, \ C_1\in \n F_q[T]$ are indefinite coefficients. What could be their degrees $\g d_2, \ \g d_1$? Unlike the case of Example 7.6 where it is clear that $\g d_1=1$, i.e. all terms of (7.6.1) have degree 4, here the situation looks mysterious. Really, the degrees in $\tau$ of $M(C_2)\ \tau^4$, resp. of $M(C_1)\ \tau^2$ are $3\g d_2+4$, resp. $3\g d_1+2$, hence the only possibility to satisfy (7.7.1) is $\g d_2=0, \ \g d_1=1$. From the first sight, we have insufficiently many indefinite coefficients in $\n F_q$ to satisfy (7.7.1). But let us calculate: 

$$M(\g p)=(\th_0+a_1\tau+a_2\tau^2+\tau^3)^2+c_1(\th_0+a_1\tau+a_2\tau^2+\tau^3)+c_0=$$ $$=(2\th_0a_2+N(a_1)+c_1a_2)\tau^2+a_1Tr(a_2)\tau^3+(a_2^2+Tr(a_1))\tau^4+Tr(a_2)\tau^5+\tau^6$$

We see that for $C_2=N(a_2)-Tr(a_1)$, $C_1=-Tr(a_2)T+[Tr(\th_0a_2^q)-N(a_1)]$ 

(7.7.1) is satisfied. 
\medskip
{\bf Example 7.8.} Let $r=2, \ d=s=3$. The only possibility for $m, \ t$ is $m=2, \ t=1$. Here we have the same phenomenon. An analog of (7.7.1) has the form $$\tau^6+M(C_1)\ \tau^3+M(\g p)=0\eqno{(7.8.1)}$$ and $C_1$ must be of degree 1. We have
$$M(\g p)=(\th_0+a_1\tau+\tau^2)^3+c_2(\th_0+a_1\tau+\tau^2)^2+c_1(\th_0+a_1\tau+\tau^2)+c_0=$$ $$=(\th_0a_1+\th_0a_1^{q^2}-\th_0^qa_1^{q^2}-\th_0^{q^2}a_1+N(a_1))\tau^3+a_1Tr(a_1)\tau^4+Tr(a_1)\tau^5+\tau^6$$
For $C_1=-Tr(a_1)T+[Tr(\th_0a_1^q)-N(a_1)]$ (7.8.1) is satisfied. 
\medskip
Maybe there exists an explicit formula for coefficients of $\Cal F$ for higher $r, \ d, \ s$? Its finding would be a good exercise for a student. 
\medskip
{\bf 7.9. Valuations over $\g p$.} 
\medskip
{\bf Theorem 7.9.1.} There are valuations $v_{\g p1}, \ v_{\g p2}, \ v_{\g p3}$ of the three leftmost fields of the diagram (7.4), defined by (7.9.3), see the below diagram. Here $e_{\g p1}, \ f_{\g p1}, \ e_{\g p2}, \ f_{\g p2}$ denote the ramification indices and degrees of inertia of the corresponding extensions. They satisfy the properties 7.9.4 below.

$$\matrix\underset{v_{\g p1}}\to{M(\n F_q(T))}\\ & \underset{e_{\g p1}, f_{\g p1}}\to{\overset{m}\to{\searrow}}\\  &&\underset{v_{\g p3}}\to{M(\n F_q(T))(\tau^s)}
\\ & \underset{e_{\g p2}, f_{\g p2}}\to{\overset{s/t}\to{\nearrow}} \\ \underset{v_{\g p2}}\to{\n F_{q}(\tau^s)} \endmatrix\eqno{(7.9.2)}$$

{\bf 7.9.3.} $v_{\g p1}$, $v_{\g p2}$ are defined uniquely by the formulas $$v_{\g p1}(M(\g p))=1 \eqno{(7.9.3.1)}$$ $$v_{\g p2}(\tau^s)=1\eqno{(7.9.3.2)}$$

{\bf 7.9.4.} We have: $v_{\g p3}$ is the only valuation over $v_{\g p2}$. Hence, $e_{\g p2} f_{\g p2}=s/t$. $v_{\g p3}$ is over $v_{\g p1}$, but not necessarily the only one over it. Hence, $e_{\g p1} f_{\g p1}\le m$. 
\medskip
We have relations between the above numbers: 
\medskip
{\bf 7.9.5.} $d\ f_{\g p1}= f_{\g p2}$ ([Pa23], (4.1.6); caution: non-concordance of indices, see the above table of correspondence of notations). 
\medskip
Really, the dimension over $\n F_q$ of the residue field of $v_{\g p1}$ is $d$, hence the dimension over $\n F_q$ of the residue field of $v_{\g p3}$ is $d\ f_{\g p1}$ (upper arrow), $f_{\g p2}$ (lower arrow). 
\medskip
{\bf 7.9.6.} $h=t\ e_{\g p1}\ f_{\g p1}$ (recall that $h$ is the height).
\medskip
Really, by the definition of the height, $$M(\g p)=b_{hd}\tau^{hd}+...+\tau^{rd}\in M(\n F_q(T))(\tau^s)$$ where $0\ne b_{hd}\in \n F_{q^s}$ is a coefficient. We have 

$$v_{\g p3}(b_{hd}\tau^{hd}+...+\tau^{rd})=e_{\g p1}$$ (from (7.9.3.1) and the definition of $e_{\g p1}$);

$$v_{\g p3}(\tau^s)=e_{\g p2}$$ (from (7.9.3.2) and the definition of $e_{\g p2}$). But 

$$\frac{v_{\g p3}(b_{hd}\tau^{hd}+...+\tau^{rd})}{v_{\g p3}(\tau^s)}=\frac{hd}s$$
i.e. $h=\frac{e_{\g p1}s}{e_{\g p2}d}$. Equalities $d\ f_{\g p1}= f_{\g p2}$ and $e_{\g p2} f_{\g p2}=s/t$ imply the result. 
\medskip
{\bf 7.10. Valuations over $\infty$.} 
\medskip
{\bf Theorem 7.10.1.} There are valuations $v_{\infty1}, \ v_{\infty2}, \ v_{\infty3}$ of the three leftmost fields of the diagram (7.4), defined by (7.10.3), see the below diagram. Here $e_{\infty1}, \ f_{\infty1}, \ e_{\infty2}, \ f_{\infty2}$ denote the ramification indices and degrees of inertia of the corresponding extensions. They satisfy the properties 7.10.4 below.

$$\matrix\underset{v_{\infty1}}\to{M(\n F_q(T))}\\ & \underset{e_{\infty1}, f_{\infty1}}\to{\overset{m}\to{\searrow}}\\  &&\underset{v_{\infty3}}\to{M(\n F_q(T))(\tau^s)}
\\ & \underset{e_{\infty2}, f_{\infty2}}\to{\overset{s/t}\to{\nearrow}} \\ \underset{v_{\infty2}}\to{\n F_{q}(\tau^s)} \endmatrix\eqno{(7.10.2)}$$

{\bf 7.10.3.} $v_{\infty1}$, $v_{\infty2}$ are defined uniquely by the formulas $$v_{\infty1}(M(\g p))=-1 \eqno{(7.10.3.1)}$$ $$v_{\infty2}(\tau^s)=-1\eqno{(7.10.3.2)}$$

{\bf 7.10.4.} We have: $v_{\infty3}$ is the only valuation over both $v_{\infty1}$, $v_{\infty2}$ (with respect to both inclusions). Hence, $e_{\infty1} f_{\infty1}=m$, $e_{\infty2} f_{\infty2}=s/t$.  . 
\medskip
{\bf 7.10.5.} We have: $f_{\infty1}=f_{\infty2}$. Really, the residue fields of both $v_{\infty1}$, $v_{\infty2}$ are $\n F_q$. 
\medskip
{\bf 7.10.6.} We have: $e_{\infty1}/e_{\infty2}=r/s$. Really, $e_{\infty1}f_{\infty1}=m=r/t$ and $e_{\infty2}f_{\infty2}=s/t$. 
\medskip
{\bf Theorem 7.11.} $\End^0(M)$ is a central division algebra over $M(\n F_q(T))(\tau^s)$, with invariants inv$(v_{\g p3})=1/t$, inv$(v_{\infty3})=-1/t$. All other invariants are 0. 
\medskip
{\bf Theorem 7.12.} For $\g l\ne\g p$ the characteristic polynomial of the Frobenius of $\n F_{q^s}$ acting on $T_\g l(M)$ is $\Cal F$. Like for the case of abelian varieties, it does not depend on $\g l$ and has coefficients in $\n F_q[T]$. 
\medskip
Analog of the formula $\#(E(\n F_p))=\Cal F(1)=1-a_p+p$: to write. See [Pa23]. 
\medskip
{\bf 7.13. Weil numbers.} See [G96], Definition 4.12.14.
\medskip
We consider the field extension $M(\n F_q(T))\hookrightarrow M(\n F_q(T))(\tau^s)$ and we extend the field of constants to $\n F_{q^s}$. Weil numbers are "abstract $\tau^s$". We use the above notations $v_{\g p1}$, $v_{\g p3}$, $v_{\infty1}$, $v_{\infty3}$. Namely, $v_{\g p3}$ and $v_{\infty3}$ are valuations of the below field $\n F_{q^s}(T)(w)$, they are over $v_{\g p1}$ and $v_{\infty1}$ respectively, and $v_{\infty3}$ is the only valuation of $\n F_{q^s}(T)(w)$ over $v_{\infty1}$. 
\medskip
{\bf  Definition 7.13.1.} Let $q,\ r, \ s$ be fixed. A Weil number corresponding to $q,\ r, \ s$ is an element $w$ of $\overline{\n F_{q^s}(T)}$ which is:
\medskip
1. $w$ is integer over $\n F_{q^s}[T]$.
\medskip
2. The dimension of $\n F_{q^s}(T)(w)$ over $\n F_{q^s}(T)$ divides $r$.
\medskip
3. $v_{\g p3}$ and $v_{\infty3}$ are the only valuations $v$ of $\n F_{q^s}(T)(w)$ such that $v(w)\ne0$.
\medskip
4. $v_{\infty3}(w)=-\frac{s}{r?}$ where the order of the residue field of $v_{\infty3}$ is $q^?$ (to write the value of ?). 
\medskip
{\bf  Theorem 7.13.2.} [G96], Theorem 4.12.15. The set of Weil numbers, up to conjugation, is the set of the Drinfeld modules of rank $r$ over $\n F_{q^s}$, up to isogeny.
\medskip
{\bf 7.14.} $\End^0(M)$ can depend on the finite field: examples. 
\medskip
{\bf 7.15.} Explicit formula for the genus of $M(\n F_q(T))(\tau^s)$ --- to write. We can apply the Riemann--Hurwitz formula for the field extension $M(\n F_q(T))\hookrightarrow M(\n F_q(T))(\tau^s)$ ( = covering of curves), because we know all ramification indices. What are singularities of this plane curve? For example, if $r$ or $d$ is 2, this curve is a hyperelliptic curve. It has only one singularity at infinity. 
\medskip
{\bf 7.16. Case of Anderson t-motives.} See [BH]. 
\medskip
Let $\iota, \ \g p, \ d, \ s$ be as above. We have an analogous definition of an Anderson t-motive of dimension $n$, rank $r$, characteristic $\g p$ over $\n F_{q^s}$: it is a map $M: \n F_q[T]\to M_{n\times n}(\n F_{q^s})\{\tau\}$ defined by $$M(T)=(\th_0+N)+A_1\tau+...+A_k\tau^k$$
where $A_i\in M_{n\times n}(\n F_{q^s})$ such that the analog of (5.2.1) holds. 
\medskip
The group of $\g p$-torsion points of $E(M)$ (considered as a set) is a vector space over $\n F_q[T]/\g p=\n F_{q^d}$. If $N=0$ then its dimension is obviously $\le r-n$. By analogy with the case of Drinfeld modules, we call ordinary such $M$ that this dimension is $r-n$. 
\medskip
{\bf Conjecture 7.16.1.} Existence of pairing between $E(M)_\g p$ and $E(M')_\g p$, as group schemes (see Section 12.1 for $M'$ --- the dual of $M$). 
\medskip
See a slightly weaker conjecture 7.17.2.
\medskip
{\bf 7.17. Reductions.} Recall that $\g p=\sum_{i=0}^dc_iT^i\in \n F_q[T]$. We consider $\g p_\th:=\sum_{i=0}^dc_i\th^i\in \n F_q[\th]$. Let $M$ be an Anderson t-motive such that all entries of $A_i$ from (5.9) belong to $\n F_q(\th)$ and are $\g p_\th$-integer. Let us consider the equation obtained from (5.9) by reduction of all entries of $A_i$ modulo
$\g p_\th$. The obtained object is called the reduction of $M$ modulo $\g p$, it is denoted by $\tilde M_\g p$. See [G96], (4.10) for the details. Roughly speaking, $M$ has a good reduction at $\g p$ if the mayor coefficient of (5.9) "does not loose its rank" after reduction (i.e. the reduced t-motive has the same $r$ as $M$ itself).
\medskip
We have also the notion of $M$ having a good ordinary reduction. We have a map of their closed $\g p$-torsion points: $$red(M, \g p): E(M)_\g p\to E(\tilde M_\g p)_\g p\eqno{(7.17.1)}$$ Its kernel has dimension $n$ over $\n F_q[T]/\g p=\n F_{q^d}$. We have an analog of Conjecture 7.16.1 for reductions (see also [GL07], Conjecture 13.4.1): 
\medskip
{\bf Conjecture 7.17.2.} $\Ker red(M, \g p)$ and $\Ker red(M', \g p)$ are dual with respect to the duality between $E(M)_\g p$ and $E(M')_\g p$.
\medskip
We reproduce here the proof of GL07], Conjecture 13.4.1, of this conjecture for a particular case: $M$ is a Drinfeld module, $\g p=T$ (hence $\g p_\th=\th$). Let a Drinfeld module $M$ be given by (5.6), $a_0=\th, \ a_r=1$: 
$$Te=(\th+a_1\tau+a_2\tau^2+\dots +a_{r-1}\tau^{r-1}+\tau^r)e$$
Condition of good ordinary reduction at $\g p=T$ means $v_\th(a_1)=0$, $v_\th(a_i)\ge0$ for $i=2,\dots, r-1$. $E(M)_T$ is the set of $x\in \p$ satisfying
$$x^{q^r}+a_{r-1}x^{q^{r-1}}+...+a_1x^q+\th x=0\eqno{(7.17.3)}$$
The dual of $M$, denoted by $M'$, is defined explicitly in Example 12.2. $E(M')_T$ is the set of $y_*=(y_1,\dots, y_{r-1})^t\in \p^{r-1}$ ($*^t$ means transposition) satisfying 
$$(\th I_{r-1}+A_1\tau+A_2\tau^2)y_*=0\eqno{(7.17.4)}$$ where $A_1, \ A_2$ are from (12.2.2). 
\medskip
The $\n F_q$-space $E(M')_T$ is $r$-dimensional. We denote its basis by $y_{1*}, \dots, y_{r*}$, and let $y_{i*}=(y_{i1},\dots, y_{i,r-1})^t$.
\medskip
Explicitly, (7.17.4) is $$\matrix \th y_1-a_1y_{r-1}^q+y_{r-1}^{q^2}=0\\ \\
\th y_2+y_1^q-a_2y_{r-1}^q=0\\ \\
\th y_3+y_2^q-a_3y_{r-1}^q=0\\
\dots \\
\th y_{r-1}+y_{r-2}^q-a_{r-1}y_{r-1}^q=0 \endmatrix \eqno{(7.17.5)}$$

The pairing between $E(M)_T$ and $E(M')_T$ is defined by the formula 
$$<x, y_*>=c(xy_{r-1}^q+x^qy_1+x^{q^2}y_2+...+x^{q^{r-1}}y_{r-1})\in \n F_q\eqno{(7.17.6)}$$
where $c\in (-\th)^{-1/(q-1)}$ (because $<x, y_*>^q=<x, y_*>$). 
\medskip
Consideration of the Newton polygon of (7.17.3) shows that $$\hbox{ for $x\in \Ker red(M, T)$ we have $v_\th(x)=\frac1{q-1}$} $$ (this is the leftmost segment of the Newton polygon of (7.17.3)). 
\medskip
Consecutive elimination of $y_1, \ y_2, \dots, y_{r-2}$ from (7.17.5) gives us the equation satisfied by $y_{r-1}$:

$$y_{r-1}^{q^{r}}-a_1^{q^{r-2}}y_{r-1}^{q^{r-1}}+\th^{q^{r-2}}a_2^{q^{r-3}}y_{r-1}^{q^{r-2}} -\th^{q^{r-2}+q^{r-3}}a_3^{q^{r-4}}y_{r-1}^{q^{r-3}}+...$$ $$...\pm \th^{q^{r-2}+q^{r-3}+...+q^2+q}a_{r-1}y_{r-1}^{q}\mp \th^{q^{r-2}+q^{r-3}+...+q^2+q+1}y_{r-1}=0 \eqno{(7.17.7)}$$

The Newton polygon of (7.17.7) has a segment of slope $-\frac1{q-1}$: $$[ \ (1, q^{r-2}+q^{r-3}+...+q+1) ; \ \ \ (q^{r-1},0)\ ]$$ The corresponding $r-1$ linearly independent elements $y_{1*}, \dots, y_{r-1,*}$ of the set of solutions of (7.17.7) have  $v_\th(y_{i, r-1})=\frac1{q-1}$ ($\forall \ i=1, \dots, r-1$). 
\medskip
Further, we have (recall that $a_0=\th$): $$y_1=\frac{a_1}{\th}y_{r-1}^{q}-\frac{\th^q}{\th^{q+1}}y_{r-1}^{q^2}$$
$$y_2=\frac{a_2}{\th}y_{r-1}^{q}-\frac{a_1^q}{\th^{q+1}}y_{r-1}^{q^2}+\frac{\th^{q^2}}{\th^{q^2+q+1}}y_{r-1}^{q^3}\eqno{(7.17.8)}$$
$$y_3=\frac{a_3}{\th}y_{r-1}^{q}-\frac{a_2^q}{\th^{q+1}}y_{r-1}^{q^2}+\frac{a_1^{q^2}}{\th^{q^2+q+1}}y_{r-1}^{q^3}-\frac{\th^{q^3}}{\th^{q^3+q^2+q+1}}y_{r-1}^{q^4}$$ etc. We get that $\forall \ y_*$ if $v_\th(y_{r-1})=\frac1{q-1}$ then $\forall \ j $ \ $v_\th(y_{j})\ge\frac1{q-1}$. Hence, these elements $y_{1*}, \dots, y_{r-1,*}$ form a basis of $\Ker red (M', T)$. 
\medskip
Since $v_\th(c)=-\frac1{q-1}$ we have that if $x\in \Ker red (M, T)$, $y_*\in \Ker red (M', T)$ (i.e. $y_*$ is a linear combination of $y_{1*}, \dots, y_{r-1,*}$), then all summands of (7.17.6) for these $x, \ y_*$ have $v_\th\ge1$, i.e. $<x,y_*>=0$. $\square$
\medskip
We can get similar calculational proofs for all $\g p$, for all t-motives defined by explicit formulas (for example, for t-motives defined by the formulas of [GL07], Section 11). 
\medskip
Clearly an analog of Conjecture 7.17.2 holds not only for the reduction map (7.17.1) on the level of $\g p$-torsion points, but also on the level of the Tate modules $T_\g p(M), \ T_\g p(M')$.
\medskip
We see that the behavior of Anderson t-motives under reduction is analogous to the behavior of abelian varieties with MIQF (is it really so? To check).
\medskip
\medskip
{\bf 8. Lattices of Drinfeld modules}
\medskip
Let $M$ be a Drinfeld module defined by (5.6).
Let us consider the below diagram where $exp=exp_M:\p\to\p$ is a map defined by the formula $$exp(z)=z+c_1z^q+c_2z^{q^2}+c_3z^{q^3}+\dots$$ where $c_i\in \p$ are some coefficients, the left vertical arrow is simply multiplication by $\th$, and the right vertical arrow is the action of $T$ on $\p$ from 6.1:

$$\matrix \p &\overset{exp}\to{\to} &\p &\\ \\ \theta\downarrow &&\downarrow & z\mapsto T(z)
\\  \\ \p &\overset{exp}\to{\to} & \p &\endmatrix \eqno{(8.1)}$$
\medskip
{\bf Remark.} For any coefficients $c_i$ the map $exp$ is $\n F_q$-linear (i.e. $exp(x+y)=exp(x)+exp(y); \ exp(cx) =c\ exp(x)$ for $c\in\n F_q$); moreover, any analytic $\n F_q$-linear $\al:\p\to\p$ is of this form. 
\medskip
{\bf Theorem 8.1.1 (Drinfeld).} For any Drinfeld module $M$ there exist numbers $c_i$ making (8.1) commutative. They can be found consecutively. They are unique.
\medskip
{\bf Example 8.2.} Finding of $c_i$ for the case of the Carlitz module $\goth C$. We have:
\medskip
$Te=(\theta +\tau)e$, hence $T(z)=\theta z+z^q$.
\medskip
In the below diagram the inner square is (8.1) for $\g C$, the external square shows the images of any $z\in \p$
(upper left) under the arrows of the inner square.

$$\matrix z&&&&\overset{exp_\g C}\to{\to} &z+c_1z^q+c_2z^{q^2}+c_3z^{q^3}+\dots \\
&\p &\overset{exp_\g C}\to{\to} &\p &&\\ \\ \downarrow & \downarrow &&\downarrow &
w\mapsto \theta w +w^q& \downarrow  \\  \\ &\p &\overset{exp_\g C}\to{\to} & \p &&
u_1=\theta(z+c_1z^q+c_2z^{q^2}+c_3z^{q^3}+\dots)\\ &&&&&+(z+c_1z^q+c_2z^{q^2}+c_3z^{q^3}+\dots)^q\\ &&&&& = \\
\theta z &&&&\overset{exp_\g C}\to{\to} & u_2=\theta z+c_1(\theta z)^q+c_2(\theta z)^{q^2}+c_3(\theta z)^{q^3}+
\dots\endmatrix $$
\medskip
\medskip
$u_1$ is $T(exp_\g C(z))$, and $u_2$ is $exp_\g C(\th z)$, hence they are equal. We get an equality of power series:
\medskip
$$\matrix &\theta z &+\theta c_1z^q&+\theta c_2z^{q^2}&+\theta c_3z^{q^3}&+\theta c_4z^{q^4}&+\dots \\ + \\ && z^q& +c_1^qz^{q^2}& +c_2^qz^{q^3}&+c_3^qz^{q^4}&+\dots\\ \\ =&\theta z &+\theta^q c_1z^q&+\theta^{q^2} c_2z^{q^2}&+\theta^{q^3} c_3z^{q^3}&+\theta^{q^4} c_4z^{q^4}&+\dots \endmatrix $$
\medskip
We get a system of equations satisfied by $c_i$:
\medskip
$$\theta c_1+1=\theta^q c_1; \ \ \ c_1=\frac{1}{\theta^q-\theta}$$
$$\theta c_2+c_1^q=\theta^{q^2} c_2; \ \ \ c_2=\frac{1}{(\theta^{q^2}-\theta^{q})(\theta^{q^2}-\theta)}\eqno{(8.2.1)}$$
$$\theta c_3+c_2^q=\theta^{q^3} c_3; \ \ \ c_3=\frac{1}{(\theta^{q^3}-\theta^{q^2})(\theta^{q^3}-\theta^q)(\theta^{q^3}-\theta)}$$
\medskip
etc. Further, $v_\infty(c_i)=iq^i$ (see Section 3 for $v_\infty$), hence $c_i\to 0$. Moreover, $\forall \ z\in \p$ we have $v_\infty(c_iz^{q^i})\to+\infty$, hence $exp_\g C(z)$ converges for all $z\in \p$.
\medskip
{\bf 8.3. Theorem (Drinfeld). 1.} For all Drinfeld modules $M$ the function $exp_M(z)$ converges for all $z\in \p$,
it is surjective, and its kernel is a lattice in $\p$ of rank $r$.
\medskip
{\bf 2.} Let us denote the above lattice by $L(M)$. The map $M \mapsto L(M)$ is a 1 -- 1 correspondence between the set of Drinfeld modules, up to isomorphism,
and the set of lattices in $\p$, up to equivalence.
\medskip
{\bf Idea of the proof:} Let $L\subset \p$ be a lattice. We associate it the following function $\wp_L: \p\to\p$ (here $L':=L-0$ ):
$$\wp_L(z)=z\prod_{\omega\in L'}(1-\frac{z}{\omega})\eqno{(8.4)}$$
\medskip
We have: $L=L(M)$ iff $exp_M=\wp_L$.
\medskip
This theorem permits us to describe the moduli space of Drinfeld modules in terms of lattices. Practically, it is the quotient space of the set of Siegel matrices of size $1\times r-1$ by the action of $GR_r(\n F_q[\th])$ defined in Section 4. We do not consider these subjects in the present survey.
\medskip
{\bf 9. Lie(M). }
\medskip
According the general formalism of Lie groups, we can define $Lie(E(M))$, or simply $Lie(M)$, as follows. We use notations of 6.3, and we define $Lie(M)$ like in 6.10. Namely, for any $\{e\}$ there exists am isomorphism $\g j_{\{e\}}: Lie(M)\to \p^n$, and the condition of concordance is that the following diagram is commutative: 

$$\matrix Lie(M)&\overset{\g j_{\{e\}}}\to{\to}&\p^n\\ \\ id\downarrow&& \downarrow W_0\\ \\ Lie(M)&\overset{\g j_{\{e'\}}}\to{\to}&\p^n\endmatrix\eqno{(9.1)}$$
Namely (notations of 6.8, 6.10) $Lie(M)$ is the quotient set of $\p^n \times \g B$ by the equivalence relation coming from (9.1): a pair $(\{x\}, \{e\})$ is equivalent to a pair $(\{x'\}, \{e'\})$ iff $\{x'\}=W_{0}(\{x\})$. 
\medskip
Since $W_0\in GL_n(\p)$ we get that $Lie(M)$ has a canonical structure of $\p$-vector space: for $x\in Lie(M)$, $c\in \p$ we have $c\cdot x:=\g j_{\{e\}}^{-1}(c\cdot \g j_{\{e\}}(x))$; it does not depend on $\{e\}$. 
\medskip
Moreover, $Lie(M)$ is a $\p[T]$-module. The multiplication by $T$ in $Lie(M)$, denoted by $T_{Lie}$, is defined by the commutative diagram ($A_0$ is  from 5.9):

$$\matrix Lie(M)&\overset{\g j_{\{e\}}}\to{\to}&\p^n\\ \\ T_{Lie}\downarrow&& \downarrow A_0\\ \\ Lie(M)&\overset{\g j_{\{e\}}}\to{\to}&\p^n\endmatrix\eqno{(9.2)}$$
Clearly $A_0$ depends on $\{e\}$, but it is checked immediately (using (9.1), (9.2)) that $T_{Lie}$ does not depend on $\{e\}$. We have $A_0=\th+N$ where $N$ (also depending on $\{e\}$) is nilpotent. We denote $N=N_{\{e\}}$ if it is necessary to emphasize its dependence on $\{e\}$. If $N=0$ then $T_{Lie}$ is simply multiplication by $\th$. 
\medskip
There exists a canonical map $exp=exp_M: Lie(M) \to E(M)$ making the following diagram commutative (recall that $T_E$ indicate the action of $T$ on $E(M)$, see (6.3.1): 
$$\matrix Lie(M)&\overset{exp}\to{\to}&E(M)\\ \\ T_{Lie}\downarrow&& \downarrow T_E\\ \\ Lie(M)&\overset{exp}\to{\to}&E(M)\endmatrix\eqno{(9.3)}$$ 

To prove its existence, we consider first $exp$ in coordinates, we denote it by $expc$ (the letter $c$ means coordinates). Namely, for any $\{e\}$ there exists a map $expc_{\{e\}}=expc_{M,\{e\}}: \p^n\to \p^n$, see (9.4b), making the following diagram commutative: 
$$\matrix \p^n&\overset{expc_{\{e\}}}\to{\to}&\p^n\\ \\ A_0\downarrow&& \downarrow \underline{A}_{\p}\\ \\ \p^n&\overset{expc_{\{e\}}}\to{\to}&\p^n\endmatrix\eqno{(9.4)}$$ 

For the case $N=0$ we have an equivalent form of (9.4): $$\matrix \p^n &\overset{expc_{\{e\}}}\to{\to} &\p^n &\\ \\ \theta\downarrow &&\downarrow & z\mapsto T(z) \\  \\ \p^n &\overset{expc_{\{e\}}}\to{\to} & \p^n &\endmatrix \eqno{(9.4a)}$$

The map $expc_{\{e\}}$ is given by the formula (here $z\in \p^n$) 

$$expc_{\{e\}}(z)= z+C_1z^q+C_2z^{q^2}+C_3z^{q^3}+\dots\eqno{(9.4b)}$$ where $C_i=C_i(\{e\})=C_i(M,\{e\})\in M_n(\p)$, and $C_0=id$ for all $M, \ \{e\}$. The matrices $C_i(\{e\})$ can be found consecutively, see [A86], (2.1.4) (see also Example 8.2). For all $M$ and $\{e\}$ the map $expc_{\{e\}}$ converges for all $z\in \p^n$.
\medskip
We define $exp$ via the following commutative diagram: 
$$\matrix Lie(M)&\overset{exp}\to{\to}&E(M)\\ \\ \g j_{\{e\}}\downarrow&& \downarrow \g i_{\{e\}}\\ \\ \p^n&\overset{expc_{\{e\}}}\to{\to}&\p^n\endmatrix\eqno{(9.5)}$$

{\bf Proposition 9.5a.} $exp$ does not depend on $\{e\}$. 
\medskip
{\bf Proof} is reduced to a careful checking of definitions. 
\medskip
Hence, diagrams (6.7), (9.3), (9.4), (9.5) are parts of the following commutative diagram: 
$$\matrix Lie(M)&&\overset{exp}\to{\to}&&E(M)&&\\ 
&\overset{\g j_{\{e\}}}\to{\searrow}&&&&\overset{\g i_{\{e\}}}\to{\searrow}&\\ 
T_{Lie}\downarrow&& \p^n&&    \underset{ \downarrow T_E}\to{\overset{expc_{\{e\}}}\to{\to}}&&\p^n\\ \\ 
Lie(M)&&\underset{\downarrow A_0}\to{\overset{exp}\to{\to}}&&E(M) && \downarrow \underline{A}_{\p}\\ 
&\overset{\g j_{\{e\}}}\to{\searrow}&&&&\overset{\g i_{\{e\}}}\to{\searrow}&\\
&& \p^n&&\overset{expc_{\{e\}}}\to{\to}     &&\p^n \endmatrix\eqno{(9.6)}$$
\medskip
%\newpage
{\bf 10. Lattice associated to an Anderson t-motive, and a principal exact sequence.}
\medskip
Let us consider an exact sequence ([G96], 5.9.16):
$$0\to \n F_q[T]\to \n F_q((T^{-1}))\to T^{-1}\cdot\n F_q[[T^{-1}]]\to0$$
We denote its terms by $A, \ K, \ K/A$ respectively, and we define a topology on $K, \ K/A$ by the condition that $T^{-i}$ tends to 0 as $i$ tends to $+\infty$; $A$ is discrete in this topology. Further, we consider an exact sequence 
$$0\to \Hom_{\n F_q}^{cont}(K/A,\p)\to \Hom_{\n F_q}^{cont}(K,\p)\to \Hom_{\n F_q}^{cont}(A,\p)\to0$$
and we denote its terms by $Z_1, \ Z_2, \ Z_3$ respectively.
\medskip
$Z_1, \ Z_2, \ Z_3$ have a natural structure of $\p[T,\tau]$-modules, see [A86], below the Lemma 2.6.4, and [G96], below Definition 5.9.21 (the reader should write down the formulas for this structure himself, there is the only one possibility for it). 
\medskip
Further, we should distinguish $T$ as an abstract variable, $T_{Lie}$ as a linear operator on $Lie(M)$, and $T_E$ as an operator on $E$. 
\medskip
{\bf Definition 10.1.} Let $\p\{T\}$ be a subset of $\p[[T]]$ formed by power series (here and below $a_i\in\p$) 
$$\sum_{i=0}^\infty a_iT^i \hbox{ such that $a_i\to0$ as $i\to\infty$; equivalently, $v_\infty(a_i)\to+\infty$}. $$

$\p\{T\}$ is a ring. It is the ring of analytic functions on $\p$ that converge on the domain $v_\infty(x)\ge0$. $\p\{T\}$ has a natural structure of $\p[T,\tau]$-module: the multiplication by $T$ is simply a multiplication of a power series by $T$, and $$\tau\cdot \sum_{i=0}^\infty a_iT^i:=\sum_{i=0}^\infty a_i^qT^i.$$

{\bf 10.1a.} We have $Z_1=\p\{T\}$ as $\p[T,\tau]$-modules. Really, let 

$$\vf: K/A\to \p, \ \ \ \vf \in Z_1$$ be a continuous map. We denote $\vf(T^{-i})$ by $a_{i-1}$, the condition of continuity means $a_i\to0$. The image if $\vf$ in $\p\{T\}$ is $\underset{i=0}\to{\overset{\infty}\to{\sum}}a_iT^i$. It is checked immediately that this is an isomorphism of $\p[T,\tau]$-modules.
\medskip
Although we shall not need this information, we mention that 
\medskip
$Z_2$ is the set of power series $\sum_{i=-k}^\infty a_iT^i$ ($k\ge 0$ is any number), all other conditions are the same;
\medskip
$Z_3$ is the set of polynomials in $T^{-1}$, with the free term 0, i.e. $Z_3=\{ \sum_{i=-k}^{-1} a_iT^i\}$; condition $a_i\to0$ does not exist for this case, the $\p[T,\tau]$-module structure is the same, the product $T\cdot T^{-1}$ is 0.
\medskip
Now we consider two exact sequences 
$$0\to \Hom_{\p[T,\tau]}(M,Z_1) \to \Hom_{\p[T,\tau]}(M,Z_2) \to \Hom_{\p[T,\tau]}(M,Z_3) \eqno{(10.2)}$$
$$0\to \Hom_{\n F_q[T]}(K/A,E(M)) \to \Hom_{\n F_q[T]}(K,E(M)) \to \Hom_{\n F_q[T]}(A,E(M)) \eqno{(10.3)}$$
Let $L=L(M)$ be the kernel of $exp$, hence we have an exact sequence 
$$0\to L(M)\overset{\zeta}\to{\to} Lie(M)\ {\overset{exp}\to{\to}}\ E(M)\eqno{(10.4)}$$

{\bf Theorem 10.5} (Anderson): Sequences (10.2) -- (10.4) are canonically isomorphic. $L$ is a free $\n F_q[T]$-module of rank $\le r$. The rank of $L$ is $r$ iff $exp$ is surjective. 
\medskip
See Theorem 12.11.2 below for a more strong form of Theorem 10.5. 
\medskip
Anderson  t-motives $M$ such that $exp$ is surjective are called uniformizable. 
\medskip
Hence, not always $exp_M$ is surjective, this is an important phenomenon that does not occur for the case of Drinfeld modules. We do not know an analog of the formula (8.4) for the case $n>1$.
\medskip
$L$ is denoted by $H_1(M)$ as well, and the rank of $L$ is denoted by $h_1(M)$. 
\medskip
The injection $\ze$ is a map of $\n F_q[T]$-modules. If $N=0$ then the action of $T_{Lie}$ on $Lie(M)$ is the multiplication by $\th$, hence in this case $L$ is a lattice in $Lie(M)$ in the meaning of Section 4. 
\medskip
To cover the case $N\ne0$, we should extend a definition of lattice as follows. Let $V=\p^n$ be a space, $N$ a nilpotent operator on $V$ and $T:=\th\cdot Id +N$, where $Id$ is the identity operator. 
\medskip
{\bf Definition 10.5.1.} A $N$-lattice $L$ of rank $r$ in $V$ is a subset of $V$ such that there exist elements $l_1, \dots, l_r\in V$ satisfying the below conditions (10.5.3), (10.5.4) (they are analogs of (4.1.1), (4.1.2)) such that $L$ is the set of the sums $$P_1(T)(l_1)+...+P_r(T)(l_r)\eqno{(10.5.2)}$$ where $\forall \ i \ P_i(T)$ is a polynomial in $T$ of some degree $\la_i$: $P_i(T)=\sum_{j=0}^{\la_i}c_{ij}T^j$ where $c_{ij}\in \n F_q$. 
\medskip
{\bf Condition 10.5.3.} $\n F_q[[T]]$-envelope of $L$ (clearly the action of $\n F_q[[T]]$ on $V$ is defined) is isomorphic to $\n R_\infty^r$ (i.e. elements of a basis of $L$ over $\n F_q[T]$ are linearly independent over $\n F_q[[T]]$) and
\medskip
{\bf Condition 10.5.4.} $\p$-envelope of $L$ is $V$.
\medskip
{\bf Definition 10.5.5.} Two $N$-lattices $L_1, \ L_2$ are isomorphic if exists an isomorphism $\ga: V\to V$ commuting with $N$ such that $\ga(L_1)=L_2$. 
\medskip
We see that $L(M)$ is a $N$-lattice in $Lie(M)$. 
\medskip
We refer to the proof of Theorems 10.5 and 12.11.2 to [A86] or [G96], Theorem 5.9.14. 
\medskip
Here we give definitions of some objects and constructions that are used for the proof, and some explicit formulas (in coordinates) of isomorphisms between terms of (10.2) -- (10.4). 
\medskip
Since $A=\n F_q[T]$ we have $\Hom_{\n F_q[T]}(A,E(M))=E(M)$. Let us prove that there exists a canonical isomorphism $\delta: \Hom_{\n F_q[T]}(K,E(M)) \to Lie(M)$, and let us describe it explicitly. 
\medskip
Let $\vf\in\Hom_{\n F_q[T]}(K,E(M))$. We denote $x_i:=\vf(T^{-i})\in E(M)$. As $T^{-i}$ tend to 0 as $i$ tends to  $+\infty$, we have $x_i$ tend to 0 as $i$ tends to $+\infty$. Let us choose a basis $\{e\}$, it identifies $Lie(M),\ E(M)$ with $\p^n$ via $\g j_{\{e\}}$, $\g i_{\{e\}}$. Since $expc_{\{e\}}$ is a local isomorphism near 0, for sufficiently large $i$ the elements $exp^{-1}(x_i)\in Lie(M)$ are well-defined. Since $\vf$ is a $\n F_q[T]$-map, we have $T_E(x_i)=x_{i-1}$ and hence $T_{Lie}(exp^{-1}(x_i))=exp^{-1}(x_{i-1})$. This means that the element 
$$T^i_{Lie}(exp^{-1}(x_i))\in Lie(M)$$ does not depend on $i$, where $i$ is sufficiently large. We denote it by $\delta(\vf)$. 
\medskip
More exactly, we have: $\g j_e(\delta(\vf))=A_0^i(expc_{\{e\}}^{-1}(\g i_{\{e\}}(x_i)))$ where $A_0$ is from (5.9) (see also (9.2)). 
\medskip
The inverse map $\delta^{-1}: Lie(M) \to \Hom_{\n F_q[T]}(K,E(M)) $ is defined as follows. Let $x\in Lie(M)$. $\forall \ i\in \n Z$ we have $$[\delta^{-1}(x)](T^{i})=exp(T^i_{Lie}(x))$$
\medskip
The proof that (10.2) and (10.3) are isomorphic is based on a trivial fact of the set theory. Let $S_1, \ S_2, \ S_3$ be 3 sets. Then 
$$\Hom(S_1\times S_2, S_3)=\Hom(S_1, \Hom(S_2, S_3))=\Hom(S_2, \Hom(S_1, S_3)).$$
Namely, we consider an exact sequence 
$$0\to \Hom(M\times K/A,\p) \to \Hom(M\times K,\p)  \to \Hom(M\times A,\p)  \eqno{(10.6)}$$ where for $f\in M$, $c\in\p$, $k\in K/A$, or $K$, or $A$, any $\vf$ belonging to these Hom's must satisfy
$$ \vf(\tau f,k)=\vf(f,k)^q; \ \ \ \vf(T f,k)=\vf(f,Tk); \ \ \ \vf(c f,k)=c\cdot \vf(f,k)^q$$
(10.6) is isomorphic to both (10.2), (10.3). 
\medskip
We denote the isomorphism from $\Hom_{\p[T,\tau]}(M,Z_1)$ to $L(M)$ by $v$, and the isomorphism from $\Hom_{\p[T,\tau]}(M,Z_1)$ to $\Hom_{\n F_q[T]}(K/A,E(M))$ by $w$.  
\medskip
{\bf 10.7.} Let us describe explicitly $\zeta\circ v: \Hom_{\p[T,\tau]}(M,Z_1)\to Lie(M)$. Let $\vf\in \Hom_{\p[T,\tau]}(M,Z_1)$ and $f\in M$. We denote $\vf(f)\in Z_1$ by $\sum_{i=0}^\infty a_iT^i$ and we denote $[w(\vf)](T^{-i})\in E(M)$ by $x_i$. Because of isomorphisms of the first terms of (10.2), (10.6), (10.3), and because of (10.1), we get that the image of $\vf$ in $\Hom_{\n F_q[T]}(K/A,E(M))$ satisfies $x_i(f)=a_{i-1}$. 
\medskip
Since $\zeta\circ v(\vf)=T_{Lie}^i(exp^{-1}(x_i))$ for all sufficiently large $i$, we have
$$[exp(T_{Lie}^{-i}(\zeta\circ v(\vf)))](f)=a_{i-1}\eqno{(10.8)}$$
\medskip
Let now $x\in Lie(M)$ and $\{e\}$ as above. We denote $[exp(T_{Lie}^{-i}(x))](\{e\})$ by $a_{*i}=(a_{1i},\dots,a_{ni})^t$. 
\medskip
%\newpage
{\bf 10.8a. Finding of $\g j_{\{e\}}(x)$ in terms of $a_{*i}$.} 
\medskip
We want to get the below formula (10.10). We consider a commutative diagram (see (6.6) for $\underline{A}_{\p}$): 
$$\matrix Lie(M)&\overset{exp}\to{\to}&E(M)&&\p^n&\overset{expc_{\{e\}}}\to{\to}&\p^n\\ \\ T_{Lie}^i\downarrow&& \downarrow T_E^i&&A_0^i\downarrow&& \downarrow (\underline{A}_{\p})^i\\ \\ Lie(M)&\overset{exp}\to{\to}&E(M)&& \p^n&\overset{expc_{\{e\}}}\to{\to}&\p^n\
\endmatrix\eqno{(10.9)}$$ 
together with the arrows (not drawn) $\g j_{\{e\}}: Lie(M)\to \p^n$ (two arrows) and $\g i_{\{e\}}: E(M)\to \p^n$ (two arrows); this is (9.6) for $T^i$ instead of $T$. 
\medskip
We consider $T_{Lie}^{-i}(x)$ in the upper $Lie(M)$, and we apply the diagram chasing to it. We have $\g j_{\{e\}}(x)$ is its image in the down-left $\p^n$. Further, $[exp(T_{Lie}^{-i}(x))](\{e\})=\g i_{\{e\}}\circ exp\circ T_{Lie}^{-i}(x)=a_{*i}$. We get that $\g j_{\{e\}}(x)=A_0^i(expc_{\{e\}}^{-1}(a_{*i}))$. For $i\to+\infty$ we have $a_{*i}\to0$. Since $expc_{\{e\}}^{-1}$ is $id$ near 0, we get that (recall that $A_0$, $A_{\p}$ depend on $\{e\}$) $$\g j_{\{e\}}(x)=\underset{i\to+\infty}\to{\lim}A_0^i(a_{*i}).\eqno{(10.10)}$$

(To check: where this formula will be used?) To look at 11.12. 
\medskip
{\bf 10.11.} Finally, we give an important definition. Let $\la\in Lie(M)$ and $f\in M$. There exists a canonically defined element $\partial_\la(f)\in \p$, see [A86], 3.3.2. Its coordinate definition is the following. 
Let $f=\sum_{i=1}^n a_{i*}e_i\in M$, where $\forall  \ i$ we have $a_{i*}=\sum_j a_{ij}\tau^j$, $a_{ij}\in \p$, and $\g j_{\{e\}}(\la)=\{x\}=(x_1,\dots,x_n)^t$. Then $$\partial_\la(f):=\sum_{i=1}^n a_{i0}x_i\eqno{(10.12)}$$

The end of this section gives a proof that this expression for $\partial_x(f)$ does not depend on a choice of a basis $\{e\}$. 
\medskip
Let $G_a$ be the additive group. First, we describe the set $\Hom(E,G_a)$ in coordinates (i.e. let a basis $\{e\}$ be fixed). We begin with the formula $\Hom(G_a,G_a)=\underset{\mu=0}\to{\overset{\infty}\to{\bigcup}} \p^{\mu+1}: $ if $\g c=\{c_0, \dots, c_\mu\}\in \p^{\mu+1}$ and $x\in G_a=\p$ then $$\g c(x)=\sum_{i=0}^\mu c_ix^{q^i}$$

For a basis $\{e\}$ we have a coordinate map $\g k_{\{e\}}: \ \Hom(E,G_a)\to      
\underset{\mu=0}\to{\overset{\infty}\to{\bigcup}} \p^{n(\mu+1)}$. We denote an element of $\p^{n(\mu+1)}$ as $\{c_{ij}\}=\{c_{*0},\dots,c_{*\mu}\}\in\p^{n(\mu+1)}$ where $i=1,\dots,n, \ \ j=0,\dots,\mu$. 
The map $\g k_{\{e\}}$ is defined as follows. Let $\be\in E(M)$ and $\g i_{\{e\}}(\be)=(x_1, \dots,x_n)^t=\{x\}\in \p^n$. Let $\vf\in \Hom(E,G_a)$. We have $\g k_{\{e\}}(\vf)=\{c_{ij}\}\in\p^{n(\mu+1)}$ iff $\forall \ \be\in E(M)$ holds $$\vf(\be)=\sum_{j=0}^\mu\sum_{i=1}^n c_{ij}x_i^{q^j}.\eqno{(10.13)}$$
\medskip
The diagram of concordance of $\g k_{\{e\}}, \ \g k_{\{e'\}}$ is the following ($W$ is from 6.8): 
$$\matrix \Hom(E,G_a)&\overset{\g k_{\{e\}}}\to{\to}&\underset{\mu=0}\to{\overset{\infty}\to{\bigcup}} \p^{n(\mu+1)}   \\ \\ id\downarrow&& \downarrow W^{-1}\\ \\ \Hom(E,G_a)&\overset{\g k_{\{e'\}}}\to{\to}&\underset{\mu=0}\to{\overset{\infty}\to{\bigcup}} \p^{n(\mu+1)} \endmatrix\eqno{(10.14)}$$ 

where the action of $GL_n(\p\{\tau\})$ on $\underset{\mu=0}\to{\overset{\infty}\to{\bigcup}} \p^{n(\mu+1)}$ is the following. Let $\{c_{ij}\}=\{c_{*0},\dots,c_{*\mu}\}\in\p^{n(\mu+1)}$ as above and $W\in GL_n(\p\{\tau\})$ from 6.8. The action is from the right: 
$$ \{c_{ij}\}^W:=\{c_{*0}W_0; \ c_{*0}W_1+c_{*1}W_0^{(1)}; \ c_{*0}W_2+c_{*1}W_1^{(1)} +c_{*2}W_0^{(2)}; \dots \}$$

For any $\vf\in \Hom(E,G_a)$ and $\la\in Lie(M)$ we can define a number $\partial_\la(\vf)\in\p$ as follows. Let $\{e\}$ be as above, $\g k_{\{e\}}(\vf)=\{c_{ij}\}\in \underset{\mu=0}\to{\overset{\infty}\to{\bigcup}} \p^{n(\mu+1)}$, $\g j_{\{e\}}(\la)=\{x\}$. We let  $\partial_\la(\vf):=\sum_{i=1}^n c_{i0}x_i$. (9.1) and (10.14) show that it does not depend on $\{e\}$.  
\medskip
In order to prove that (10.12) does not depend on $\{e\}$, we consider the above $f=\sum_{i=1}^n a_{i*}e_i\in M$. It defines an element $\vf_f$ of $\Hom(E,G_a)$ such that $\g k_{\{e\}}(\vf_f)=\{a_{ij}\}$. We have that $\partial_\la(\vf)=\partial_\la(f)$, hence it does not depend on $\{e\}$. 
\medskip

\medskip
{\bf 11. $\th$-shift.} 
\medskip
Let $\sum_{i=0}^\infty a_iT^i\in \p[[T]]$. Let us consider its $\th$-shift, i.e. the result of substitution $T=N+\th$. We define $$S_k(i)=S_k(i)[\sum_{j=0}^\infty a_jT^j]:=\sum_{j=0}^k (-1)^{k-j}\binom{k}{j} a_{i+j}\th^{i+j}\in \p.$$

{\bf Example.} $S_0(i)=a_{i}\th^{i}$.
\medskip
$S_1(i)=a_{i+1}\th^{i+1}-a_{i}\th^{i}$.
\medskip
$S_2(i)=a_{i+2}\th^{i+2}-2a_{i+1}\th^{i+1}+a_{i}\th^{i}$.
\medskip
$S_3(i)=a_{i+3}\th^{i+3}-3a_{i+2}\th^{i+2}+3a_{i+1}\th^{i+1}-a_{i}\th^{i}$ etc. 
\medskip
Since $S_{k+1}(i)=S_k(i+1)-S_k(i)$, we have: 
\medskip
If for some $k_0$ we have $S_{k_0}(i)\to0$ as $i\to\infty$, then 
\medskip
First, $\forall \ k>k_0$ we have $S_k(i)\to0$ as $i\to\infty$;
\medskip
Second, there exists $\underset{i\to+\infty}\to{\lim}S_{k_0-1}(i)$.
\medskip
Let $F(T)=\sum_{i=0}^\infty a_iT^i\in \p[[T]]$ such that $\exists \ k_0$ such that $S_{k_0}(i)\to0$ as $i\to\infty$. We denote by $k$ the minimal number such that $S_k(i)\to0$. 
\medskip
{\bf Lemma 11.1.} Let $F(T)$, $k$ be as above. 
Then after the substitution $T=N+\theta$ we have: $F(T)$ becomes a well-defined series
\medskip
$G(N)=\sum_{i=-k}^\infty d_iN^i$, and (here $k>0$) 
$$d_{-k}=(-\th)^k \underset{i\to+\infty}\to{\hbox{ lim }}S_{k-1}(i).\eqno{(11.2)}$$
{\bf Proof.} Elementary calculations. $\square$
\medskip
In order to get a formula for $d_{-k+1}$ we apply (11.2) for $F(T)-d_{-k}(-\th+T)^{-k}\in\p[[T]]$, etc. for all $d_{-i}$.
\medskip
We need the following lemma. Let $\bar N$ be an operator on $\p^n$ satisfying $\bar N^\g m=0$ and $y\in \p^n$ any element (a vector column). Let us consider $\sum _{k=1}^\infty (\th+\bar N)^{-k}(y)T^{k-1}\in \p[[T]]^n$. 
\medskip
{\bf Lemma 11.3.} (Algebraic lemma.) The $\th$-shift of $\sum _{k=1}^\infty (\th+\bar N)^{-k}(y)T^{k-1}$ is given by 
$$\sum _{k=1}^\infty (\th+\bar N)^{-k}(y)T^{k-1}=-\sum _{l=1}^{\g m}\frac{\bar N^{l-1}(y)}{N^{l}}$$

{\bf Trivial example.} Let $n=1$, $y=1$, $\bar N=0$, hence $\g m=1$. The equality becomes $\sum _{k=1}^\infty \th^{-k}T^{k-1}=-\frac1{N}$ where $N=T-\th$. 
\medskip
{\bf Proof.} Immediate; this is (up to details) [A86], formula (3.3.5). We have 
$$(\th+\bar N)^{-k}(y)=\sum_{\al=0}^\infty (-1)^\al\binom{k+\al-1}{\al}\th^{-(k+\al)}\bar N^\al(y)$$
$$\sum _{k=1}^\infty (\th+\bar N)^{-k}(y)T^{k-1}=\sum _{k=1}^\infty\sum_{\al=0}^\infty (-1)^\al\binom{k+\al-1}{\al}\th^{-(k+\al)}\bar N^\al(y)T^{k-1}\eqno{(11.4)}$$
$$N^{-l}=(-1)^l\sum _{\be=0}^\infty \binom{l+\be-1}{\be} \th^{-(l+\be)}T^\be$$
$$-\sum _{l=1}^{\g m}\frac{\bar N^{l-1}(y)}{N^{l}}=-\sum _{l=1}^{\g m}\sum _{\be=0}^\infty (-1)^l\binom{l+\be-1}{\be} \th^{-(l+\be)}T^\be \bar N^{l-1}(y)\eqno{(11.5)}$$
Substituting $\be=k-1$, $l=\al+1$, we get that the terms of the right hand side of (11.5) become the terms of the right hand side of (11.4). $\square$

\medskip

Let us fix a basis $\{e_*\}$ of $M$ over $\p\{\tau\}$, and let $x\in Lie(M)$. We denote $\g i_{\{e\}}(exp(T_{Lie}^{-k}(x)))$ by $\mu_k$, it belongs to $\p^n$. Let us consider $\sum_{k=1}^\infty \mu_kT^{k-1}\in (\p[[T]])^n$.
\medskip
{\bf Lemma 11.6.} The $\th$-shift of $\sum_{k=1}^\infty \mu_kT^{k-1}$ is well-defined, it has the form $$\sum_{k=1}^\infty \mu_kT^{k-1}=\sum_{l= -\g m}^\infty d_lN^l\eqno{(11.7)}$$
and for $l\in[1, \dots, \g m]$ we have $d_{-l}=-\g j_{\{e\}}(N_{Lie}^{l-1}(x))$. 
\medskip
{\bf Proof.} $\mu_k=\g i_{\{e\}}(exp(T_{Lie}^{-k}(x)))=expc_{\{e\}}(\g j_{\{e\}}(T_{Lie}^{-k}(x)))=$
\medskip
$=expc_{\{e\}}(\th+N_{\{e\}})^{-k}(\g j_{\{e\}}(x)))$ (for $N_{\{e\}}$ see the line below (9.2) ).
\medskip
Let $C_l$ be from (9.4b). We have $$expc_{\{e\}}((\th+N_{\{e\}})^{-k}(\g j_{\{e\}}(x)))=(\th+N_{\{e\}})^{-k}(j_{\{e\}}(x))+$$ $$+\sum_{l=1}^\infty C_{l}\cdot [(\th+N_{\{e\}})^{-k}(j_{\{e\}}(x))]^{(l)}.\eqno{(11.8)}$$

Let $l\ge1$ be fixed. We consider the sum $$\sum_{k=1}^\infty C_{l} \cdot [(\th+N_{\{e\}})^{-k}(j_{\{e\}}(x))]^{(l)}T^{k-1}\in (\p[[T]])^n.$$ 
Its $\th$-shift has no terms of negative powers of $N$, because $\underset{i\to+\infty}\to{\lim}S_{0}(i)=0$. Really, while $k$ grows, the minimum of $v_\infty$'s (see Section 3 for the valuation $v_\infty$) of entries of $C_{l}\cdot [(\th+N_{\{e\}})^{-k}(j_{\{e\}}(x))]^{(l)}$ grows as $kq^l+const$, hence $$v_\infty (\th^{k-1}\cdot\hbox{ (entries of }C_{i} [(\th+N_{\{e\}})^{-k}(j_{\{e\}}(x))]^{(i)}\ )\to \ + \ \infty.$$
This means that only the term $(\th+N_{\{e\}})^{-k}(j_{\{e\}}(x))$ of (11.8) contributes to $d_l$ of (11.7) for $l<0$. According to Lemma 11.3, we have $d_{-l}=-N_{\{e\}}^{l-1}(\g j_{\{e\}}(x))$. We have $N_{\{e\}}^{l-1}(\g j_{\{e\}}(x))=\g j_{\{e\}}(N_{Lie}^{l-1}(x))$, hence the lemma. $\square$
\medskip
Let, like in 10.7,  $\vf\in \Hom_{\p[T,\tau]}(M,Z_1)=L(M)$ and $f\in M$. We denote $\vf(f)\in Z_1$ by $\sum_{i=0}^\infty a_iT^i$. 
\medskip
{\bf Theorem 11.9.} ([A86], Th. 3.3.2). The $\th$-shift of $\sum_{i=0}^\infty a_iT^i$ is well-defined. We denote it by 
$$\sum_{i=0}^\infty a_iT^i=\frac{z_k}{N^k}+...+\frac{z_1}{N}+\sum_{j=0}^\infty y_jN^j$$

We have $k=\g m$ and for $i\ge0$ we have $$\partial_{N^i(\zeta(\vf))}(f)=-z_{i+1}\eqno{(11.10)}$$
\medskip
{\bf Remark.} This formula will be used in the proof of the equality (12.4.13), see below.
\medskip
{\bf Proof.} Let $f=\{e_*\}$ (a basis). We have $\partial_{\zeta(\vf)}(e_*)=\g j_{\{e\}}(\ze(\vf))$ (vector column). Let us consider the column series 
$$\vf(e_*)=\sum_{i=0}^\infty a_{i*}T^i$$
(here $a_{i*}$ are column vectors). According (10.8), we have 
$$exp(T^{-(i+1)}(\ze(\vf)))(\{e_*\})=a_{i*}\eqno{(11.11)}$$

For $x\in E$ we have $x(\{e_*\})=\g i_{\{e_*\}}(x)$, hence (11.11) means 
$$a_{i*}=\g i_{\{e_*\}}(exp(T^{-(i+1)}(\ze(\vf))))$$
The theorem for $f=$ one of $e_i$ follows from Lemma 11.6 applied for $\ze(\vf)$. For this case we have $k=\g m$. 
\medskip
By linearity, to finish the proof, it is sufficient to prove (11.10) for $f=\tau^ke_j$. We have $\partial_{\ze(\vf)}(\tau^k(e_j))=0$. Hence, we must prove that the $\th$-shift of $\vf(\tau^k e_j)$ has no terms of negative powers of $N$, i.e. $\underset{i\to+\infty}\to{\lim}S_{0}(i)[\vf(\tau^k e_j)]=0$. 

Let $\vf(e_j)=\sum_{\be=0}^\infty a_\be T^\be$. Since $\vf$ commute with $\tau$, we have 
$$\vf(\tau^k(e_j))=\sum_{\be=0}^\infty a_\be^{q^k}T^\be.$$
We have $S_\g m(i)[\sum_{\be=0}^\infty a_\be T^\be]\to0$. This implies that $v_\infty(\th^\be a_\be)$ are bounded from the down. Really if $v_\infty(\th^\be a_\be)$ are not bounded from the down, then there exists a sequence $\be_1, \ \be_2, \dots$ having the properties:
\medskip
(a) $\forall \ \ga=1,2,\dots$, $\forall \ \be<\be_\ga$ we have $v_\infty(\th^{\be} a_{\be})>v_\infty(\th^{\be_\ga} a_{\be_\ga})$;
\medskip
(b)  $v_\infty(\th^{\be_\ga} a_{\be_\ga})$ tends to $-\infty$ as $\ga\to\infty$. 
\medskip
Property (a) implies that $v_\infty(S_\g m(\be_\ga-\g m)[\sum_{\be=0}^\infty a_\be T^\be])=v_\infty(\th^{\be_\ga} a_{\be_\ga})$. 
\medskip
This implies that $v_\infty(S_\g m(\be_\ga-\g m)[\sum_{\be=0}^\infty a_\be T^\be])$ tends to $-\infty$ --- a contradiction to the condition $S_\g m(i)[\sum_{\be=0}^\infty a_\be T^\be]\to0$. 
\medskip
This implies that $\underset{i\to+\infty}\to{\lim}S_0(i)[\sum_{\be=0}^\infty a_\be^{q^k}T^\be]=0$. Really, 
\medskip
$S_0(i)[\sum_{\be=0}^\infty a_\be^{q^k}T^\be]=\th^ia_i\cdot a_i^{q^k-1}$, 
\medskip
$v_\infty (\th^ia_i)$ are bounded from the down and $a_i^{q^k-1}\to0$. 
\medskip
Hence, $\underset{i\to+\infty}\to{\lim}S_{0}(i)[\vf(\tau^k e_j)]=0$. $\square$
\medskip
{\bf Remark 11.12.} The above results can be reformulated as follows. Let us fix a basis $\{e_*\}$ of $M$ over $\p\{\tau\}$, and let $l\in L(M)$. From one side, we consider $l$ as an element of $\Hom_{\p[T]}(M,Z_1)^\tau$. We denote $l(\{e_*\})$ by $A=\sum A_iT^i\in Z_1^n$, $A_i\in \p^n$. From another side, $l\in Lie(M)$, and we have $\g j_\{e_*\}(l)\in \p^n$; we denote it by $\{b_*\}$.
\medskip
What is a relation between $A$ and $\{b_*\}$?
\medskip
{\bf Lemma 11.12.1.} $\{b_*\}$ is the coefficient at $N^{-1}$ of the $\th$-shift of $A$.
\medskip
{\bf Proof} follows immediately from the above formulas. $\square$
\medskip
To check: 10.8a: is this the same or not? Is it true if $N\ne0$?

\medskip
%\newpage
{\bf 12. Dualities and pairings.}
\medskip
{\bf 12.1. Dual of a t-motive.} Let $M$ be an Anderson t-motive having $N=0$. Its dual $M'$ is defined by (5.19.2), (5.19.1). Let us recall here this definition: 
$$M':= \Hom_{\p[T]}(M,\g C)$$
where the action of $\tau$ on $\Hom_{\p[T]}(M,\g C)$ is defined by the formula 
$$\tau(\vf)(m_1):=\tau(\vf(\tau^{-1}(m_1)))$$
where $\vf\in \Hom_{\p[T]}(M_1,M_2)$, $m_1\in M_1$. 
\medskip
For the case of Drinfeld modules the definition of dual was given in [Tag95], for any t-motive in [GL07]. A version of the definition of dual for more general objects is given in [HJ]. Hartl and Juschka consider $M^*_{HJ}:=\Hom_{\p[T]}(M,${\bf1}) where {\bf1} is a 1-dimensional $\p[T]$-module with the trivial $\tau$-action, hence $M'=M^*_{HJ}\underset{\p[T]}\to{\otimes}\g C$. 
\medskip
Let $M$ be a t-motive with $N=0$ having the dual $M'$. Let $\{f\}=\{f_1, \dots,
f_r\}^t$ be a basis of $M$ over $\p[T]$ (see 5.11), and let $\{e\}=e_1$ the only element of a basis of $\g C$ over $\p[T]$. They define the dual basis $\{f'\}=\{f'_1, \dots,
f'_r\}^t$ of $M'$ over $\p[T]$ by the formulas 
$$f'_i(f_j)=\delta^i_je$$
$\{f'\}$ is defined up to multiplication by an element of $\n F_q^*$; the same is true for some other below objects.
\medskip
Properties of $M'$ are the following:
\medskip
{\bf 12.1.1.} The matrix $Q'$ of $M'$ (see 5.11) in $\{f'\}$ is $(T-\th)(Q^t)^{-1}$ where $Q$ is the matrix for $M$ over $\{f\}$.
\medskip
{\bf 12.1.2.} Dimension of $M'$ is $r-n$ (the rank is the same, i.e. $r$).
\medskip
{\bf Example 12.2: the dual of a Drinfeld module.} Let a Drinfeld module $M$ be given by (5.6), $a_0=\th, \ a_r=1$: 

$$Te=(\th+a_1\tau+a_2\tau^2+\dots +a_{r-1}\tau^{r-1}+\tau^r)e$$
The basis of $M$ over $\p[T]$ is $f_*=(f_1,\dots,f_r)^t:=(e,\tau e, \tau^2e,\dots, \tau^{r-1}e)^t$. The matrix $Q$ in $f_*$ is given in (5.12.1). The matrix of $M'$ in the dual basis $f'_*$ is 
$$(T-\th)Q^{t-1}=\left(\matrix a_1&T-\th&0&0&\dots&0&0
\\ a_2&0&T-\th&0&\dots&0&0
\\  \dots&\dots&\dots&\dots&\dots&\dots&\dots
\\ a_{r-1}&0&0&0&\dots&0&T-\th
\\ 1 & 0&0&0&\dots&0&0\endmatrix\right)\eqno{(12.2.1)}$$ We consider a dual basis $e'_*:=(e'_1,\dots, e'_{r-1})^t:=(f'_2,\dots, f'_{r})^t$. We have $f'_1=\tau(f'_r)$ and the action of $T$ on $e'_*$ is given by the formula $$Te'_*=\th e'_*+A_1\tau e'_* + A_2\tau^2 e'_*\eqno{(12.2.2)}$$ where 
$A_1=\left(\matrix  0&0&\dots&0&0&-a_1 \\ 1&0&\dots&0&0&-a_2 \\ 0&1&\dots&0&0&-a_3\\ \dots&\dots&\dots&\dots&\dots&\dots\\ 0&0&\dots&1&0&-a_{r-2} \\ 0&0&\dots&0&1&-a_{r-1} \endmatrix \right), \ \ \ A_2=\left(\matrix 0&0&\dots&0&1 \\ 0&0&\dots&0&0 \\   0&0&\dots&0&0 \\  \dots&\dots&\dots&\dots&\dots \\  0&0&\dots&0&0 \\  \endmatrix \right)$,

See [Tag95], Section 5, p. 578. Similar explicit formulas for duals of some other t-motives (called standard t-motives) are given in [GL07], Section 11.
\medskip
{\bf 12.3. Definitions, isomorphisms and pairings for $H_1(M)$, $H_1(M')$, $H^1(M)$, $H^1(M')$.}
\medskip
{\bf 12.3.1. Case of $H_1(M)$.} We can interpret $L(M)=H_1(M)$ as follows. Let 
$$M_{\{T\}}:=\Hom_{\p[T]}(M,\p\{T\}).\eqno{(12.3.1.1)}$$ We have: $\tau$ acts on $M_{\{T\}}$ by the standard action of operator on Hom group, i.e. $(\tau(\vf))(m):=\tau(\vf(\tau^{-1}(m)))$ (we neglect here that $\tau^{-1}(m)$ maybe does not exist --- this is a model example). We have (Theorem 10.5):
$$L(M)=H_1(M)=\Hom_{\p[T,\tau]}(M,\p\{T\})={M_{\{T\}}}^\tau=(\Hom_{\p[T]}(M,\p\{T\}))^\tau$$

The explicit description of $H_1(M)$ is the following. Recall that $f_*:=(f_1, \dots, f_r)^t$ is a basis of $M$ over $\p[T]$ and $Q\in M_{r\times r}(\p[T])$ the matrix of multiplication by $\tau$ in $f_*$, i.e. $$\tau f_* = Q f_*.$$ Let $$l\in H_1(M)=\Hom_{\p[T]}(M,\p\{T\})^\tau.$$ We denote $$X_*=X_*(l):=l(f_*)\in M_{r\times 1}(\p\{T\}).\eqno{(12.3.1.2)}$$ Condition that $l$ is $\tau$-invariant is equivalent to the condition $$QX_*=X_*^{(1)}.\eqno{(12.3.1.3)}$$
\medskip
Let now $M$ be uniformizable. Let $l_1,\dots, l_r$ be a basis of $H_1(M)$, and $X_1, \dots, X_r$ are the above $r\times 1$ matrices corresponding to $l_1,\dots, l_r$. We consider their columnwise union, i.e. a matrix $\Psi\in M_{r\times r}(\p\{T\})$ such that for $i=1,\dots, r$ the $i$-th column of $\Psi$ is $X_i$. By definition, we have $$\Psi_{ij}=l_j(f_i)\eqno{(12.3.1.4)}$$ where $\Psi_{ij}$ is the $(i,j)$-th entry of $\Psi$. 
\medskip
$\Psi$ is called the scattering matrix ([A86], Section 3.1) with respect to the bases $l_*, \ f_*$. It belongs to $ GL_r(\p\{T\})$, because $X_1, \dots, X_r$ are linearly independent. (12.3.1.3) becomes $$Q\Psi=\Psi^{(1)}\eqno{(12.3.1.5)}$$
\medskip
{\bf 12.3.2. Case of $H^1(M)$.} 
\medskip
Let us consider the tensor product instead of Hom in (12.3.1.1), i.e.
we define
$$ M\{T\}:=M\otimes_{\p[T]}\p\{T\}\eqno{(12.3.2.1)}$$
We have: $\tau$ acts on $M\{T\}$ by the standard formula of the action of an operator on tensor products (i.e. $\tau(a\otimes b)=\tau(a)\otimes \tau(b)$ ), and we define
$$H^1(M):=M\{T\}^\tau\eqno{(12.3.2.2)}$$

$H^1(M)$ is also free $\n F_q[T]$-module, its rank is denoted by $h^1(M)$. Like for the case of $H_1(M)$, we have $h^1(M)\le r$ (see, for example, Theorem 16.11.2). 
\medskip
Let $c\in H^1(M)=(M\underset{\p[T]}\to{\otimes}\p\{T\})^\tau$, $c=\sum_{i=1}^r f_i\otimes y_i$. We denote $$Y_*=Y_*(c):=(y_1, \dots, y_r)\in M_{1\times r}(\p\{T\})\eqno{(12.3.2.3)}$$ a row matrix. Condition that $c$ is $\tau$-invariant is equivalent to the condition $$Y_*^{(1)}Q=Y_*.\eqno{(12.3.2.4)}$$

Let us consider $M$ such that $h^1(M)=r$. Let $c_1,\dots, c_r$ be a basis of $H^1(M)$ over $\n F_q[T]$, and $Y_i=Y_*(c_i)$ the corresponding row matrices. We consider their rowwise union, i.e. a matrix $\Phi\in M_{r\times r}(\p\{T\})$ such that for $i=1,\dots, r$ the $i$-th row of $\Phi$ is $Y_i$. It satisfies $$\Phi^{(1)}Q=\Phi\eqno{(12.3.2.5)}$$
\medskip
{\bf Remark.} Both (12.3.1.5), (12.3.2.5) are particular cases of a system of affine equations, see Section 16. 
\medskip
{\bf Example 12.4. Explicit calculations for $\g C$.} 
\medskip
{\bf Case of $H^1(\g C)$.} Recall that $\g C$ is the Carlitz module. Let as above $e$ be the only element of a basis of $\g C$ over $\p[T]$, it satisfies $\tau e=(T-\th)e$. Elements of $\g C\{T\}$ have the form
$$(d_0+d_1T+d_2T^2+...)e\eqno{(12.4.1)}$$ where $d_i\in \p$, $d_i\to0$. (12.3.2.4) for it is an affine equation (see Section 16)
$$(d_0^q+d_1^qT+d_2^qT^2+...)(T-\th)=d_0+d_1T+d_2T^2+...\eqno{(12.4.2)}$$ Multiplying we get a system of constituent affine equations
$$\th d_0^q+d_0=0$$
$$\th d_1^q+d_1-d_0^q=0\eqno{(12.4.3)}$$
$$\th d_2^q+d_2-d_1^q=0$$
etc. Let us solve them consecutively choosing (for all equations) the root corresponding to the leftmost segment of the Newton polygon, i.e. the root with the maximal $v_\infty$. We get that for these roots $d_0, \ d_1, \dots$ we have $$v_\infty(d_i)=\frac{q^i}{q-1}\to + \infty.$$ 
\medskip
{\bf Definition 12.4.4.} This element $\sum_{i=0}^\infty
d_iT^i\in H^1(\g C)$ is denoted by $\Xi$, see [G96], p. 171, (*) and p. 172, line 1. It is defined up to
multiplication by $\n F_q^*$. It satisfies
$$\Xi=(T-\theta)\Xi^{(1)}, \ \ \ \lim_{i\to\infty}d_i=0, \ \ \ |d_0|>|d_i| \ \ \forall
i>0\eqno{(12.4.5)}$$ 
\medskip
{\bf Case of $H_1(\g C)$.} Let $X=x_0+x_1T+x_2T^2+...$ be a solution of (12.3.1.3) for $\g C$, i.e. it satisfies  $$(T-\theta)X=X^{(1)}.\eqno{(12.4.6)}$$ The system (12.4.3) for it has the form

$$x_0^q+\th x_0=0$$
$$x_1^q+\th x_1-x_0=0\eqno{(12.4.7)}$$
$$x_2^q+\th x_2-x_1=0$$
etc. As above, we choose the simplest solution of (12.4.6) choosing $x_i$ corresponding to the leftmost segment of the Newton polygon. We have for them $v_\infty(x_i)=i-\frac1{q-1}$. We denote this solution by $X_0$. It is defined up to
multiplication by $\n F_q^*$. 
\medskip
{\bf Remark 12.4.8.} We see that the growth of $v_\infty(d_i)$, $v_\infty(x_i)$ is exponential for $\Xi$ and linear for $X_0$. 
\medskip
{\bf Lemma 12.4.9.} We have $X_0=\Xi^{-1}$. 
\medskip
{\bf Proof.} Comparing (12.4.5) and (12.4.6) we get that $(X_0\Xi)^{(1)}=X_0\Xi$, hence $X_0\Xi\in \n F_q[[T]]$. We have $X_0\Xi=d_0x_0+(d_1x_0+d_0x_1)T+...$. Since $v_\infty(d_i)=\frac{q^i}{q-1}$ and $v_\infty(x_i)=i-\frac1{q-1}$ we get that $v_\infty(d_0x_0)=0$, $v_\infty$ of other coefficients of $X_0\Xi$ is $> 0$. Hence, $X_0\Xi=1$. $\square$
\medskip
Let us consider the lattice of the Carlitz module. Its exponential map $\exp_\g C: \p\to \p$ is given by the formula $$\exp_\g C(z)=z+c_1z^q+c_2z^{q^2}+...$$ where $c_i$ are from (8.2.1). We have $L(\g C)=\Ker \exp_\g C$. It is a 1-dimensional $\n F_q[\th]$-submodule of $\p$. We denote its generator by $\pi$, it is defined up to
multiplication by $\n F_q^*$. It is an analog of $2\pi i\in \n C$ (which is defined up to $\pm1$). 
\medskip
The set of $X$ satisfying (12.4.6) is in canonical (up to
multiplication by $\n F_q^*$) isomorphism $\ze$ with $L(\g C)$. 
\medskip
{\bf Lemma 12.4.10.} $X_0$ corresponds to $\pi$ with respect to $\ze$. 
\medskip
{\bf Proof.} The word "corresponds" should be understood in the meaning of Remark 11.12. We consider $X_0$ as an element of $\Hom_{\p[T]}(\g C,Z_1)^\tau$, where the basis $\{e_*\}=e_1$ of $\g C$ over $\p\{\tau\}$ is fixed. We need to show that $z_1$ (see below) is $-\pi$. 
\medskip

It is possible to show that $X_0$ is a generator of the set of solutions of (12.4.6) over $\n F_q[T]$. But it is simpler to use Theorem 11.9 (or Remark 11.12). We have (notations of Theorem 11.9) $\vf=X_0$, $f=e$ the only element of the basis of $\g C$ over $\p\{\tau\}$, $k=1$, $N=T-\th$. Let $$\frac{z_1}{N}+\sum_{j=0}^\infty y_jN^j\eqno{(12.4.11)}$$ be the $\th$-shift of $X_0$. (11.10) for $i=0$ is $$\partial_{\ze(X_0)}(e)=-z_1.$$
According the definition of $\partial$, we have $\partial_{\ze(X_0)}(e)=\ze(X_0)$. 

Further, (11.2) gives us $z_1=-\th \underset{i\to+\infty}\to{\hbox{ lim }}x_i\th^i$. Since $v_\infty(x_i)=i-\frac1{q-1}$ we get that $v_\infty(z_1)=-\frac{q}{q-1}$. Considering the Newton polygon of $\exp_\g C$ we have that $v_\infty(\pi)=-\frac{q}{q-1}$, hence the lemma. $\square$

Let $c:=(-\th)^{\frac{-q}{q-1}}$ (defined up to
multiplication by $\n F_q^*$). We define (see [AT90]) 

$$\oo=c(1-\th^{-q}T)(1-\th^{-q^2}T)(1-\th^{-q^3}T)...$$ 

{\bf Lemma 12.4.12.} $\oo^{(-1)}=\Xi$.
\medskip
{\bf Proof.} $\oo^{(-1)}$ satisfies $\oo^{(-1)}=(T-\th)\oo$. Let $\oo^{(-1)}=k_0+k_1T+k_2T^2+...$. We have $$v_\infty(k_i)=v_\infty(c^{1/q}\th^{-1-q-...-q^{i-1}})=\frac{q^i}{q-1}.$$ This means that $\oo^{(-1)}=\Xi$. $\square$
\medskip
Finally, we have (recall that $T-\th=N$) $$(T-\th)X_0=X_0^{(1)}=[\Xi^{(1)}]^{-1}=\oo^{-1}; \ \ \ (T-\th)X_0=z_1+\sum_{j=0}^\infty y_jN^{j+1}.$$ Substituting $T\mapsto \th$ we get (because $z_1=\pi$ mod $\n F_q^*$):
$$\pi=\oo^{-1}(\th)=(-\th)^{\frac{q}{q-1}}(1-\frac{\th}{\th^q})^{-1}(1-\frac{\th}{\th^{q^2}})^{-1}(1-\frac{\th}{\th^{q^3}})^{-1}\dots\eqno{(12.4.13)}$$ (formula of Carlitz?)
\medskip
{\bf 12.5.} For $M$ of rank $r$ we have 
\medskip
{\bf Theorem 12.5.1.} $h^1(M)=r \iff h_1(M)=r$.
\medskip
For a relation with a theorem of Anderson ([A86], Theorem 4; [G96], Theorem 5.9.14) see below, Theorem 12.11.2. Also, the proof gives us a relation between $\Psi$ and $\Phi$.
\medskip
{\bf Proof.} See, for example, [GL21], Remark 1.7.5. Let $h^1(M)=r$. The lines of $\Phi$ are linearly independent over $\n F_q[T]$. According the Lemma 1.5.4.1 of [GL21], they are linearly independent over $\p\{T\}$. This means that $\det \Phi\ne0$, $\det \Phi\in \p\{T\}$. 
\medskip
Taking the determinant of (12.3.2.5) and taking into consideration Lemma 5.11.2 we get $$\det\Phi^{(1)}\cdot c \cdot (T-\th)^n=\det\Phi.$$ 
Since all solutions of the equation  
$$Z=(T-\th)Z^{(1)}.$$
where $Z\in \p\{T\}$ is an unknown have the form $U\cdot\Xi$ where $U\in \n F_q[T]$ we get that $$\det \Phi\in\bar \n F_p\cdot \n F_q[T]\cdot\Xi^n.$$
\medskip
Now we define $\Psi_0:=\Phi^{-1}$ (a preliminary value of $\Psi(M)$). Since $\Xi^{-1} \in \p\{T\}$, we get: $$\Psi_0\in \n F_q(T)\cdot GL_r(\p\{T\}).$$ It satisfies (12.3.1.5). After multiplication $\Psi_0$ by a non-zero element of $\n F_q[T]$ we get $\Psi\in GL_r(\p\{T\})$, it continue to satisfy (12.3.1.5). This means that the columns of $\Psi$ form $r$ column vectors satisfying (12.3.1.3), hence this $\Psi$ is really $\Psi(M)$, i.e. $h_1(M)=r$.
\medskip
These considerations are convertible, i.e. $h_1(M)=r\implies h^1(M)=r$. $\square$
\medskip
{\bf 12.6. Pairings and isomorphisms.} Here we give definitions and explicit formulas for pairings and isomorphisms between $\n F_q[T]$-modules $H_1(M)$, $H^1(M)$, $H_1(M')$, $H^1(M')$. 
\medskip
For pairings, we have a diagram $$ \matrix  H_1(M)& \overset{\rho}\to{\longleftrightarrow} & H_1(M') \\ \\ \pi\updownarrow &&\updownarrow \\ \\ H^1(M) & \longleftrightarrow & H^1(M')\endmatrix\eqno{(12.6.1)}$$ where $\longleftrightarrow$ indicate existence of the canonical pairing with values in $\n F_q[T]$. 
\medskip
Further, there are canonical isomorphisms $$H_1(M)= H^1(M'), \ \ H_1(M')=H^1(M).$$

{\bf 12.6.2. Definition of the pairing $\pi$. } We have 
$$\pi: H_1(M)\underset{\n F_q[T]}\to{\otimes}H^1(M)\to\n F_q[T]$$

Let $y=\sum m_i\otimes Z_i\in H^1(M)$, where $m_i\in M$, $Z_i\in \p\{T\}$ and $\vf: M\to\p\{T\}$ belongs to $H_1(M)$.  By definition, $$\pi(y\otimes\vf)=\sum \vf(m_i)Z_i.$$ 

In coordinates, $\pi$ is given by the below formula (which one?). It is checked immediately that $\pi(y\otimes\vf)\in \n F_q[T]$, because $(YX)^{(1)}=YQ^{-1}QX=YX$, i.e. $YX\in \n F_q[T]$.
\medskip
{\bf 12.6.3. Definition of the isomorphism $H_1(M)\to H^1(M')$.} 
\medskip
It is given in coordinates. We choose a basis $f_*=(f_1,\dots,f_r)$ of $M$ over $\p[T]$. Hence, we have the dual basis $f'_*=(f'_1, \dots, f'_r)$ of $M'$ over $\p[T]$. We identify $H_1(M), \ H^1(M')$ with $X_*, \ Y_*$ from (12.3.1.2) for $M$, (12.3.2.3) for $M'$. 
\medskip
Since $Q(M')=(T-\th)Q(M)^{t-1}$ we have 
$$X_*\hbox { satisfies (12.3.1.3) $\iff \ \Xi X_*^t$ satisfies (12.3.2.4) for } M'$$

This means that the map $X_*\mapsto \Xi X_*^t$ (see (12.7.7) below) gives us the isomorphism $H_1(M)\to H^1(M')$.
It is easy to check that it does not depend on the choice of a basis $f_*=(f_1,\dots,f_r)$, i.e. it is canonical (up to multiplication by $\n F_q^*$). 
\medskip
Other pairings and isomorphisms come from these ones. 
\medskip
{\bf Remark 12.6.4.} This isomorphism exists for all $M$, not necessarily uniformizable. Particularly, $\forall \ M$ we have 
$$h_1(M)=h^1(M') \hbox{ and hence } h^1(M)=h_1(M')$$
From another side, unlike existence of pairings of (12.6.1), there exist $M$ such that $h^1(M)\ne h_1(M)$. Namely, for a t-motive $M(A)$ of Example 5.13 where $n=2$, $q=2$ and $A=\left(\matrix \th & \th^6\\ \th^{-2}&0 \endmatrix \right)$ we have $h^1(M(A))=0$, $h_1(M(A))=1$ ([GL21], Section 4).
\medskip
%\newpage
{\bf 12.7. List of formulas for pairings and isomorphisms.}
\medskip
For the reader's convenience, we give here a list of formulas for pairings and isomorphisms between $H_1(M)$, $H^1(M)$, $H_1(M')$, $H^1(M')$. Recall that they are valid for not necessarily uniformizable $M$, $M'$. 
\medskip

The matrices $Q, \ Q'\in \p[T]$ are matrices of multiplication by $\tau$ in $f_*, \ f'_*$:
$$\tau f_*=Qf_*, \ \ \ \tau f'_*=Q'f'_*$$
The below formulas are valid for any, i.e. non necessarily uniformizable, $M$ having dual (here for any matrix $A$ we denote $A^{t-1}:=(A^t)^{-1}$ \ ): $$Q'=(T-\th)Q^{t-1}, \ \ \ \ Q=(T-\th){Q'}^{t-1}$$
For $H_1(M)$, respectively $H_1(M')$ we have (see 12.3.1.3) $$QX=X^{(1)}, \hbox{ respectively } Q'X'={X'}^{(1)} \eqno{(12.7.1)}$$
For $H^1(M)$, respectively $H^1(M')$ we have (see 12.3.2.4) $$Y^{(1)}Q=Y, \hbox{ respectively } {Y'}^{(1)}Q'=Y'\eqno{(12.7.2)}$$

Pairing $\pi$ between $H_1(M)$ and $H^1(M)$: $$<Y, X>=YX\eqno{(12.7.3)}$$
Pairing $\rho$ between $H_1(M)$ and $H_1(M')$, where $M$ is considered as the initial t-motive: $$<X, X'>=\Xi \ {X'}^tX\eqno{(12.7.4)}$$
If we consider $M'$ as the initial t-motive, the result is the same. 
\medskip
Pairing between $H^1(M')$ and $H_1(M')$: $$<Y', X'>=Y'X'\eqno{(12.7.5)}$$
Pairing between $H^1(M)$ and $H^1(M')$, where $M$ is considered as the initial t-motive: $$<Y, Y'>=\Xi^{-1}Y{Y'}^t\eqno{(12.7.6)}$$
If we consider $M'$ as the initial t-motive, the result is the same. 
\medskip
Isomorphism $H_1(M)\to H^1(M')$: $$X\mapsto \Xi \ X^t\eqno{(12.7.7)}$$
Isomorphism $H^1(M)\to H_1(M')$: $$Y\mapsto \Xi^{-1} Y^t\eqno{(12.7.8)}$$

\newpage
{\bf 12.8. Connections between $\Psi(M), \ \Phi(M), \ \Psi(M'), \ \Phi(M')$.}
\medskip
Let us consider the case of uniformizable $M$, $M'$. In this case the pairing $\rho$ from (12.6.1) is perfect (to give a proof), hence for a basis $l_1,\dots, l_r$ of $H_1(M)$ over $\n F_q[T]$ there exists its dual $l'_1,\dots, l'_r$ of $H_1(M')$ over $\n F_q[T]$. 
\medskip
Matrices $\Psi(M), \ \Phi(M), \ \Psi(M'), \ \Phi(M')$ are defined in bases $f_1,\dots, f_r$, $f'_1,\dots, f'_r$, $l_1,\dots, l_r$, $l'_1,\dots, l'_r$. Their relations are: 

$$\Psi(M')=\Xi^{-1}(\Psi(M)^{t-1}\eqno{(12.8.1)}$$
$$\Phi(M)=\Psi(M)^{-1}\eqno{(12.8.2)}$$

Further, the matrix $\Psi^{-1}$ satisfies 
$$\Psi^{-1(1)}Q=\Psi^{-1}\eqno{(12.8.3)}$$
{\bf 12.8.4.} The lines of $\Psi^{-1}$ correspond to $X$-elements of a basis of $H^1(M)$. 
\medskip
(12.7.8) and (12.8.4) imply $$\Psi'=\Xi^{-1}\Psi^{t-1}.\eqno{(12.8.5)}$$
\medskip
{\bf 12.9. Duality of lattices.} 
\medskip
The simplest definition of the dual lattice is the following. Let
$L \subset \p^n$ be a lattice of rank $r$ and $S$ its Siegel matrix in a basis. By definition, its dual lattice $L'$ is a lattice whose Siegel matrix is $S^t$.
\medskip
If $S$ is a Siegel matrix of a lattice, not necessarily that $S^t$ is a Siegel matrix of a lattice (see Section 13), i.e. not all lattices have dual. (Most likely we should modify Definition 4.1 in order to exclude this phenomenon). 
\medskip
This means that $L'$ is a lattice of rank $r$ in $\p^{r-n}$. 
It is easy to show that this definition is well-defined (i.e. does not depend on a choice of a Siegel matrix). 
\medskip
Let us give an invariant definition (see, for example, [GL07], Definition 2.3; it is a simple version of (14.1) for the case $N=0$). 

Let $L_T$ be an abstract free module over $\n F_q[T]$ of dimension $r$ and $L'_T:=\Hom _{\n F_q[T]}(L_T,\n F_q[T])$ its dual. 
We have:
$ L_T\underset{\n F_q[T]}\to{\otimes}\p$ is isomorphic to $\p^r$.

Let now $L\overset{\ze}\to{\hookrightarrow} Lie(M)=\p^n$ be a lattice. $\ze$ gives us an exact sequence $$0\to \Ker \ze\to L_T\underset{\n F_q[T]}\to{\otimes}\p\overset{\ze}\to{\to} \p^n\to0\eqno{(12.9.1)}$$
where $\ze$ gives us an isomorphism from $L_T\otimes 1$ to $L\hookrightarrow Lie(M)=\p^n$. 
\medskip
We have $\dim_{\p}\Ker \ze=r-n$. Let us consider its dual $(\Ker \ze)'\hookrightarrow L'\underset{\n F_q[T]}\to{\otimes}\p$ with respect to the natural duality between $ L_T\underset{\n F_q[T]}\to{\otimes}\p$ and $ L'_T\underset{\n F_q[T]}\to{\otimes}\p$.
We have an exact sequence $$0\to (\Ker \ze)'\to L'\underset{\n F_q[T]}\to{\otimes}\p\to \p^{r-n}\to0\eqno{(12.9.2)}$$ which defines us the dual lattice $L'_T$ and hence the lattice $L'\hookrightarrow \p^{r-n}$.
\medskip
{\bf Proposition 12.9.3.} The Siegel matrix of $L'$ is $-S^t$. 
\medskip
{\bf Proof.} Let $l_1, \dots, l_r$ be a $\n F_q[T]$-basis of $L$, and $l'_1, \dots, l'_r$ the dual $\n F_q[T]$-basis of $L'$. Recall that 
$\forall \ i=1,\dots, r-n$ we have $l_{n+i}=\sum_{k=1}^n s_{ik}l_k$, i.e. elements 
$$v_i:=l_{n+i}-\sum_{k=1}^ns_{ik}l_k,\ \ \ \ i=1,\dots, r-n$$ form a basis of $\Ker \ze$. 

Let us consider elements $$w_j:=l'_j+\sum_{m=1}^{r-n}s_{mj}l'_{n+m},\ \ \ \ j=1,\dots, n$$ We have $\forall \ i, \ j \ \ <v_i, w_j>=0$, hence the desired. $\square$
\medskip

{\bf Theorem 12.10} - The duality theorem: The functors of lattice and duality commute. Namely, if $M$ is uniformizable then $M'$ is uniformizable, and $L(M')=(L(M))'$.
\medskip
{\bf Proof.} Recall that here we consider the case $N=0$. For the case $N\ne0$ see Section 14.  The proof follows [GL07]. We consider bases $\{l_*\}$, $\{f_*\}$, $\{l'_*\}$, $\{f'_*\}$. We can rearrange elements of $\{l_*\}$ by a manner that the $\ze$-images of $l_1, \dots, l_n$ form a $\p$-basis \footnotemark \footnotetext{For the case when $\ze$-images of $l_1, \dots, l_n$ do not form a $\p$-basis of $Lie(M)$ see Section 19.} of $Lie(M)$.\footnotemark \footnotetext{We follow the proof of [GL07], Section 5.4. But in [GL07] analogs of $l_1, \dots, l_n$ are the {\it last} elements of the set $l_1, \dots, l_r$. This fact causes some changes in the below formulas. Sorry for such inconvenience.} The Siegel matrix $S\in M_{(r-n)\times n}(\p)$ is defined as follows: $$\left(\matrix \ze(l_{n+1}) \\ \dots \\ \ze(l_{r}) \endmatrix \right)=S\left(\matrix \ze(l_{1}) \\ \dots \\ \ze(l_{n}) \endmatrix \right)\eqno{(12.10.1)}$$ We shall show in 12.10.8 that in this case the $\ze$-images of $l'_{n+1}, \dots, l'_r$ form a $\p$-basis of $Lie(M')$, and we consider the corresponding Siegel matrix $S'\in M_{n\times(r-n)}(\p)$ defined as follows: $$\left(\matrix \ze(l'_{1}) \\ \dots \\ \ze(l'_{n}) \endmatrix \right)=S'\left(\matrix \ze(l'_{n+1}) \\ \dots \\ \ze(l'_{r}) \endmatrix \right)\eqno{(12.10.2)}$$ The statement of (12.10) is equivalent to $S'=-S^t$, so we shall prove this fact. 
\medskip
Now we need an important relation between $\Psi$ and $S$, this is formula (12.10.5) below. 
\medskip
We denote by $D=\{d_{ij}\}$, respectively $D'=\{d'_{ij}\}\in M_{r\times r}(\p)$ the coefficients at the degree $-1$ of the $\th$-shift of $\Psi$, respectively of $\Psi'$. (12.10.1) implies 
$$\forall \ k, j \ \ \ \partial_{\ze(l_{n+k})}(f_j)=\sum_{i=1}^ns_{ki}\partial_{\ze(l_{i})}(f_j)\eqno{(12.10.3)}$$ Because of (11.10), this is equivalent to 
$$\forall \ k, j \ \ \ -d_{j,n+k}=\sum_{i=1}^ns_{ki}d_{ji}\eqno{(12.10.4)}$$
Let us denote by $D_\la\in M_{r\times n}$, respectively $D_\rho\in M_{r\times (r-n)}$ (here $\la, \ \rho$ mean "left; right") a submatrix of $D$ formed by all lines and the first $n$ columns, respectively by all lines and the last $r-n$ columns of $D$. We define analogous submatrices   $D'_\la, \ D'_\rho$ of $D'$. (12.10.4) becomes 
\medskip
{\bf Relation between $S$ and $\Psi$}
$$-SD_\la^t=D_\rho^t\eqno{(12.10.5)}$$
\medskip
The $\th$-shift of $\Xi^{-1}$ starts from the term $c_{-1}N^{-1}$ (I think that $c_{-1}=1$; to check). Since $\Psi'\Psi^t=\Xi^{-1}$, we get that $D'D^t=0$. 
\medskip
We have $$0=D'D^t=D'_\la D^t_\la+D'_\rho D^t_\rho=D'_\la D^t_\la-D'_\rho SD_\la^t=(D'_\la -D'_\rho S)D_\la^t\eqno{(12.10.6)}$$
Since $D_{\la}^t$ is a $n\times r$-matrix of rank $n$, it is not a zero-divisor from the right, so $$D'_{\la}=D'_{\rho}S, \hbox{ i.e. } S^t{D'_\rho}^t={D'_{\la}}^t\eqno{(12.10.7)}$$ Since the rank of $D'$ is $r-n$ and $D'_{\rho}$ is a
$r\times (r-n)$ matrix, (12.10.7) implies that 
\medskip
{\bf 12.10.8.} Columns of $D'_{\rho}$ are linearly independent. Hence, the $\ze$-images of $l'_{n+1}, \dots, l'_r$ form a $\p$-basis of $Lie(M')$. 
\medskip
We see that (12.10.7) is an analog of (12.10.5) for $M'$, i.e. $-S^t$ is a Siegel matrix of
$M'$. $\square$
\medskip
{\bf 12.11. Theorem of Anderson, complete form.} 
\medskip
Theorem 10.5 is a non-complete form of a theorem of Anderson. To give its complete form, we need one more notion. Let us consider $M\{T\}$ and $H^1(M)=M\{T\}^\tau$ from (12.3.2.1), (12.3.2.2). 
\medskip
{\bf Definition 12.11.1.} $M$ is called rigid-analytically trivial if $$H^1(M)\underset{\n F_q[T]}\to{\otimes}\p\{T\}=M\{T\}$$ ([A86], (2.3.3); [G96], Definition 5.9.13). 
\medskip
{\bf Theorem 12.11.2.} ([A86], Theorem 4; [G96], Theorem 5.9.14). Let $M$ be a t-motive of rank $r$. The following four conditions are equivalent: 
\medskip
{\bf 1.} $h_1(M)=r$.
\medskip
{\bf 2.} $exp_M$ is surjective.
\medskip
{\bf 3.} $M$ is rigid-analytically trivial.
\medskip
{\bf 4.} $h^1(M)=r$.
\medskip
We refer to the proofs to [A86], [G96] (to write some lacking details!) We give only the following remarks. First, {\bf (3)} trivially implies {\bf (4)}, but the converse implication is not trivial. Further, {\bf (1) $\iff$ (4)}, this is Theorem 12.5.1. Finally, the implication {\bf (2) $\implies$ (1)} is simple (its proof is of few lines long), but the proofs of the remaining implications {\bf (3) $\implies$ (2)}, {\bf (1) $\implies$ (3)} are more complicated. 
\medskip
{\bf Remark.} To check where Anderson and Goss prove that $L(M)$ satisfies 4.1.2 (the linear independence over $\r$).
\medskip
{\bf 13. The lattice functor.}
\medskip
Recall that in characteristic 0 the lattice functor from abelian varieties to lattices satisfied the Riemann condition is an isomorphism. For the characteristic $p$ the situation is worse. Let us describe what is known.
\medskip
First of all, not all Anderson t-motives $M$ are uniformizable, i.e. not for all of them the dimension of $L(M)$ is $r$ as it should be. An example of a non-uniformizable t-motive is given in [A86], 2.2 ( = [G96], 5.9.9). The proof of its non-uniformizability given in [A86] is "artificial". The same example is treated in [EGL] by a direct calculation, see [EGL], Example 6.3. It is shown that the example of [A86], 2.2 is, in some sense, the "simplest" non-uniformizable t-motive.
\medskip
Let us consider notations of (5.9). It is known that if $A_k$ is invertible and for $i=1,\dots, k-1$ the entries of $A_i$ are sufficiently small then $M$ are uniformizable. An explicit estimate is given in [GL17] (end of page 383 and Proposition 2): 
\medskip
{\bf Theorem 13.1.} Let $M$ be given by (5.9) with $k=2$, $A_0=\th I_n$, $A_2=I_n$. If $v_\infty$ of all entries of $A_1$ is $> \frac q{q^2-1}$ then $M(A)$ is uniformizable. 
\medskip
Clearly this estimate is too weak and can be improved, but we do not know the best result. 
\medskip
As we indicated above (Theorem 8.3, (2); see [D76] for a proof), the lattice functor is 1 -- 1 on Drinfeld modules. 
\medskip
Further, it is necessary to have in mind that not all lattices have dual. Really, let us consider the case $n=1$, any $r$, i.e. the case of lattices $L=(\n F_q[\th])^r$ in $\p$. A matrix $S=\left(\matrix s_1\\ \dots \\ s_{r-1}\endmatrix\right)\in M_{(r-1)\times 1}(\p)$ is a Siegel matrix of such lattice iff 
$$1, s_1, \dots , s_{r-1}\hbox{ are linearly independent over } \r\eqno{(13.2)}$$ while its transposed $S^t=\left(\matrix s_1& \dots & s_{r-1}\endmatrix\right)\in M_{1\times (r-1)}(\p)$ is a Siegel matrix of a lattice of rank $r$ in $\p^{r-1}$ iff 
$$\exists \ i \hbox{ such that } s_i\not\in\r\eqno{(13.3)}$$
For $r>2$ (13.3) is weaker than (13.2), i.e. not all lattices of rank $r$ in $\p^{r-1}$ have duals. 
\medskip
{\bf Remark 13.4.} The same phenomenon occurs for Siegel matrices of other sizes. For example, let $S=\{s_{ij}\}\in M_{2\times 2}(\p)$ be a matrix such that $\vi(a_{11})\not\in \n Z$,  $\vi(a_{21})\not\in \n Z$,  $\vi(a_{11})-\vi(a_{21})\not\in \n Z$. Then obviously $S$ is a Siegel matrix of a lattice $(\n F_q[\th])^4\subset\p^2$. But if $\ \ 1, \ a_{11}, \ a_{12}$ are linearly dependent over $\r$ then $S^t$ is not a Siegel matrix of a lattice. 
\medskip
We think that it is necessary to modify the Definition 4.1 of a lattice, in order to exclude this phenomenon. 
\medskip
The theory of duality ([GL07]) gives us more information. Namely, for $n=r-1$ we have an immediate corollary of Theorem 8.3:
\medskip
{\bf Theorem 13.5.} ([GL07], Corollary 8.4). All pure t-motives of rank $r$ and dimension $r-1$ over $\p$ having $N=0$ are uniformizable. There is a 1 -- 1 correspondence between their set, and the set of lattices of rank $r$ in $\p^{r-1}$ having dual.
\medskip
Further, we have an important theorem ([HJ], Theorem 3.34 (b)):
\medskip
{\bf Theorem 13.6.} The lattice functor $$M\mapsto \{\ L(M)\hookrightarrow Lie(M)\ \}$$ from the category of pure uniformizable t-motives $M$ to the category of lattices is fully faithful , i.e. the lattice functor is injective. 
\medskip
The condition of mixedness is necessary for injectivity. It is proved in [GL24] that different $M$ that enter in (5.17.4) are non-isomorphic, but their lattices are isomorphic. 
\medskip
{\bf 13.7. Surjectivity of the lattice functor.} We have a result of local surjectivity ([GL17], Proposition 2 (p. 384)). We consider, from one side, Anderson t-motives given by (5.13.1) such that the matrix $A$ is in a neighborhood of 0. From another side, we consider lattices whose Siegel matrix is in a neighborhood of $\omega I_n$ (here $\omega\in \n F_{q^2}-\n F_q$; a Siegel matrix of an Anderson t-motive given by (5.13.1) with $A=0$ is exactly $\omega I_n$). We prove that the lattice map is a surjection in a neighborhood of a Siegel matrix $\omega I_n$.
\medskip
We think that explicit calculations similar to the ones of [GL24] will give us an answer to the problem of surjectivity of the lattice map for the case $r=3, \ n=2$. If a lattice has dual (i.e. its Siegel matrix satisfies (13.2)) then it is the lattice of a pure t-motive. But if it satisfies (13.3) but not (13.2) the question is open. 
\medskip
{\bf 13.8.} Finally, let us consider the problem of explicit calculation of $h^1(M), \ h_1(M)$. Some cases are treated in [EGL], [GL21]. These are explicit calculations similar to the ones of Section 12.4, but more complicated. For example, formulas [GL21], (3.9) permit us to find $h^1(M)$ for any $M(A)$ of the form (5.13.1), case $n=2$. Paper [EGL] gives a complete answer to the problem of finding of $h^1(M), \ h_1(M)$ for $M(A)$ of the form (5.13.1) where $A=\left(\matrix 0&a_{12}\\ a_{21}&0\endmatrix \right)$.
\medskip
Results of some computer calculations, as well as hand calculations made by students of the course given by the second author in POSTECH, suggest that for "almost all" $M$ we have $h^1(M)=h_1(M)=r$. But it is necessary to indicate what is "almost all"!
\medskip
Hence, there is a natural problem to continue calculations of [EGL], to extend them to all $M$ of the form (5.13.1). Advance in this activity is a good (and easy) research problem for students. It can be compared with the problems of catalogizations of all stars in the sky in astronomy, or of all species of insects in biology: it is necessary to consider a lot (probably infinite quantity) of cases. Like any biologist who comes to Manaus jungle (surrounding the home university of the second author) will find, without doubt, some new species of insects, any mathematician who will consider new types of matrices $A$ will get new results in calculation of $h^1$, $h_1$ of the corresponding $M(A)$.
\medskip
{\bf 14. Case $N\ne0$.}
\medskip
Here we consider the case $N\ne0$. Let $M$ be an Anderson t-motive. Recall (see beginning of Section 5) that $T-\th$ is nilpotent on $M/\tau M$. Recall that $T$ acts on $E(M)$ and $\Lie(M)$, and on $\Lie(M)$ we have $T=N+\th  I_n$ where $N$ is nilpotent. Let $\g m$ be the minimal number such that $N^\g m=0$. The definition of $N$-lattice was given in 10.5.1.
\medskip
Let us consider the case of uniformizable $M$, and let as earlier $L(M)$ be the kernel of $exp$. We have: $L(M)$ is a lattice of $\g V=\Lie(M)$ in the meaning of Definition 10.5.1.
\medskip
The next important object is due to Pink [P], see also [GL07], [GL18], [HJ]. There is a natural inclusion $\p[T]\to\p[[N]]$: \ \ $T\mapsto N+\th$. It defines an inclusion $\n F_q[T]\hookrightarrow\p[[N]]$. Hence, we can consider $L(M)\underset{\n F_q[T]}\to{\otimes}\p[[N]]$. An inclusion $L(M)\hookrightarrow \Lie(M)$ defines a map $L(M)\underset{\n F_q[T]}\to{\otimes}\p[[N]] \to \Lie(M)$ which is a surjection because of (10.5.4). Its kernel is denoted by $\g q=\g q(M)$. It is an invariant of $M$ which is more "thin" than $L(M)$ for the case $N\ne0$. Hence, we have an exact sequence

$$0\to \g q\to L(M)\underset{\n F_q[T]}\to{\otimes}\p[[N]]\to \Lie(M)\to0\eqno{(14.1)}$$

We have: $\p[[N]]$ is a local Dedekind ring, the theory of finitely generated modules over these rings is very simple. $L(M)\underset{\n F_q[T]}\to{\otimes}\p[[N]]$ is free of rank $r$ over $\p[[N]]$, and $\Lie(M)$ is isomorphic to $\p^n$. (14.1) gives us that $\g q$ is isomorphic to $\p[[N]]^r$ as a $\p[[N]]$-module, i.e. it is a lattice in $L(M)\underset{\n F_q[T]}\to{\otimes}\p[[N]]$. It is checked immediately that $\g q\supset N^\g m (L(M)\underset{\n F_q[T]}\to{\otimes}\p[[N]])$.
\medskip
{\bf Remark 14.2.} We can define an abstract Hodge-Pink structure as follows (see, for exanple, [HJ]). Let $V$ be any free $\n F_q[T]$-module\footnotemark \footnotetext{Do not confuse this $V$ and the above $\g V$, they play different roles.} of rank $r$ and $\g q\subset V\underset{\n F_q[T]}\to{\otimes}\p[[N]]$ a $\p[[N]]$-sublattice of rank $r$. A homomorphism between two such objects $(V_1, \g q_1)$, $(V_2, \g q_2)$ is a $\n F_q[T]$-linear map $V_1\to V_2$ inducing a map $\g q_1 \to \g q_2$.
\medskip
Really, a definition of Hodge-Pink structure of [HJ] is a version of the above one. First, they consider a weight filtration on $V$. Second, they consider not a free $\p[[N]]$-module but a $\p((N))$-vector space; particularly, their $\g q$ is not necessarily contained in $V\underset{\n F_q[T]}\to{\otimes}\p[[N]]$ but only is commensurable with it.
\medskip
$\{\g q=\g q_M, L(M)\}$ from (14.1) are called the Hodge-Pink structure associated to $M$ (recall that $M$ is uniformizable).
\medskip
In order to describe explicitly the inclusion $\g q\subset V\underset{\n F_q[T]}\to{\otimes}\p[[N]]$, we use an analog of a Siegel matrix. Namely, for the case $\g m=1$ we get exactly a Siegel matrix. Really, let $e_1, \dots, e_r$ be a basis of $V$ over $\n F_q[T]$ such that images of $e_1, \dots, e_n$ form a basis of $V\underset{\n F_q[T]}\to{\otimes}\p[[N]]/\g q$. There exists a matrix $S=\{s_{ij}\}$ such that elements
$$Ne_1, \dots, Ne_n, \ \ \ e_{n+i}-\sum_{j=1}^nS_{ij}e_j, \hbox { where } i=1, \dots, r-n$$
form a $\p[[N]]$-basis of $\g q$. Clearly this $S$ is exactly a Siegel matrix of $L(M)$ if $\g q$ comes from $L(M)$, as explained above.
\medskip
If $\g m>1$ then an analog of Siegel matrix is a family of matrices called a Siegel object, it is defined for example in Sections 19.6.5, 19.6.6, or in [GL18], (3a.5.4.1). Roughly speaking, we consider the $\p$-vector spaces
\medskip
$((N^i\cdot V\underset{\n F_q[T]}\to{\otimes}\p[[N]] )\cap \g q ) / ((N^{i+1}\cdot V\underset{\n F_q[T]}\to{\otimes}\p[[N]]) \cap \g q )$, $i=0,\dots, \g m-1$,
\medskip
their bases, lifts of these bases to $\g q$ etc, see [GL18] for details. The family of these matrices is parametrized by integer points of a tetrahedron, i.e. it depends on 3 integer parameters. These matrices are denoted by $S_{uvyz}$ in [GL18] and Section 19.6.6, where $u,\ v, \ y,\ z$ are integer parameters.\footnotemark \footnotetext{This 4-parameter notation is used for convenience for a proof of duality theorem. Really, always $v=u-1$, i.e. practically we have 3 parameters.} 
\medskip
{\bf Remark 14.3.} The set of Siegel objects is the set of elements in the maximal Schubert cell of a generalized flag variety, see Definition 19.6.4.6. See also [GL18], 3a.1.2 for a definition of this generalized flag variety in terms of bases of $N$-lattices, and [GL18], 3a.6 for a definition in terms of Hecke cosets. 
\medskip
Let $M$ be an Anderson t-motive and $\g q\subset L(M)\underset{\n F_q[T]}\to{\otimes}\p[[N]]$ its Hodge-Pink structure. Like for the case $N=0$, there is a formula for a Siegel object of $\g q$ defined in terms of $\th$-shift of a scattering matrix $\Psi$ of $M$. Namely, we denote $$\th\hbox{-shift of }\Psi=\sum_{u=-\g m}^{\infty}D_{-u}N^u $$ where $D_1, \dots, D_\g m$ are matrix coefficients at negative powers of $N$. Formula (3.38) of [GL18] gives relations between some submatrices of $D_i$ and matrices of Siegel object of $\g q(M)$.
\medskip
There are natural definitions of the tensor product and Hom of Hodge-Pink structures. Namely, let
\medskip
$HP_1=\{\g q_1\subset V_1\underset{\n F_q[T]}\to{\otimes}\p[[N]]\ \}$, $HP_2=\{\g q_2\subset V_2\underset{\n F_q[T]}\to{\otimes}\p[[N]]\ \}$
\medskip
be two Hodge-Pink structures. Their tensor product is
\medskip
$\g q_1\underset{\p[[N]]}\to{\otimes} \g q_2 \subset (V_1 \underset{\n F_q[T]}\to{\otimes} V_2) \underset{\n F_q[T]}\to{\otimes}\p[[N]]$.
\medskip
The definition of Hom is similar, see [HJ] (it is defined only for Hodge-Pink structures in the meaning of [HJ]).
\medskip
 Let us give the definition of the $\g m$-dual structure --- an important particular case of Hom. We denote by $\g C^\g m$ (the $\g m$-th tensor power of the Carlitz Hodge-Pink structure) a structure $N^\g m \p[[N]] \subset \p[[N]]$.
\medskip
The Hodge-Pink structure $\Hom(HP_1, \g C^\g m)$ --- the $\g m$-dual of $HP_1$ --- is defined as follows. Its $V$ is $V_1':=\Hom(V_1, \n F_q[T])$ and
$$\g q'=\{\ \ \vf: V_1 \underset{\n F_q[T]}\to{\otimes} \p[[N]]\to \p[[N]]\ \ \  | \ \ \  \vf(\g q)\subset N^\g m \p[[N]]\ \ \}$$
\medskip
{\bf Theorem.} The functor of Hodge-Pink structure from uniformizable Anderson t-motives to Hodge-Pink structures commutes with tensor products and Hom's.
\medskip
For tensor products this was proved by Anderson in ~ 2000 (non-published), for tensor products of t-motives having $N=0$ in [GL07], for the general tensor products and Hom's in [HJ], 3.5b. For the case of duals and $N=0$ the theorem was proved in [GL07], for $\g m$-duals, case $N\ne0$, in [GL18] by means of explicit calculation of Siegel objects of $L(M)$ and of $L({M'}^\g m)$ (recall that ${M'}^\g m$ is the $\g m$-dual of $M$, see 5.19.3).
\medskip

{\bf 15. L-functions.}
\medskip
There are various types of L-functions of t-motives, see, for example, [B02], Definition 15; [B05], Definition 2.19; [G92], (3.4.2a); [G96], Section 8; [GLZ22], (0.1.3A,B); [L],  upper half of page 2603; [Tae], p. 371, below Remark 2. 
\medskip
We give here a definition of the simplest version of L-function of t-motives. See, for example, [TW], Section 2. 
\medskip
{\bf 15.1.} First, we consider the case of one Euler factor, i.e. the case when $M$ is defined over a finite field $\n F_{q^d}$. Let $f\in \n F_q[T]$ be an irreducible polynomial (in most steps of the proof we shall not use the fact that $f$ is irreducible). We denote the $f$-adic completion of $\n F[T]$ by $\n F_q[T]_f$, i.e. $$\n F_q[T]_f:=\underset{\underset{n\to\infty}\to{\longleftarrow}}\to{\hbox{lim}} \n F_q[T]/f^n\n F_q[T]$$
it is the function field analog of $p$-adic numbers $\n Z_p$. 
\medskip
The below theorem is an analog of the Theorem 1.5.1, i.e. we shall prove that (roughly speaking) the characteristic polynomial of the action of Frobenius on $T_f(M)$ (called the first characteristic polynomial $C_{1,M,f}(U)$ ), coincides with the characteristic polynomial of the action of Frobenius on $M$ itself (called the second characteristic polynomial $C_{2,M}(U)$ ). Particularly, this will imply that $C_{1,M,f}(U)$, which a priori belongs to $\n F_q[T]_f[U]$, in reality belongs to $\n F_q[T][U]$ and does not depend on $f$. 
\medskip
Let us define $C_{2,M}(U)$. We consider $Q$ from (5.11.1). We have $Q\in M_{r\times r}(\n F_{q^d}[T])$. 
\medskip
{\bf Definition 15.2.} $C_{2,M}(U):=\det(I_r-Q^{[d]}U)$. 
\medskip
{\bf Remark 15.3.} $C_{2,M}(U)\in \n F_q[T][U]$, because $Q^{(d)}=Q$ implies $$C_{2,M}(U)^{(1)}=\det(I_r-Q\cdot Q^{[d-1]}\cdot Q^{[d-2]}\cdot...\cdot Q^{[2]}\cdot Q^{[1]}U).$$ We have $\det(I_r-AB)=\det (I_r-BA)$, hence $C_{2,M}(U)^{(1)}=\det(I_r-Q^{[d]}U)=C_{2,M}(U)$. 
\medskip
\medskip
{\bf Theorem 15.4.} $C_{1,M,f}(U)=C_{2,M}(U)$. 
\medskip
{\bf Remark 15.4.1.} $C_{1,M,f}(U)$ comes from the action of Frobenius on $E=E(M)$, and $E$ is defined by the formula (5.9) --- the action of $T$ on a $\p[\tau]$-basis of $M$. From another side, $C_{2,M}(U)$ comes from a matrix $Q$ from (5.11.1), it describes the action of $\tau$ on a $\p[T]$-basis of $M$. Hence, $C_{1,M,f}(U)$ and $C_{2,M}(U)$ are of different nature, and their coincidence does not follow immediately from definitions. See Remark 15.4.6 below.
\medskip
Below in 15.5 we give an example of a "direct" calculation that shows that $C_{1,M,f}(U)=C_{2,M}(U)$ in a particular case. We see that it is non-trivial. 
\medskip
{\bf Proof.} See [G96], Section 5.6, especially Proposition 5.6.9. Let us reproduce it explaining some details. We start from a t-motive $M$ over $\p[T,\tau]$ and any element $f\in \n F_q[T]$. Let as earlier $E=E(M)$. Let 

$$W_f:=\n F_q[T]/ f \ \n F_q[T]$$

We need its dual $W_f^*:=\Hom_{\n F_q}(W_f, \n F_q)$ (it is equal to $f^{-1}\Omega / \Omega$ where $\Omega=\n F_q[T]\ dT$, see [A86], p. 466 - 467, and [G96], p. 153, but we do not need this fact). Further, let $$V_f:=\Hom_{\n F_q}(W_f, \p)=W_f^*\underset{\n F_q}\to{\otimes}\p$$
($V_f$ is denoted by $V$ in [G96], p. 153). 
\medskip
$V_f$ has a natural structure of a left $\p[T,\tau]$-module (see [A86], proof of Proposition 1.8.3, and [G96], p. 153; it is the only such natural structure, and the reader can easily define it himself). 
\medskip
We have: 
$$E_f=\Hom_{\p[\tau]}(M/fM, \p)\eqno{(15.4.2)}$$
$$E_f=\Hom_{\n F_q}(\ (M/fM)^\tau, \n F_q)\eqno{(15.4.3)}$$
$$E_f=\Hom_{\p[T,\tau]}(M/fM, V_f)\eqno{(15.4.4)}$$
$$E_f=\Hom_{\n F_q[T]}( (M/fM)^\tau, W_f^*)\eqno{(15.4.5)}$$

{\bf Remark 15.4.6.} The non-formal meaning of these equalities is the following. We start from $\p[\tau]$-modules and come to $\n F_q[T]$-modules. This will give us a transition from $C_{1,M,f}(U)$ defined in terms of the $\p[\tau]$-action on $M$, to $C_{2,M}(U)$ defined in terms of the $\p[T]$-action on $M$.
\medskip
To prove (15.4.2), we apply the functor $* \mapsto \Hom_{\p[\tau]}(*, \p)$ to the exact sequence 
$$0\to M\overset{f}\to{\to}M\to M/ fM\to0$$ and we use Definition 6.3.
\medskip
To deduce (15.4.3) from (15.4.2), we apply the Lang's theorem to $M/fM$, see the formulas of [G96], proof of Proposition 5.6.3 for details; \footnotemark \footnotetext{The last line of [G96], proof of Proposition 5.6.3 contains a typographical error: must be $M/fM$ instead of $M/\tau M$.}
\medskip
Deduction (15.4.3) $\to$ (15.4.4) is [G96], Lemma 5.6.5 and [A86], a line in the proof of Proposition 1.8.3 (both are given as an exercice for the reader); 
\medskip
Deduction (15.4.4) $\to$ (15.4.5) is again the application of the Lang's theorem to $M/fM$, see [A86], Proposition 1.8.3, and [G96], Theorem 5.6.6.
\medskip
Recall that if $f$ is irreducible then the Tate module $T_f(M):=\underset{\underset{n\to\infty}\to{\longleftarrow}}\to{\hbox{lim}} E_{f^n}$. We have
$$T_f(M)=\Hom_{\n F_q[T]_f}( (\underset{\underset{n\to\infty}\to{\longleftarrow}}\to{\hbox{lim}} M/f^nM)^\tau, \underset{\underset{n\to\infty}\to{\longleftarrow}}\to{\hbox{lim}}\ W_{f^n}^*)\eqno{(15.4.7)}$$ this is a limit of (15.4.5). Obviously $\underset{\underset{n\to\infty}\to{\longleftarrow}}\to{\hbox{lim}}\ (M/f^nM)^\tau= (\underset{\underset{n\to\infty}\to{\longleftarrow}}\to{\hbox{lim}} M/f^nM)^\tau$. 
\medskip
Now let us consider the Galois actions. Let $M$ be defined over a field $L_1$. This means that $M$ is a left module over $L_1[T,\tau]$ satisfying (5.2.1 -- 5.2.3). We have: $A_i$ from (5.9) belong to $M_{n\times n}(L_1)$, $Q$ from (5.11.1) belongs to $M_{r\times r}(L_1[T])$. 
\medskip
Let $L$ be an algebraic closure of $L_1$, and $G:=\Gal(L/L_1)$. We consider $\overline{M}:=M\underset{L_1}\to{\otimes}L$. In the above formulas (15.4.2) --- (15.4.5), (15.4.7) we can substitute $\p$ by $L$, and $M$ by $\overline{M}$. 
\medskip
$G$ acts on $\overline{M}$ via $L$, i.e. it acts on coefficients. This action commutes with multiplication by $T$ and $\tau$, hence $G$ acts on $\overline{M}/ f^n\overline{M}$, $(\overline{M}/ f^n\overline{M})^\tau$. 
\medskip
According the definition 6.3, we have $E=\Hom_{L[\tau]}(\overline{M}, L)$. Hence, we can define an action of $G$ on $E$, as usual for the action of groups on sets of Hom's: $$g\vf(m):=g(\vf(g^{-1}(m)))$$ Further, the action of $G$ is well-defined on $E_f$ and on $T_f(M)$. It is straightforward to check that formulas (15.4.2) --- (15.4.5), (15.4.7) are compatible with this action of $G$. 
\medskip
Now let us consider the case when $L_1=\n F_{q^d}$ is a finite field. We denote by $g_0$ the Frobenius of $G$, i.e. $g_0(x)=x^{q^d}$. It is useful to distinguish three "Frobenius" endomorphisms of $\overline{M}$, we denote them by $\g F_1, \ \g F_2, \ \g F_3$ respectively, although we shall need only $\g F_1$ (see Remark 15.4.12). They are defined as follows (here as earlier $e_*$, respectively $f_*$ are bases of $M$ over $\n F_{q^d}[\tau]$, respectively over $\n F_{q^d}[T]$, and $l_*\in \overline{\n F_q}$):
$$\g F_1(\sum_{i, j}\ l_{ij}\tau^je_i):=\sum_{i, j}\ l_{ij}^{q^d}\tau^je_i; \ \ \  \g F_1(\sum_{i, j}\ l_{ij}T^jf_i)=\sum_{i, j}\ l_{ij}^{q^d}T^jf_i;$$
$$\g F_2(\sum_{i, j}\ l_{ij}\tau^je_i):=\sum_{i, j}\ l_{ij}\tau^{j+d}e_i; $$
$$\g F_3(\sum_{i, j}\ l_{ij}\tau^je_i):=\sum_{i, j}\ l_{ij}^{q^d}\tau^{j+d}e_i; $$
We have: all $\g F_1, \ \g F_2, \ \g F_3$ commute with multiplication by $T$ and $\tau$, $\g F_1$ is the above action of $g_0$ on $\overline{M}$; 
\medskip
$\g F_2$ is a $\overline{\n F_q}[T,\tau]$-linear endomorphism of $\overline{M}$, i.e. $\g F_2$ commutes with the multiplication by any $l \in \overline{\n F_q}$;
\medskip
$\g F_3(x)=\tau^d\cdot x$ for $x\in \overline{M}$, and $\g F_3=\g F_1\circ \g F_2=\g F_2\circ \g F_1$. 
\medskip
Since $\g F_1, \ \g F_2, \ \g F_3$ commute with multiplication by $T$ and $\tau$, they act on $(\overline{M}/ f^n\overline{M})$, $(\overline{M}/ f^n\overline{M})^\tau$, $ (\underset{\underset{n\to\infty}\to{\longleftarrow}}\to{\hbox{lim}} \overline{M}/f^n\overline{M})^\tau$. 
\medskip
The first characteristic polynomial $C_{1,M,f}(U)$ is defined in terms of the action of $g_0$ on $T_f(\overline{M})$. Applying (15.4.7) to $\overline{M}$ we get that $C_{1,M,f}(U)$ is defined in terms of the action of $g_0^{-1}$ on $(\underset{\underset{n\to\infty}\to{\longleftarrow}}\to{\hbox{lim}} \overline{M}/f^n\overline{M})^\tau$, because the action of $G$ on $\underset{\underset{n\to\infty}\to{\longleftarrow}}\to{\hbox{lim}}\ W_{f^n}^*$ is trivial. 
\medskip
Let $h_*=(h_1,\dots,h_r)^t$ be a basis of $(\underset{\underset{n\to\infty}\to{\longleftarrow}}\to{\hbox{lim}} \overline{M}/f^n\overline{M})^\tau$ over $\n F_q[T]_f$. Let $Z\in M_{r\times r}(\n F_q[T]_f)$ be the matrix of the action of $g_0$ ( = the action of $\g F_1$ ) on $h_*$. We see that $$C_{1,M,f}(U)=\det(I_r-Z^{-1}U).\eqno{(15.4.8)}$$
\medskip
Recall that the matrix of $\tau$ in $f_*$ is denoted by $Q\in M_{r\times r}(\n F_{q^d}[T])$. We consider the same basis $f_*$ as a basis of $\underset{\underset{n\to\infty}\to{\longleftarrow}}\to{\hbox{lim}} \overline{M}/f^n\overline{M}$ over $\overline{\n F_q}[T]_f$, and the image of $Q$ in $M_{r\times r}(\overline{\n F_q}[T]_f)$. Let $X\in GL_r(\overline{\n F_q}[T]_f)$ be the (transposed) matrix of change of basis from $f_*$ to $h_*$, i.e. $h_*=X \ f_*$. 
\medskip
Hence, we have two bases $h_*$ and $f_*$, two operators $g_0=\g F_1$ and $\tau$, and three matrices $Z, \ X, \ Q$ such that 
$$\tau f_*=Qf_*, \ \ \ \ g_0h_*=Zh_*, \ \ \  \ h_*=Xf_*$$
They satisfy relations: 
$$X^{(1)}Q=X \hbox{    (equality in } \overline{\n F_q}[T]_f)\eqno{(15.4.9)}$$
(because $h_*$ is $\tau$-invariant), and
$$g_0(h_*)=X^{(d)}g_0(f_*)=X^{(d)}f_*; \ \  g_0(h_*)=Zh_*=ZXf_* \implies  X^{(d)}=ZX.\eqno{(15.4.10)}$$
(15.4.9) implies $$X^{(d)}Q^{[d]}=X \hbox{    (equality in } \overline{\n F_q}[T]_f)..\eqno{(15.4.11)}$$
(15.4.10), (15.4.11) imply
$$Z=X(Q^{[d]})^{-1}X^{-1}.$$
This means that $\det(I_r-Z^{-1}U)=\det(I_r -Q^{[d]}U)$, hence --- because of (15.4.8) --- we get the desired. $\square$
\medskip
{\bf Remark 15.4.12.} Unlike the number field case (see Idea of the proof of Theorem 1.5.1), the endomorphism $\g F_1$ used in the above proof is not an algebraic endomorphism of $\overline{M}$: it does not commute with multiplication by $l\in \bar \n F_q$. We should not confuse it with $\g F_2$ which is an algebraic endomorphism of $\overline{M}$. By the way, like in the number field case, for almost all Drinfeld modules $M$ over $\n F_{q^d}$ we have $\End(M)=\n F_q[T][\g F_2]$, see [G96], Section 4. 
\medskip

{\bf 15.5. Example: Explicit calculations.}
\medskip
In order to show that the Theorem 15.4 is non-trivial, we give here its calculational proof in two simple cases. 
\medskip
{\bf Case A.} $M$ is a Drinfeld module of rank 2, $d=1$, and $f=T$.
\medskip
Let $a_0, a_1, a_2$ be from (5.6), hence $a_0=\th$, and let $a_2=1$. The matrix $Q$ is given in (5.12.1), $Q^{[1]}=Q$, and $$C_{2,M}(U)=1+a_1U+(a_0-T)U^2\eqno{(15.5.1)}$$
\medskip
The $T$-Tate module $T_T(M)$ is the set of concordant elements $x_0, \ x_1, \dots\in \overline{\n F_q}$. 
Conditions of concordance are 
$$x_0^{q^2}+a_1x_0^q+a_0x_0=0\eqno{(15.5.2)}$$ 
$$x_1^{q^2}+a_1x_1^q+a_0x_1-x_0=0\hbox{ etc.}\eqno{(15.5.3)}$$
We associate it an element $y:=x_0+x_1T+x_2T^2+...\in \overline{\n F_q}[[T]]$. Conditions (15.5.2 --- 15.5.3) etc. are equivalent to the condition $$y^{(2)}+a_1y^{(1)}+(a_0-T)y=0\eqno{(15.5.4)}$$
$T_T(M)$ is the set of solutions of (15.5.4). It is a free $\n F_q[[T]]$-module of dimension 2. Let $y_1, \ y_2$ be its basis, and $y_*:=\left(\matrix y_1\\ y_2\endmatrix \right)$.
\medskip
Let us find the characteristic polynomial of the action of Frobenius $y\mapsto y^{(1)}$ in this basis. 
 Let $C\in M_{2\times 2}(\n F_q[[T]])$ be its matrix in $y_*$, i.e. $y_*^{(1)}=Cy_*$. Hence, $y_*^{(2)}=C^2y_*$ and the equation (15.5.4) implies $C^2+a_1C+(a_0-T)I_2=0$, i.e. (the Cayley–Hamilton theorem) $\det(U\cdot I_2 -C)=U^2+a_1U+a_0-T$. It is equivalent to (15.5.1). 
\medskip
{\bf Case B.} $M$ is a Drinfeld module of rank 2, $d=2$, and $f=T$. We use the above notations. In this case $a_1,  a_0\in \n F_{q^2}$, and the Frobenius is $y\mapsto y^{(2)}$. As above we consider $\det(U\cdot I_2-Q^{[2]})$ instead of $\det(I_2-Q^{[2]}U)$; this is the same. We have
\medskip
$\tilde Q^{[2]}=\left(\matrix 0&1\\ T-a_0^q&-a_1^q\endmatrix \right)\left(\matrix 0&1\\ T-a_0&-a_1\endmatrix \right)=\left(\matrix T-a_0&-a_1\\ -a_1^qT+a_1^qa_0&T-a_0^q+Na_1\endmatrix \right)$ 
\medskip
and (here $N$, $Tr$ are the norm and the trace of $\n F_{q^2}/\n F_q$): $$\det(U\cdot I_2-\tilde Q^{[2]})=U^2+(-2T+Tr(a_0)-N(a_1))U+(T^2-Tr(a_0)T+ N(a_0))\eqno{(15.5.5)}$$
We have $y_*^{(2)}=Cy_*$, $C \in M_{2\times 2}(\n F_q[[T]])$. We show that $\det(U\cdot I_2-C)$ is (15.5.5) only for the $T$-free terms. Let $y_i=x_{i0}+x_{i1}T+x_{i2}T^2+...$, where $i=1, \ 2$, and $c_{ij}=c_{ij0}+c_{ij1}T+c_{ij2}T^2+...$, where $i,j=1, \ 2$, $c_{***}\in \n F_q$. Formulas for the $T$-free terms are 

$$x_{10}^{q^2}+a_1x_{10}^q+a_0x_{10}=0 \hbox{ and } x_{10}^{q^2}=c_{110}x_{10}+c_{120}x_{20},$$ hence
$$x_{10}^q=-\frac{a_0+c_{110}}{a_1}x_{10}-\frac{c_{120}}{a_1}x_{20}\eqno{(15.5.6)}$$ 
and analogously 
$$x_{20}^q=-\frac{c_{210}}{a_1}x_{10}-\frac{a_0+c_{220}}{a_1}x_{20}\eqno{(15.5.7)}$$ Hence, 
$$x_{10}^{q^2}=-\frac{a_0^q+c_{110}}{a_1^q}x_{10}^q-\frac{c_{120}}{a_1^q}x_{20}^q=\eqno{(15.5.8)}$$ 
$$=-\frac{a_0^q+c_{110}}{a_1^q}(-\frac{a_0+c_{110}}{a_1}x_{10}-\frac{c_{120}}{a_1}x_{20}) -\frac{c_{120}}{a_1^q}(-\frac{c_{210}}{a_1}x_{10}-\frac{a_0+c_{220}}{a_1}x_{20})=\eqno{(15.5.9)}$$ 
$$=\frac{(a_0^q+c_{110})(a_0+c_{110})+c_{120}c_{210}}{N(a_1)}x_{10}+\frac{(a_0^q+c_{110})c_{120}+(a_0+c_{220})c_{120}}{N(a_1)}x_{20}\eqno{(15.5.10)}$$ 
Comparing with $x_{10}^{q^2}=c_{110}x_{10}+c_{120}x_{20}$ we get 
$$\frac{(a_0^q+c_{110})(a_0+c_{110})+c_{120}c_{210}}{N(a_1)}=c_{110}\eqno{(15.5.11)}$$ 
$$\frac{(a_0^q+c_{110})c_{120}+(a_0+c_{220})c_{120}}{N(a_1)}=c_{120}\eqno{(15.5.12)}$$
Let $C_0$ be the $T$-free term of $C$. (15.5.12) gives us $Tr(C_0)=N(a_1)-Tr(a_0)$ and (15.5.11) gives us $\det(C_0)=N(a_0)$. We get a concordance with the $T$-free part of (15.5.5). 
\medskip
The above method is a naive school-style elimination. The equivalent "scientific" method for the previous calculation is the following (it will give us equality $C_{1,M,T}(U)=C_{2,M}(U)$ for all terms, unlike the equality for the $T$-free terms as above). 
\medskip
Let $y_1, \ y_2$ be as above. We consider the set of $\{F_1, \ F_2\}$, where $F_i\in \n F_{q^2}[[T]][\tau]$ such that $F_1y_1+F_2y_2=0$. It is a left $\n F_{q^2}[[T]][\tau]$-module; we denote it by $\g M$. 
\medskip
$\g M$ is generated by 4 polynomials

$$P_1=\{\tau^2+a_1\tau+A_0; \ \ \ 0\} \hbox { where $A_0$ is the above } a_0-T$$
$$P_2=\{0;\ \ \ \tau^2+a_1\tau+a_0\} \  \ \ \ \ \ \ \ \ \ \ \ \ \ \ \ \ \  \ \ \ \ \ \ \  \ \ \ \ \ \ \ \ \ \ \  \ \  $$
$$P_3=\{\tau^2-c_{11};\ \ \  -c_{12}\} \hbox { here and below } c_{ij}\in \n F_q[[T]]$$
$$P_4=\{-c_{21};\ \ \ \tau^2 -c_{22}\} \  \ \ \ \ \ \ \ \ \ \ \ \ \ \ \ \ \  \ \ \ \ \ \ \  \ \ \ \ \ \ \ \ \ \ \  \ \  $$
We have: $P_1-P_3=\{a_1\tau+A_0+c_{11}; \ \ \ \ c_{12}\}$ (this is (15.5.6));
\medskip
$P_2-P_4=\{c_{21}, \ \ \ a_1\tau+A_0+c_{22}\}$ (this is (15.5.7));
\medskip
$a_1\tau (P_1-P_3)=\{N(a_1)\tau^2+a_1(A_0^{(1)}+c_{11})\tau; \ \ \ \ a_1c_{12}\tau\}$ (this is (15.5.8); here $N(a_1):=a_1\cdot a_1^{(1)}$ and analogously for other similar notations); 
\medskip
$(A_0^{(1)}+c_{11})(P_1-P_3)=\{a_1(A_0^{(1)}+c_{11})\tau+(A_0^{(1)}+c_{11})(A_0+c_{11}); \ \ \ (A_0^{(1)}+c_{11})c_{12}\}$   
\medskip
$c_{12}(P_2-P_4)=\{c_{12}c_{21}; \ \ \ c_{12}( a_1\tau+A_0+c_{22})\}$
\medskip
$(a_1\tau - (A_0^{(1)}+c_{11}))(P_1-P_3)-c_{12}(P_2-P_4)=$
\medskip
$=\{N(a_1)\tau^2-(A_0^{(1)}+c_{11})(A_0+c_{11})-c_{12}c_{21};\ \ \  -c_{12}(A_0^{(1)}+c_{11}+A_0+c_{22})\}$
\medskip
The final step: 
\medskip
$\Cal F:=(a_1\tau - (A_0^{(1)}+c_{11}))(P_1-P_3)-c_{12}(P_2-P_4)-N(a_1)P_3=$
\medskip
$=\{-(A_0^{(1)}+c_{11})(A_0+c_{11})-c_{12}c_{21}+c_{11}N(a_1);\ \ \  -c_{12}(-N(a_1)+A_0^{(1)}+c_{11}+A_0+c_{22})\}$
\medskip
We get an element $\Cal F=\{\Cal F_1, \ \Cal F_2\}\in \g M$ such that both $\Cal F_1, \ \Cal F_2\in \n F_{q^2}[[T]]$. This means that both $\Cal F_1, \ \Cal F_2$ are 0. We get the same relations for the characteristic polynomial of $C$, now for all terms, not only for the $T$-free terms. 
\medskip
{\bf 15.6. Global case.}
\medskip
For simplicity, we consider here only Euler factors corresponding to primes of good reduction, and we omit some technical details.
\medskip
Let $M$ be defined over
$\n F_q(\theta)$. Let $\goth P$ be an irreducible polynomial in $ \n F_q[\theta]$ of degree $d$. Let us assume that $M$ has a good reduction at $\g P$. This implies that there exists a $\n F_q(\theta)[T]$-basis of $M$ such that all entries
of $Q$ are integer at $\goth P$. The reduction of $M$ at $\g P$ is denoted by $\tilde M_\g P$. The set of bad primes is denoted by $S$.
\medskip
The local $\goth P$-factor $L_\goth P(M,U)$ is, by definition, 
$$L_\g P(M,U):=(C_{2, \tilde M_\g P}(U^d))^{-1}\in \n F_q[T][[U^d]]\eqno{(15.6.1)}$$
(we need to use $U^d$ instead of $U$ for the below convergence), and the global $L$-function (factors of good reduction) is 
$$L_S(M,U):=\prod_{\goth P\not\in S} L_\goth P(M,U)\in \n F_q[T][[U]]\eqno{(15.6.2)}$$

{\bf 15.6.3. Example: the Carlitz module over $\n F_q$.} We have $Q=T-\th$, $S=\emptyset$. Let
$\goth P=a_0+a_1\th+a_2\th^2 +...+ a_d\th^d$, where $a_i\in \n F_q$.
\medskip
We have $\tilde Q^{[d]}=a_0+a_1T+a_2T^2 +...+ a_dT^d$ (exercise for the reader).
\medskip
Hence, the local $\goth P$-factor of $L(\g C,U)$ is $(1-(a_0+a_1T+a_2T^2 +...+ a_dT^d)U^d)^{-1}\in \n F_q[T][[U]]$.
For example, for $q=2$ we have a table of
$L_\goth P(\g C,U)$ for
the first small $\goth P$:
\settabs 6 \columns
\medskip
\centerline{{\bf Table 15.6.4.}}
\medskip
\+ $\goth P$ & $L_\goth P(\g C,U)$\cr
\medskip
\+ $\th$ & $1+TU+T^2U^2+T^3U^3+...$\cr
\medskip
\+ $\th+1$ & $1+(T+1)U+(T^2+1)U^2+(T+1)^3U^3+...$\cr
\medskip
\+ $\th^2+\th+1$ & $1+(T^2+T+1)U^2+(T^2+T+1)^2U^4+...$\cr
\medskip
\+ $\th^3+\th+1$ & $1+(T^3+T+1)U^3+(T^3+T+1)^2U^6+...$\cr
\medskip
\+ $\th^3+\th^2+1$ & $1+(T^3+T^2+1)U^3+(T^3+T^2+1)^2U^6+...$\cr
\medskip
and $L(\g C,U)$ is their product.
\medskip
{\bf Theorem 15.6.5.} For $q=2$ we have $L(\g C,U)=1+U$, for $q>2$ we have $L(\g C,U)=1$.
\medskip
{\bf Proof.} An analog of the classical equality of Euler
$$\sum_{n=1}^\infty \ \frac1{n^s}=\prod_{p\hbox{ prime}}\frac1{1-p^{-s}}$$
holds for this case: 
$$L(\g C,U):=\prod_{\g P}\ \frac1{1-\g PU^d}= \sum_{\g N}\ \g N\ U^d$$
where $\g P$ runs over the set of all monic irreducible polynomials in $\n F_q[T]$, $\g N$ runs over the set of all monic polynomials in $\n F_q[T]$, and $d$ is the degree of $\g P, \ \g N$. The coefficient at $U^d$ of $\sum_{\g N}\ \g N\ U^d$ is the sum of all monic polynomials in $\n F_q[T]$ of degree $d$ which is obviously 0 except the cases $d=0$ and $q=2, \ d=1$. $\square$
\medskip
There exists a formula for $L(M,U)$ for a slightly other object --- a sheaf $\Cal F$ on a curve $X$ instead of a $\p[T]$-module $M$, see Section 16. This formula (a version of the Lefschetz trace formula) gives us $L(\Cal F,U)$ in terms of Frobenius action on $H^0(X, \Cal F)$, $H^1(X, \Cal F)$. A statement of this formula is given in [L], p. 2603, or in [GL16], (3.4). For a proof see [A00] (the original proof), or [B12], Section 9.
\medskip
As a corollary, we get that $L(M,U)\in \n F_q[T](U)$. Particularly, the order of zero of $L(M,U)$ at $U=1$ is defined. Is it a good analog of the analytic rank of abelian varieties? 
\medskip
{\bf 15.6.6.} As an example, we can apply this formula to twists of Carlitz modules. The corresponding theory is developed in [GL16], [GLZ22]. A generalization of this theory to other types of t-motives (for example, Drinfeld modules) is a research subject for beginners. It is much simpler than other research subjects related to Drinfeld modules.
\medskip
{\bf 15.6.7.} For a definiton of local $L$-factors at points of bad reduction, and for an analog of the Neron-Ogg-Shafarevich criterion see, for example, [Ga].
\medskip
{\bf 15.7. $v_\infty$ of eigenvalues of pure t-motives.} 
\medskip
We consider the valuation infinity $v_\infty$ on the field $\n F_q(T)$ defined by $v_\infty(T)=-1$, and we choose and fix its continuation to $\overline{\n F_q(T)}$.
\medskip
Let $M$ be a pure t-motive of rank $r$ and dimension $n$, defined over $\n F_{q^d}$. Its Euler factor $\det(I_r-Q^{[d]}U)$ belongs to $(\n F_q[T])[U]$, hence its roots $\la_i$ (inverse eigenvalues) belong to $\overline{\n F_q(T)}$. 
\medskip
{\bf Theorem 15.7.1.} $\forall \ i $ $v_\infty(\la_i)=-\frac{nd}{r}$ ([G96], Theorem 5.6.10, page 154). 
\medskip
This is an analog of (1.6.2) for the number field case. 
\medskip
{\bf Idea of the proof.} Let as earlier $S=T^{-1}$. We have: $\overline{\n F_q}((S))$ is a complete normed field with the valuation $v_\infty$, it is extended uniquely to $\overline{\n F_q((S))}$. Let (the simplest model case) $B\in GL_r(\ \overline{\n F_q[[S]]}\ )$ and $\mu_i\in \overline{\n F_q((S))}$ the roots of $\det(I_r-BU)=0$. Then $v_\infty(\mu_i)=0$. Really, the free term of $\det(I_r-BU)$ (considered as a polynomial in $U$) is 1, the leading coefficient is $\pm\det(B)\in  \overline{\n F_q[[S]]}^{*}$, and other coefficients belong to $\overline{\n F_q[[S]]}$. This means that the Newton polygon of $\det(I_r-BU)$ is the segment $[(0,0), \ (r,0)]$, hence the affirmation. 
\medskip
{\bf15.7.1.1. } Since $M$ is pure, we have $C\in GL_r(\ \overline{\n F_q((S))}\ )$ from (5.15.2). Let us consider a model example when $Q\in M_{r\times r}(\n F_q[T])$, $C\in GL_r(\ \n F_q((S))\ )$. In this case $Q^{[r]}=Q^r$, $Q^{[d]}=Q^d$, $C^{(r)}=C$ and (5.15.2) is 
$$ S^n\ C\ Q^r\ C^{-1}\in GL_r(\n F_q[[S]])$$ We get that in this case $v_\infty(\la_i)=-\frac{nd}{r}$. For the proof of the theorem we must reduce the case of general $Q$ and $C$ to this simple case. 
\medskip
{\bf Proof.}\footnotemark \footnotetext{ The authors are grateful to Urs Hartl for some help.} First, let us show that coefficients of the entries of $C$ belong to a finite field. The corresponding lemma is formulated entirely in terms of lattices in $\overline{\n F_q}((S))$, we do not need any action of $\tau$. 
\medskip
{\bf Lemma 15.7.1.2.} Let $V$ be the coordinate space $[\ \overline{\n F_q}((S))\ ]^r$ and $W\subset V$ a complete $\overline{\n F_q}[[S]]$-lattice. There exists $k$ such that there exists a basis of $W$ having a matrix $C\in \n F_{q^k}[[S]]$. 
\medskip
{\bf Proof of the lemma.} Let $W_0$ be the initial lattice $[\ \overline{\n F_q}[[S]]\ ]^r\subset V$. Multiplying $W$ by $S^n$ if necessary we can assume that $W\subset W_0$. There exists a sequence of lattices $$W_0\supset W_1\supset \dots \supset W_l=W$$ such that $\forall \ i$ we have $W_i/W_{i+1}=\overline{\n F_q}$ (this is a sequence of edges in the Bruhat-Tits building of $GL_r(\ \overline{\n F_q}((S))\ )$ joining the nodes corresponding to $W_0$ and $W$). Since $C$ is the product of matrices corresponding to $W_i\supset W_{i+1}$, it is sufficient to prove the lemma only for $W=W_1$. 
\medskip
We have $W_1\supset SW_0$, and $W_0/SW_0= \overline{\n F_q}^r$, $W_1/SW_0$ is a hyperplane in it. This hyperplane depends on finitely many parameters, hence $\exists \ k$ such that these parameters belong to $\n F_{q^k}$. This proves the lemma. 
\medskip
{\bf Remark.} We see that for a t-motive $M$ defined over $\n F_{q^d}$ we have a notion "pure over $\n F_{q^k}$". Maybe in reality all t-motives defined over $\n F_{q^d}$ are pure over $\n F_{q^d}$?

\medskip
So, let $C$ and $k$ be from the lemma. We denote $S^nC^{(r)}Q^{[r]}C^{-1}$ from (5.15.2) by $B$. We have $B\in GL_r(\ \n F_{q^k}[[S]]\ )$ and
$$B^{(kdr-r)}\cdot  B^{(kdr-2r)}\cdot  ... \cdot  B^{(2r)}\cdot  B^{(r)}\cdot  B= S^{nkd}C^{(kdr)}Q^{[kdr]}C^{-1}$$ $$=S^{nkd}CQ^{[kdr]}C^{-1}\in GL_r(\ \n F_{q^k}[[S]]\ ).$$
Let $\mu_i$ be the roots of $\det(I_r-Q^{[kdr]}U)=0$. Arguing as in (15.7.1.1) we get $v_\infty(\mu_i)=-nkd$. 
Further, we have $$Q^{[kdr]}=(Q^{[d]})^{(kdr-d)}\cdot  (Q^{[d]})^{(kdr-2d)}\cdot  ... \cdot  (Q^{[d]})^{(2d)}\cdot  (Q^{[d]})^{(d)}\cdot  Q^{[d]}=(Q^{[d]})^{kr}.$$ This means that $\mu_i=\la_i^{kr}$ hence the theorem. $\square$
\medskip
{\bf 15.8. Formula for $v_\g p(\la_i)$.} This is an analog of 1.6.3 for the function field case. Let us consider a numerical example of a Drinfeld module of rank $r$ over $\n F_{q^d}$. It comes from a map $$\iota: \n F_q[T]\to \n F_{q^d}, \ \ \ \ T\mapsto a_0$$ where $a_0\in \n F_{q^d}$ is $\th$. The characteristic is $\Ker \iota= < (T-a_0)(T-a_0^q)\dots (T-a_0^{q^{d-1}})>$. We denote it by $\g p$, its valuation by $v_\g p$ and its continuations to $\overline{\n F_q(T)}$ by $v_{\g p,i}$. 
\medskip
The matrix $Q$ is given in (5.12.1) ($\th$ is $a_0$). We have $\det(Q)=\pm (T-a_0)$, hence $\det(Q^{[d]})=\pm \g p$. Since $\det(Q^{[d]})$ is the leading term of $C_{2,M}(U)$, the rightmost vertex of the Newton polygon for \{$C_{2,M}(U)$, valuation $v_\g p$\} is $(r,1)$. The leftmost vertex is (0,0).
\medskip
All other vertices are $(i,0)$, $i=1,\dots, r-1$. To check: either this is true for all Drinfeld modules, or only for almost all? See [G96], Section 4. 
\medskip
The convex hull is \{$(0,0), \ (r-1,0), \ (r,1)$\}, the field (fields ? ) $\n F_q(T)(\la_i)$ has degree $r$ over $\n F_q(T)$ (is it Galois?) 
\medskip
We have: $v_{\g p,i}(\la_j)=\delta_j^i$, $i,j\in [1,\dots,r]$ (compare with (1.6.3)). See [G96] for details, and ? for the case of t-motives. I think that for (almost all ?) t-motives we have: $n$ eigenvalues $\la_i$ have $v_\g p=1$, remaining $r-n$ eigenvalues have $v_\g p=0$. 
\medskip
{\bf Analogy with abelian varieties with MIQF.} Let $A$ be an abelian variety with MIQF, of signature $(n,r-n)$, defined over a finite field. We can expect that (to continue).

\medskip
{\bf 16. Affine equations.}
\medskip
Sources: [GL21], [GL23]. We follow the exposition of [GL23]. 
\medskip
The ring $\p[[T]]$ is a left module over $\p[[T]]\{\tau\}$, where the multiplication by $\tau$ is defined by the formula 
$$\tau(\sum_{i=0}^\infty a_iT^i) :=(\sum_{i=0}^\infty a_iT^i)^{(1)} :=  \sum_{i=0}^\infty a_i^qT^i$$
An affine equation is an equation $$PX=0\eqno{(16.1)}$$ where $X=\sum_{i=0}^\infty x_iT^i\in \p[[T]]$ is an unknown and $$P=\sum_{\ga=0}^{r_0} a_\ga\tau^\ga + \sum_{\be=1}^\infty\ (\sum_{\ga=0}^{\vk_\be} b_{\be\ga} \tau^\ga ) \ T^\be\in\p[[T]]\{\tau\}\eqno{(16.2)}$$ is a coefficient, see [GL21], (2.2.2). 
\medskip
Also, we can consider a system of $\la$ affine equations in $\mu$ variables 
$$\g P\g X=0\eqno{(16.3)}$$
where $\g P=\{P_{ij}\}\in M_{\la\times\mu}\p[[T]]\{\tau\}$ is a matrix and $\g X=\{X_{i}\}\in M_{\mu\times1}\p[[T]]$ is a matrix column. 
\medskip
{\bf Example 16.4.} Equations (12.3.1.3), (12.3.2.4) defining $H_1(M), \ H^1(M)$ are systems of affine equations, $\la=\mu=r$: for (12.3.1.3), respectively (12.3.2.4) we have 
$$\g P=I_r \tau-Q, \hbox{ respectively $\g P=Q^t \tau-I_r$ (for } \g X_*=Y_*^t). $$

{\bf Example 16.5.} For any monic irreducible polynomial $\g L$ in $\n F_q[T]$, the equations for the $T_\g L(M)$ --- the $\g L$-Tate module of $M$, see (6.13), are a system of affine equations. 
\medskip
Let us consider the case $\g L=T$. We have first, the system of affine equations defining $H_1(M)$, it has $\la=\mu=r$, and the system of affine equations defining $T_T(M)$, it has $\la=\mu=n$. We have
\medskip
{\bf Proposition 16.5.1.} For any $M$ non-commutaive determinants (see 16.9) of the above two systems coincide. 
\medskip
For a Drinfeld module $M$ a proof by means of explicit calculations is given in [GL23]. For a proof in the general case see [NP]. (to give more details. Non-commutaive determinants along what column or linear form?)
\medskip
The explicit form of equations (16.1) is the following (we use coefficients of $P$ from (16.2)). The equality of the $T$-free term of (16.1) is 
$$a_rx_0^{q^r}+a_{r-1}x_0^{q^{r-1}}+...+a_1x_0^q+a_0x_0=0\eqno{(16.6)}$$

The set of $x_0$ satisfying this equation is a vector space over $\n F_q$ of dimension $r$ contained in $\p$. We denote it by $W_0$.
\medskip
The equality of the coefficient at $T$ in (16.1) is 
$$a_rx_1^{q^r}+a_{r-1}x_1^{q^{r-1}}+...+a_1x_1^q+a_0x_1+\eqno{(16.7h)}$$ $$ +b_{1,\vk_1}x_0^{q^{\vk_1}}+ b_{1,\vk_1-1}x_0^{q^{\vk_1-1}}+...+b_{11}x_0^q+b_{10}x_0=0\eqno{(16.7t)}$$

If we fix $x_0\in W_0$ then the value of the left hand side of (16.7t) is fixed, and the set of $x_1$ which are roots of (16.7h) is an affine space over $W_0$, this justifies the terminology "affine equation". The part (16.7h) of (16.7) is called the head of the equation, and the part (16.7t) is called its tail, hence the notations. 

The equality of the coefficient at $T^2$ in (16.1) is 
$$a_rx_2^{q^r}+a_{r-1}x_2^{q^{r-1}}+...+a_1x_2^q+a_0x_2+\eqno{(16.8h)}$$ $$ +b_{1,\vk_1}x_1^{q^{\vk_1}}+ b_{1,\vk_1-1}x_1^{q^{\vk_1-1}}+...+b_{11}x_1^q+b_{10}x_1+\eqno{(16.8t.1)}$$ $$ +b_{2,\vk_2}x_0^{q^{\vk_2}}+ b_{2,\vk_2-1}x_0^{q^{\vk_2-1}}+...+b_{21}x_0^q+b_{20}x_0=0\eqno{(16.8t.2)}$$

The parts (16.8t.1), (16.8t.2) are called respectively the first and second parts of the tail of (16.8). 
As above, if we fix $x_0, \ x_1$ satisfying (16.7) then the set of $x_2$ which are roots of (16.8h) is also an affine space over $W_0$. 
\medskip
And so on. We distinguish two cases: if $P\in \p[T]\{\tau\}$ is a polynomial in $T$ of degree say $n$, then the length of the tail ( = the quantity of the tail parts) stabilizes at $n$. These equations are called the equations of bounded tail. For example, equations obtained from (12.3.1.3), (12.3.2.4), and from $T_\g L(M)$ after elimination of unknowns, see 16.9, are of bounded tail. 
\medskip
If $P\in \p[[T]]\{\tau\}$ is a power series in $T$, then the length of the tail grows to infinity. 
\medskip
Let is give properties of the sets of solutions of (the systems of) affine equations. First, we can reduce the case of a system of affine equations to the case of one equations.  
\medskip
{\bf 16.9. Elimination of variables.} Let us consider a system (16.3) where $\la=\mu=r$. We can eliminate unknowns $X_2, \dots, X_r$ getting one affine equation $$det_{1,c}(\g P)\cdot X_1=0$$
where $det_{1,c}(\g P)\in  \p[[T]]\{\tau\}$ is called "the non-commutative determinant of $\g P$ along the first column". See [GL23], Proposition 1.5 for a proof. 
\medskip
%\newpage
{\bf 16.10. Set of solutions of an affine equation.} 
\medskip
{\bf Definition 16.10.1.} Let $\g W$ be a $\n F_q[[T]]$-submodule of $\p[[T]]$. $\g W$ is called a $T$-divisible submodule if for any $y\in \g W\cap T\cdot \p[[T]]$ we have $y/T\in \g W$.
\medskip
{\bf Proposition 16.10.2.} ([GL21], Proposition 2.3.2). The set of solutions of (16.1) is a $T$-divisible submodule of $\p[[T]]$. Conversely, for any $T$-divisible submodule $W$ of $\p[[T]]$ of dimension $r_0$ there exists the only affine equation (16.1) we have (notations of (16.2)):
$$r=r_0, \ \ \ a_r=1, \ \ \ a_0\ne0, \ \ \ \forall \ \be\hbox{ we have }\vk_\be< r\eqno{(16.10.3)}$$
such that $W$ is the set of its roots. 
\medskip
{\bf Remark.} It would be interesting to find a criterion in terms of $W$: whether its affine equation $PX=0$ is of finite tail, or not? 
\medskip
{\bf 16.11. Small solutions.} 
\medskip
{\bf Definition 16.11.1.} A solution of (16.1), respectively (16.3) is called small if it belongs to $\p\{T\}$, respectively to $\p\{T\}^\mu$.
\medskip
Since we apply the theory of affine equations to solutions of (12.3.1.3), (12.3.2.4), we are interested in small solutions. The set of small solutions is denoted by $W_s$, it is a $\n F_q[T]$-module. We have
\medskip
{\bf Theorem 16.11.2.} Let $f_1, \dots,f_n\in W_s$ be linearly independent over $\n F_q[T]$. Then they are linearly independent over $\n F_q[[T]]$.
\medskip
For a proof see [GL21], Proposition 2.3.4 (the proof is not immediate). As a corollary of this theorem, we get that $h^1, \ h_1\le r$. 
\medskip
{\bf 16.12. Holonomic sequences. } The analogy 5.20 between the Anderson ring and a ring of differential operators permits us to define the classical notion of holonomic sequences for the function field case. Really, $F$ from (5.20.2) is an analog of $P\in\p[T,\tau]$, and a differential equation $F(g)=0$ is an analog of (16.1). 
\medskip
Let us recall the definitions and properties of holonomic sequences in $\n C$, see [K], [KP]. We follow [GL23]:
\medskip
{\bf Definition 16.12.1.} A sequence $x_0, \ x_1, \dots$, where $x_i\in \n C$, is called a (classical) holonomic sequence, if two equivalent conditions hold: 
\medskip
{\bf 16.12.1.1.} $\exists \ r_1>0$ and polynomials $p_0(x),\dots,p_{r_1}(x)\ne0$ such that $\forall \ k\ge0$ we have
$$p_0(k)a_k+p_1(k)a_{k+1}+...+p_{r_1}(k)a_{k+r_1}=0$$

{\bf 16.12.1.2.} $\exists \ r_2 >0$ and polynomials $q_0(T),\dots,q_{r_2}(T)\ne0$ such that $f(T):=\underset{i=0}\to{\overset{\infty}\to{\sum}}x_iT^i\in \n C[[T]]$ satisfies $P(f)=0$, where $P=\underset{i=0}\to{\overset{r_2}\to{\sum}}q_iD^i$ is a differential operator, i.e. $$q_0f+q_1f'+...+q_{r_2}f^{(r_2)}=0$$

Let $n_1$, resp. $n_2$ be the maximal degree of $p_i(x)$, resp. $q_i(T)$. We have ([KP], Th. 7.1):
\medskip
If (16.12.1.1) holds for $x_0, \ x_1, \dots$ with some $r_1$, $n_1$, then (16.12.1.2) holds for $x_0, \ x_1, \dots$ with $r_2\le n_1$, $n_2\le r_1+n_1$;
\medskip
If (16.12.1.2) holds for $x_0, \ x_1, \dots$ with some $r_2$, $n_2$, then (16.12.1.1) holds for $x_0, \ x_1, \dots$ with $r_1\le r_2+ n_2$, $n_1\le r_2$. 
\medskip
{\bf 16.12.2. Properties of holonomic sequences.} We have a theorem for classical holonomic sequences: 
\medskip
Let $x_0, \ x_1, \dots$ and $y_0, \ y _1, \dots$ be holonomic sequences. Then
\medskip
1. $x_0+y_0, \ x_1+y_1, \dots$;
\medskip
2. $x_0y_0, \ x_1y_1, \dots$;
\medskip
3. $x_0y_0, \ x_0y_1+x_1y_0, \dots, \underset{i=0}\to{\overset{n}\to{\sum}}x_iy_{n-i}, \dots$
\medskip
are holonomic sequences. 
\medskip
See [K], Theorem 4 for a statement of these properties and for some other similar properties. 
\medskip
{\bf Definition 16.12.3.} A sequence $x_0, \ x_1, \dots$, where $x_i\in \p$, is called a (characteristic $p$) holonomic sequence, if $\exists \  P\in \p[T,\tau]$ such that $X:=\underset{i=0}\to{\overset{\infty}\to{\sum}}x_iT^i$ is a solution of the affine equation $PX=0$. 
\medskip
Here it is important that $P\in \p[T,\tau]$ (case of finite tail) and not $P\in \p[[T]][\tau]$, because $\forall \ X \ \exists \  P\in \p[[T]][\tau]$ such that $PX=0$. 
\medskip
We can expect that analogs of the properties 16.12.2, (1) -- (3), and of many other properties of classical holonomic sequences, hold for the characteristic $p$ holonomic sequences. A particular case of 16.12.2.1 is proved in [GL23], Section 6.
\medskip
There exist generalizations of classical holonomic sequences for high-dimensional cases, when we consider operators of partial derivation. It would be desirable to extend these definitions for the case of characteristic $p$ holonomic sequences. 
\medskip
{\bf 16.13. Ideal annulator. } Let $X=\underset{i=0}\to{\overset{\infty}\to{\sum}}x_iT^i$ be a characteristic $p$ holonomic sequence. The set of $P\in \p[T,\tau]$ such that $PX=0$ is a left ideal in $\p[T,\tau]$ denoted by $Ann(X)$. It is not necessarily principal, see [GL23], end of (6.3). 
\medskip
What are possible degrees in $T$ and $\tau$ of generators of $Ann(X)$? [KP], Fig. 7.1, p. 143 gives us an example of these degrees for a classical holonomic sequence. 
\medskip
We can mention that for a holonomic sequence $X$ satisfied $PX=0$ for $$P=\tau^2+a\tau+\th+\tau T$$ for the only $P_0\in \p[[T]][\tau]$ satisfying (16.10.3) such that $P_0X=0$ we have: $P_0\not\in\p[T,\tau]$ ([GL23], Example 6.4).
\medskip
{\bf 16.14. $W$ as Galois modules.} Let coefficients of $P$ or of entries of $\g P$ of (16.1), (16.3) belong to $\n F_q(\th)$. Then $\Gal(\n F_q(\th))$ acts on $W$. Hence, $W$ is an object of the category Rep$_{A_\pi}(G_L)$ (see [G95], Appendix) where $A_\pi=\n F_q[[T]]$, and $L= \n F_q(\th)$.
\medskip
According [G95], Proposition, p. 184, attached to $W$ is a $\vf$-module over $\n F_q(\th)[[T]]$ denoted by $D(W)$. If $\g P$ comes from (12.3.1.3) or (12.3.2.4) of a t-motive $M$ over $\n F_q(\th)$ or if $P$ is its non-commutative determinant, see 16.9, then $D(W)$ should be related with $M$. 
\medskip
What is $D(W)$ if $P$ does not come from any $M$?
\medskip
{\bf 17. The Drinfeld "upper half plane".}
\medskip
The characteristic $p$ analog of the upper half plane $\Cal H$ is the Drinfeld upper half plane $\oo:=\p-\r$. We start giving a detailed description of its covering spaces $D_{(n,x)}$ ([GR], p. 33 - 35). 

We use $\pi:=\th^{-1}$ as an uniformizer, i.e. $v_\infty(\pi)=1$. Sets $D_{(n,x)}\subset \oo$ of [GR] will be denoted by $D(n,x)$. We need to show that $\oo=\cup D(n,x)$, and to find intersections of $D(n,x)$. 
\medskip
{\bf Definition 17.1.} Let $n\in \n Z$. The subset $D(n)\subset \oo$ is the union of 3 sets 

$$z\in \oo \ | \ n<v_\infty(z)<n+1\eqno{(17.2)}$$
$$z\in \oo  \ | \  \vi(z)=n \hbox { and $ \forall \ c\in \n F_q $ we have } \vi(z-c\pi^n)=n\eqno{(17.3)}$$
$$z\in \oo  \ | \  \vi(z)=n+1 \hbox { and $ \forall \ c\in \n F_q $ we have } \vi(z-c\pi^{n+1})=n+1.\eqno{(17.4)}$$
The conditions (17.3), (17.4) mean that if $\vi(z)$ is a integer $m$, where $m=n$ or $n+1$, then $\forall \ c\in \n F_q $ we have: $z$ is "sufficiently far" from $c\pi^m$: if it is close to $c\pi^m$ then $\vi(z-c\pi^m)>m$.
\medskip
The geometric interpretation of the reduction of $D(n)$ is the following. Let $z$ satisfy (17.2) or (17.3). Since $\vi(z\pi^{-n})\ge0$, the reduction $\widetilde{z\pi^{-n}}\in \overline{\n F_q}\subset P^1(\overline{\n F_q})$ is defined. If $z$ satisfies (17.2) then $\widetilde{z\pi^{-n}}=0$; if $z$ satisfies (17.2) then $$\widetilde{z\pi^{-n}}\in P^1(\overline{\n F_q})-\infty-(\n F_q-0)=P^1(\overline{\n F_q})-P^1(\n F_q)\cup0.\eqno{(17.5)}$$

Let now $z$ satisfy (17.2) or (17.4). The reduction $\widetilde{z\pi^{-(n+1)}}\in P^1(\overline{\n F_q})$: if $z$ satisfies (17.2) then $\widetilde{z\pi^{-(n+1)}}=\infty\in P^1(\overline{\n F_q})$. If $z$ satisfies (17.4) then $$\widetilde{z\pi^{-(n+1)}}\in P^1(\overline{\n F_q})-\n F_q=P^1(\overline{\n F_q})-P^1(\n F_q)\cup\infty.\eqno{(17.6)}$$

We see that $\widetilde{D(n)}$ is the union of the sets (17.5) and (17.6) glued by the points 0 from (17.5) and $\infty$ from (17.6). 
\medskip
Let us define $D(n,x)$. We fix $n\in \n Z$ and we consider the set of $x\in F_q(\th)$ such that 
$$x=c_n\pi^n+c_{n-1}\pi^{n-1}+ ... +c_{n-k}\pi^{n-k}\eqno{(17.7)}$$ for some $k$, where $c_{*}\in \n F_q$. We see that $\vi(x)\le n$. By definition ($x$ is from 17.7),
$$D(n,x)=z\in\p \  |\ z-x \in D(n).$$
We have: if $x\ne0$ then
$$\forall \ z\in D(n,x)\hbox { we have }\vi(z)=\vi(x)=n-k$$
This is clear if $k>0$; it is easy to see that this is also true for $k=0$. 
\medskip
Let $z\in \oo$. Let us find all $(n,x)$ such that $z\in D(n,x)$. 
\medskip
First, if $\vi(z)\not\in\n Z$, namely, $n<\vi(z)<n+1$, then $z\in D(n,0)$ and to no other $D(k,x)$. 
\medskip
If $\vi(z)=n$ and $\forall \ c_n\in\n F_q$ we have $\vi(z-c_n\pi^n)=n$, then $z\in D(n-1,0)$, $z\in D(n,c_n\pi^n)$ ($\forall \ c_n\in\n F_q$), and to no other $D(k,x)$. 
\medskip
If $\exists \ c_n\in\n F_q$ such that $\vi(z-c_n\pi^n)>n$, then we denote $z_1:=z-c_n\pi^n$ and we apply the above procedure to $z_1$. 
\medskip
Namely, if $\vi(z_1)\not\in\n Z$: $n_1<\vi(z_1)<n_1+1$, where $n_1\ge n$, then $z\in D(n_1,\ c_n\pi^n)$ and to no other $D(k,x)$. 
\medskip
If $\vi(z_1)\in\n Z$, namely, $\vi(z_1)=n_1$ (here $n_1>n$) and 
\medskip
$\forall \ c_{n_1}\in\n F_q$ we have $\vi(z_1-c_{n_1}\pi^{n_1})=n_1$, 
\medskip
then $z\in D(n_1-1,\ c_n\pi^n)$, $z\in D(n_1,\ c_n\pi^n+d_{n_1}\cdot \pi^{n_1})$ ($\forall \ d_{n_1}\in\n F_q$), 
and to no other $D(k,x)$. 
\medskip
If $\exists \ c_{n_1}\in\n F_q$ such that $\vi(z_1-c_{n_1}\pi^{n_1})>n_1$, then we denote $z_2:=z_1-c_{n_1}\pi^{n_1}$ and we apply the above procedure to $z_2$, etc.
\medskip
This procedure is finite: if we have an infinite sequence $z-\sum_i  c_{n_i}\pi^{n_i}$ that tends to 0, then this means that $z\in \n R_\infty$. Hence, two possibilities of the end of this procedure can occur. Namely, we shall come to 
$$z_m:=z-c_{n}\pi^{n}-c_{n_1}\pi^{n_1}-...-c_{n_{m-1}}\pi^{n_{m-1}}$$
First possibility: $\vi(z_m)\not\in \n Z$:
$$n_m<\vi(z_m)<n_m+1$$
In this case $$z\in D(n_m,\ c_n\pi^n+c_{n_1}\pi^{n_1}+...+c_{n_{m-1}}\pi^{n_{m-1}})$$ and to no other $D(k,x)$. 

Second possibility: $\vi(z_m)=n_m\in \n Z$, and $\forall \ c_{n_m}\in\n F_q$ we have $$\vi(z_m-c_{n_m}\pi^{n_m})=n_m.$$ 
\medskip
In this case $$z\in D(n_m-1,\ c_n\pi^n+c_{n_1}\pi^{n_1}+...+c_{n_{m-1}}\pi^{n_{m-1}});$$
$$z\in D(n_m,\ c_n\pi^n+c_{n_1}\pi^{n_1}+...+c_{n_{m-1}}\pi^{n_{m-1}}+d_{n_{m}}\pi^{n_{m}})$$ ($\forall \ d_{n_m}\in\n F_q$)
and to no other $D(k,x)$. 
\medskip
{\bf 17.8. Sets $D(n,x)$ and the Bruhat-Tits building $\Cal T$ of $\oo$.}
\medskip
Let us consider the set of all $D(n,x)$, where $n\in \n Z$ and $x$ is from 17.7. They are in 1 -- 1 correspondence with the edges of the Bruhat-Tits building $\Cal T$ of $\oo$, and the vertices attached to an edge are two irreducible components (17.5) and (17.6). 
\medskip
It is easy to see that if we fix one of these irreducible components (i.e. a vertice) then there are exactly $q+1$ edges adjacent to it, they correspond to elements of $P^1(\n F_q)$. The reduction map $\oo \to \Cal T(\n Q)$ coincides with the map defined by norms. 
\medskip

{\bf 18. Relations between $\vf$-modules and Galois modules.}
\medskip
Here we explain the contents of [G95], Appendix (based, in its turn, on [Tag96]). There is a simple fact: for a Drinfeld module $M$ its $H_1(M)$ essentially coinsides with its $T_T(M)$ --- the $T$-Tate module of $M$. We give here an explicit calculation that an analog of this coincidence holds for any prime ideal $\pi$ of $\n F_q[T]$, not only for the prime ideal $T\cdot \n F_q[T]$, i.e. that $T_\pi(M)$ --- the $\pi$-Tate module of $M$ --- essentially coincides with a $\pi$-analog of $H_1(M)$. 
\medskip
Therefore, we fix an irreducible polynomial = a prime ideal $\pi$ of $\n F_q[T]$ of degree $d$. We consider $\n F_q[T]_\pi$ --- the completion of $\n F_q[T]$ at $\pi$. It is the set of series 
$$A_0+A_1\pi+A_2\pi^2+...\eqno{(18.1)}$$
where $A_i\in \n F_q[T]$ are polynomials in $T$ of degree $<d$ --- representatives of $\n F_q[T]/\pi\n F_q[T]$. 
\medskip
Further, let $L$ be a (finite) extension of $\n F_q(\th)$. We consider $\n F_q[T]_\pi\underset{\n F_q}\to{\otimes}L$ (denoted by {\bf A}$_{L,\pi}$ in [G95]); it is a semi-local ring. A $\vf$-module over $\n F_q[T]_\pi\underset{\n F_q}\to{\otimes}L$ is a module $M$ over $\n F_q[T]_\pi\underset{\n F_q}\to{\otimes}L$ and a skew linear map $\vf: M \to M$, where skew linear means that $$\vf[(a\otimes l)(m)]=(a\otimes l^q)(\vf(m))$$ ($a\in \n F_q[T]_\pi, \ l\in L$).
\medskip
{\bf Example 18.2.} Let $M$ be an Anderson t-motive over $L$. For any $\pi$ it can be extended to a $\vf$-module over $\n F_q[T]_\pi\underset{\n F_q}\to{\otimes}L$ ($\vf$ is $\tau$), we denote it by $M(\pi)$. For example, if $\pi=T$ and $L=\p$ then $$M(T)=M\otimes_{\p[T]}\p[[T]]$$
(compare with $M\{T\}$ from (12.3.2.1)). 
\medskip
There is a functor $V$ from the category of $\vf$-modules over $\n F_q[T]_\pi\underset{\n F_q}\to{\otimes}L$ to the category of $\n F_q[T]_\pi$-modules with $\Gal(L)$-action. It is defined as follows. Let $\bar L$ be an algebraic closure of $L$.  We consider  $$\n F_q[T]_\pi\otimes\bar L:=(\n F_q[T]_\pi\underset{\n F_q}\to{\otimes}L)\underset{L}\to{\otimes}\bar L=\n F_q[T]_\pi\underset{\n F_q}\to{\otimes}\bar L$$ and $\bar M:=M\underset{L}\to{\otimes}\bar L$. The action of $\vf$ can be skew-continued to $\bar M$. By definition, $V(M)$ is $\bar M^\vf$ --- the set of $\vf$-invariant elements. It is obviously a $\n F_q[T]_\pi$-module, and $\Gal(L)$ acts on it by the action on $\bar L$. 
\medskip
{\bf Example 18.3.} Let as above $M$ be an Anderson t-motive over $L$, $\pi=T$ and $L=\p$. If $M$ is uniformizable then we have $V(M(\pi))=H^1(M)\underset{\n F_q[T]}\to{\otimes}\n F_q[[T]]$. If $M$ is non-uniformizable then we have $V(M(\pi))\supset H^1(M)\underset{\n F_q[T]}\to{\otimes}\n F_q[[T]]$.
\medskip
(maybe $=H_1(M)$? Few likely, because both $V(M(\pi)), \ H^1(M)$ are covariant with respect to $M$, and elements of $H^1(M)$ are exactly $\tau$-invariant elements of $M\{T\}$). 
\medskip
We refer to [G95], p. 182 - 184 for the definition of an \'etale $\vf$-module. 
\medskip
Let us recall the definition of the inverse functor $D$ from the category of $\n F_q[T]_\pi$-modules with $\Gal(L)$-action to the category of $\vf$-modules over $\n F_q[T]_\pi\underset{\n F_q}\to{\otimes}L$. Namely, let $W$ be a $\n F_q[T]_\pi$-module with a $\Gal(L)$-action. We consider $\bar L[T]_\pi\underset{?}\to{\otimes}W$ (utochnit', nad chem). (Raznica mezhdu $\n F_q[T]_\pi\otimes \bar L$ i $\bar L[T]_\pi$). The action of $\Gal(L)$ is defined as follows: $g(l\otimes w):=g(l)\otimes g(w)$. The set $D(W)$ is the set of $\Gal(L)$-invariant elements of  $\bar L[T]_\pi\underset{?}\to{\otimes}W$. It is a $\n F_q[T]_\pi$-module. The action of $\tau$ is given by the formula $\tau(P\otimes w)=P^{(1)}\otimes w$, here $P\in \bar L[T], \ w\in W$. 
\medskip
{\bf Example 18.4.} Let  $W$ be a $\n F_q[T]_\pi$-module with a $\Gal(L)$-action. Then $V(D(W))=W$. 
\medskip
{\bf Proof.} Let $h_*=(h_1,\dots,h_r)^t$ be a basis of $W$ over $\n F_q[T]_\pi$. The action of $\Gal(L)$ on $W$ is given by a map $\ga: \Gal(L)\to GL_r(\n F_q[T]_\pi)$. 
\medskip
It is known (reference!) that $D(W)$ has a basis over $\n F_q[T]_\pi\underset{\n F_q}\to{\otimes}L$ consisting of $r$ elements $Z_1, \dots, Z_r\in \bar L[T]_\pi\underset{?}\to{\otimes}W$. We can denote $Z_i=\sum_{j=1}^r Z_{ij}\otimes h_j$. 
\medskip
Further, we denote $Z=\left(\matrix Z_1\\ \dots \\ Z_r\endmatrix \right)$ (a matrix column). We denote the $r\times r$ matrix $\{Z_{ij}\}$ by $Z_m$ (the subscript $m$ means matrix). The fact that all $Z_i$ are $\Gal(L)$-invariant means that $$\forall \ g\in \Gal(L) \hbox{ we have } g(Z_m)\ga(g)=Z_m\eqno{(18.4.1)}$$ 

The action of $\tau$ on $Z$ is given by a matrix $Q\in M_{r\times r}(\n F_q[T]_\pi\underset{\n F_q}\to{\otimes}L)$ satisfying $QZ_m=Z_m^{(1)}$. 
\medskip
Let us consider the functor $V$ applied to $D(W)$. There exists a matrix $Y\in  \bar L[T]_\pi\underset{?}\to{\otimes}W$ having the property: its lines $Y_i:=(Y_{i1}, \dots, Y_{ir})$ give us $\tau$-invariant elements of $D(W)\underset{?}\to{\otimes}\bar L[T]$. More exactly, this means that $\forall \ i $ we have: 
$\sum_j Y_{ij}Z_j$ is $\tau$-invariant. 
\medskip
This is equivalent to $$Y^{(1)}Q=Y,$$
i.e. $$Y=Z_m^{-1}\eqno{(18.4.2)}$$ (why $Z_m$ is invertible?) 
\medskip
Now, we need to consider the action of $g\in \Gal(L)$ on $Y_i$. Let $\de=\de(g)$ be a matrix describing this action: $g(Y_i)=\sum_j \de_{ij}Y_j$. This is equivalent to $$g(Y)=\de Y\eqno{(18.4.3)}$$ 
\medskip
Comparing (18.4.1) and (18.4.3), and taking into consideration (18.4.2), we get that $\de(g)=\ga(g)$. This means that $V(D(W))=W$. $\square$
\medskip
{\bf Theorem 18.8.} Let $M$ be a Drinfeld module. The $\n F_q[T]_\pi$-module with $\Gal(L)$-action $V(M(\pi))^*$ is isomorphic to $T_\pi(M)$.
\medskip
{\bf Proof.} This is [G95], p. 186, we give here explicit calculations. Let $\pi=T^d+c_{d-1}T^{d-1}+...+c_0\in \n F_q[T]$ and $M$ be defined by (5.6), $a_0=\th$, $a_r=1$:
$$Te=(\th+a_1\tau+a_2\tau^2+\dots +a_{r-1}\tau^{r-1}+\tau^r)e$$
A basis of $\bar M$ over $\n F_q[T]_\pi\otimes\bar L$ is $f_*=(e, \tau e, ...,\tau^{r-1}e)^t$. The matrix of the action of $\tau=\vf$ in this basis is $Q$ from (5.12.1). We use notations of (12.3.1): $X_*$ is the matrix column of coefficients of a $\tau$-stable element of $\bar M$ in this basis, and (12.3.1.3): $QX_*=X_*^{(1)}$ is the condition of $\tau$-stability. We write 
$$X_*=(X_0,\dots,X_{r-1})^t\in (\n F_q[T]_\pi\otimes\bar L)^r, $$ $$ X_i=\sum_{j=0}^\infty X_{ij}\pi^j\in \n F_q[T]_\pi\otimes\bar L,\ \ \ \ X_{ij}=\sum_{k=0}^{d-1}x_{ijk}T^k,\ \ \ \ x_{ijk}\in \bar L,$$
like in (18.1). Condition (12.3.1.3) is $$X_0^{(1)}=X_1, \  \  \ X_1^{(1)}=X_2, \dots, X_{(r-2)}^{(1)}=X_{r-1}, $$ $$ X_{(r-1)}^{(1)}=(T-\th)X_0-a_1X_1-...-a_{r-1}X_{r-1},$$  i.e. $$X_0^{(r)}=(T-\th)X_0-a_1X_0^{(1)}-...-a_{r-1}X_0^{(r-1)}\eqno{(18.8.1)}$$
We have $$X_0^{(\al)}=\sum_{j=0}^\infty \ ( \ \sum_{k=0}^{d-1}  x_{0jk}^{q^\al}T^k\ )\ \pi^j.$$ Let us write first the explicit form of (18.8.1) for the $\pi$-free terms. We have: 
$$\th X_0=\th x_{000}+\th x_{001}T+...+\th x_{00,d-1}T^{d-1}\hbox { (only $\pi$-free terms);}\eqno{(18.8.2a)}$$
$$a_1 X_0^{(1)}=a_1 x_{000}^q+a_1 x_{001}^qT+...+a_1 x_{00,d-1}^qT^{d-1}\hbox { (only $\pi$-free terms);}$$
$$\dots$$
$$a_{r-1} X_0^{(r-1)}=a_{r-1} x_{000}^{q^{r-1}}+a_{r-1} x_{001}^{q^{r-1}}T+...+a_{r-1} x_{00,d-1}^{q^{r-1}}T^{d-1}\hbox { (only $\pi$-free terms);}$$
$$X_0^{(r)}= x_{000}^{q^{r}}+ x_{001}^{q^{r}}T+...+ x_{00,d-1}^{q^{r}}T^{d-1}\hbox { (only $\pi$-free terms);}$$
$$TX_0=x_{000}T+ x_{001}T^2+...+x_{00,d-2}T^{d-1}-$$ $$ -c_0x_{00,d-1}-c_1x_{00,d-1}T-...-c_{d-1}x_{00,d-1}T^{d-1}+ $$ $$+x_{00,d-1}\pi+\hbox{ other $\pi$-terms}\eqno{(18.8.2b)}$$
(because $T^d=-c_0-c_1T-...-c_{d-1}T^{d-1}+\pi$). 
For $x\in \p$ we denote $$T_E(x):=\th x +a_1 x^q+...+a_{r-1}x^{q^{r-1}}+x^{q^r}$$ (the action of $T$ on $E(M)$, see (6.0)). Equality (18.8.1) becomes (taking into consideration (18.8.2.*))

$$T_E(x_{00,d-1})-x_{00,d-2}+c_{d-1}x_{00,d-1}=0 \eqno{(18.8.3a)}$$ (equality of coefficients at $T^{d-1}$ of the $\pi$-free term in (18.8.1));
$$T_E(x_{00,d-2})-x_{00,d-3}+c_{d-2}x_{00,d-1}=0 \eqno{(18.8.3b)}$$ (equality of coefficients at $T^{d-2}$ of the $\pi$-free term in (18.8.1)); etc; 
$$\dots$$
$$T_E(x_{001})-x_{000}+c_{1}x_{00,d-1}=0 \eqno{(18.8.3c)}$$ (equality of coefficients at $T$ of the $\pi$-free term in (18.8.1)); 
$$T_E(x_{000})+c_{0}x_{00,d-1}=0 \eqno{(18.8.3d)}$$ (equality of the free term in (18.8.1)).

(18.8.3a) implies that $$x_{00,d-2}=(T_E+c_{d-1})(x_{00,d-1}). \eqno{(18.8.3.1)}$$ Further, (18.8.3b) implies that $$x_{00,d-3}=T_E(x_{00,d-2})+c_{d-2}x_{00,d-1}=(T_e^2+c_{d-1}T_E+c_{d-2})(x_{00,d-1}). \eqno{(18.8.3.2)}$$
Continuing we get $$T_\pi(x_{00,d-1})=0.$$

{\bf 18.8.4.} Let us consider now equalities of coefficients at $\pi^j$-terms of (18.8.1), $j\ge1$. First, for $j=1$ we have the same equations (18.8.3*) where $x_{00k}$ is replaced by $x_{01k}$, and the term $x_{00,d-1}\pi$ of (18.8.2b) will be added with the coefficient $-1$ to (18.8.4d): 
$$T_E(x_{01,d-1})-x_{01,d-2}+c_{d-1}x_{01,d-1}=0 \eqno{(18.8.4a)}$$ (equality of coefficients at $T^{d-1}\pi$ in (18.8.1));
$$T_E(x_{01,d-2})-x_{01,d-3}+c_{d-2}x_{01,d-1}=0 \eqno{(18.8.4b)}$$ (equality of coefficients at $T^{d-2}\pi$ in (18.8.1)); etc; 
$$\dots$$
$$T_E(x_{011})-x_{010}+c_{1}x_{01,d-1}=0 \eqno{(18.8.4c)}$$ (equality of coefficients at $T\pi$ in (18.8.1)); 
$$T_E(x_{010})+c_{0}x_{010}-x_{00,d-1}=0 \eqno{(18.8.4d)}$$ (equality of coefficients at $\pi$ in (18.8.1)).
Acting like in (18.8.3.1), (18.8.3.2) etc. we get $$T_\pi(x_{01,d-1})=x_{00,d-1}. \eqno{(18.8.4.1)}$$
Further, the analog of the term $x_{00,d-1}\pi$ of (18.8.2b) will be the term $x_{01,d-1}\pi^2$, hence for $j=2$ the analogs of (18.8.4d), (18.8.4.1) are $$T_E(x_{020})+c_{0}x_{020}-x_{01,d-1}=0.$$ (equality of coefficients at $\pi^2$ in (18.8.1)), and
$$T_\pi(x_{02,d-1})=x_{01,d-1}. $$ Continuing we get that for any $j$ $$T_\pi(x_{0,j+1,d-1})=x_{0,j,d-1} $$hence the elements $x_{0,j,d-1}$ ($j=0,1,\dots$) are an element of $T_\pi(E(M))$. $\square$
\medskip
{\bf 18.9.} Hence, we have an identification $\Hom_{\p[T]}(M,\n F_q[T]_\g p\otimes \p)^\tau$ and $T_\g p(M)$. If $\g p=T$ then it is reasonable to consider small elements in $\Hom_{\p[T]}(M,\n F_q[T]_T\otimes \p)^\tau=\Hom_{\p[T]}(M,\p\{T\})^\tau$. Their dimension is $h_1(M)$ --- an important invariant of $M$. We have a natural question: what the dimension of the space of small elements in $T_\g p(M)$?
\medskip
Answer: it is the same. 
\medskip
{\bf 19. Flag varieties related to lattices of t-motives.}
\medskip
Let us describe the contents of this chapter. We work exclusively with lattices in $\p^n$, we do not consider t-motives at all. We describe here an analog of a Siegel matrix of an abelian variety with MIQF. Recall that a Siegel matrix of an abelian variety with MIQF is a $(r-n)\times n$ complex matrix satisfying (1.8.1). For the function field case, its analog for a lattice is an element of the Grassmannian $Gr(n,r)$ (case $N=0$) or of a kind of a flag variety (denoted by tFV) for the case $N\ne0$. We describe the Schubert cells containing these analogs of Siegel matrices, we define modular forms, we show how to apply this construction to the properties of duality of lattices. In future, this construction can be applied to the notion of tensor product of lattices and of t-motives. 
\medskip
Results of Sections 19.6.1 -- 19.6.6 are borrowed from [GL18], results of 19.2 -- 19.5, and of 19.6.7 are new. 
\medskip
{\bf 19.1. Definitions related to Grassmannians.} For a Grassmannian $Gr(d,n)$ we use notations of [B], Section 1. Let us recall them. Namely, let $e_1, \dots, e_n$ be a fixed (standard) basis of the space $\p^n$;
\medskip
$I=\{i_1, \dots, i_d\}$, where $i_1< ... <i_d$, a subset of $\{1,\dots, n\}$;
\medskip
$E_I$ the $d$-dimensional subspace of $\p^n$ spanned on $e_{i_1},\dots, e_{i_d}$. 
\medskip
Further, let $U\subset GL_n(\p)$ be the group of upper unitriangular matrices of size $n$ and $U^I$ the subgroup of $U$ consisting of elements $(a_{ij})$ such that 
$$\hbox{ For $i<j$ we have: $a_{ij}\ne 0 \ \implies i\not\in I$ and $ j\in I $}.\eqno{(19.1.1)}$$
\medskip
{\bf Example.} Let $n=7$, $d=4$, $I=\{2,4,6,7\}$. $U^I$ is the set of matrices where possible non-0 entries are marked by asterisks:
$$\left(\matrix 1&*&0&*&0&*&*\\ 
0&1&0&0&0&0&0\\
0&0&1&*&0&*&*\\
0&0&0&1&0&0&0\\
0&0&0&0&1&*&*\\
0&0&0&0&0&1&0\\
0&0&0&0&0&0&1\endmatrix \right)$$
The set of Schubert cells of $Gr(d,n)$ is in 1 -- 1 correspondence with the set of subsets $I\subset \{1,\dots, n\}$ of order $d$. The Schubert cell corresponding to $I$ is denoted by $C_I$. The set of elements of $C_I$ is in 1 -- 1 correspondence with $U^I$. Namely, let $g\in U^I$. The element of $C_I$ corresponding to $g$ is $gE_I$ (recall that an element of $Gr(d,n)$ is a $d$-dimensional subspace in $\p^n$). Clearly the sets $C_I$ depend on a choice of a basis $\{e_1,\dots, e_n\}$ (recall that we consider it as a fixed basis). 
\medskip
If $I=\{n-d+1,\dots,n\}$ is the maximal set, then $U^I$ is the set of the $\{(n-d),\ d\}$-block unitriangular matrices $g=\left(\matrix I_{n-d}&S\\ 
0&I_d\endmatrix \right)$, $C_I$ is the maximal Schubert cell, and the $\{(n-d)\times d\}$-matrix $S$ is called the Siegel matrix of $g$, or, equivalently, of $gE_I\in C_I$. 
\medskip
{\bf 19.1.2.} By analogy, for the case of any $I$ and any $g\in U^I$, we shall call the set of elements $a_{ij}$ of $g$ for $i,\ j$ satisfying (19.1.1), by the Siegel set of $g$ (or of $gE_I\in Gr(d,n)$). 
\medskip
{\bf 19.1.3.} We denote the set of the above $I$ ( = the set of Schubert cells) by $\g I$, and we denote the map from $Gr(d,n)$ to $\g I$ associating to any element its Schubert cell, by $\s$. 
\medskip
{\bf 19.2. Map from the set of bases to the Grassmannian.} 
\medskip
Let $l_*=\{l_1, \dots, l_r\}$ be elements of $\p^n$ such that they are a $\n F_q[\th]$-basis of a lattice $L$ (see Definition 4.1); we call $l_*$ simply a L-basis. Two L-bases $l_{1*}=\{l_{11}, \dots, l_{1r}\}$ and $l_{2*}=\{l_{21}, \dots, l_{2r}\}$ are called equivalent if $\exists \ \vf\in GL(\p^n)$ such that $\vf(l_{1i})=l_{2i}$. The set of L-bases up to equivalence is denoted by $\Cal B=\Cal B(n,r)$. 
\medskip
Recall (see [Sh], or Section 1.8) that for the case of abelian varieites with MIQF an $O_K$-basis of a lattice $L\subset \n C^r$ defines a Siegel matrix $S\in Mat_{n,r-n}(\n C)$ satisfying (1.8.1), i.e. an element of the affine space $\n C^{n(r-n)}$. In the function field case, instead of $Mat_{n,r-n}(\n C)$, the target of the corresponding map (denoted by $\de$) is the Grassmannian $Gr(r-n,r)$; it contains $Mat_{n,r-n}(\n C)$ as the maximal Schubert cell. 
\medskip
{\bf 19.2.1. Definition of the map $\de: \cb \to Gr(r-n,r)$.} Let $V=\p^r$ and $e_1,\dots, e_r$ the standard basis of $V$. Let $l_*$ be a L-basis. A map $\sigma=\si(l_*): V\to \p^n$ is defined as follows: $\sigma (e_i)=l_i$. 
\medskip
{\bf Definition 19.2.1a.} $\de(l_*)=\Ker \sigma$. 
\medskip
Obviously $\de$ is well-defined: if $\vf\in GL(\p^n)$, then $\si(\vf(l_*))=\vf\circ\si(l_*)$, hence $\Ker \sigma(\vf(l_*))=\Ker \sigma(l_*)$, i.e. $\de(\vf(l_*))=\de(l_*)$.
\medskip
{\bf 19.2.2.} Obviously $\de$ is an inclusion. Really, let $l_*, \ l'_*$ be two L-bases such that $\de(l_*)=\de(l'_*)$. Let $x\in \p^n$. $\vf(x)$ is defined as $\sigma'(y)$ where $y\in V$ is any element of $\sigma^{-1}(x)$. Since $\Ker \sigma=\Ker \sigma'$, $\vf(x)$ does not depend on a choice of $y$. Clearly $\vf(l_i)=l'_i$. 
\medskip
If $l_1, \dots, l_n$ is a $\p$-basis of $V$, then $\de(l_*)$ belongs to the maximal Schubert cell of $Gr(n,r)$, and $\de(l_*)$ is defined via its Siegel matrix.
\medskip
{\bf 19.2.3.} Let us describe the Schubert cell of $\de(l_*)$, i.e. $\s\circ\de(l_*)$. Let $l_*$ be a L-basis. There exists a unique sequence of numbers $\al_*:=(\al_1, \dots, \al_n)$ such that $1=\al_1<\al_2< ... < \al_n\le r$ defined as follows: 
\medskip
$\forall \  k\in [1,\dots, n]$ we have: 
\medskip
(19.2.4) The dimension of the $\p$-linear envelope of $l_1, l_2, \dots, l_{\al_k}$ is $k$;
\medskip
(19.2.5) $\forall \  \be\in [\al_k,\dots, \al_{k+1}-1]$ we have: the dimension of the $\p$-linear envelope of $l_1, \dots, l_{\be}$ is $k$ (i.e. $\al_{k+1}$ is the minimal number such that the dimension of the $\p$-linear envelope of $l_1, \dots, l_{\al_{k+1}}$ becomes $k+1$). 
\medskip
{\bf Remark 19.2.6.} Since $l_1\ne0$, we have alwais $\al_1=1$. 
\medskip
Particularly, for $i=1,\ 2, \dots, \al_2-2$ we have $$l_{1+i}=s_{1i1}l_1;$$
For $i= 1,  2, \dots, \al_3-\al_2-1$ we have $$l_{\al_2+i}=s_{2i1}l_1+s_{2i2}l_{\al_2};$$
etc.: for $i=1, 2, \dots, \al_{j+1}-\al_j-1$ we have $$l_{\al_j+i}=s_{ji1}l_1+s_{ji2}l_{\al_2}+...+s_{jij}l_{\al_j}\eqno{(19.2.7)}$$
(if $j=n$ then we let $\al_{j+1}-1=r$). 
\medskip
Obviously $\al_*$ and $s_{kij}$ depend only on an element of $\cb$ (i.e. not on its representative $l_*$). 
\medskip
Let $b_*\in \cb$. We shall call the set of $s_{kij}$ from (19.2.7) as the Siegel set of $b_*$. This definition is concordant with the Siegel set of $\de(b_*)$ from (19.1.2), see the below 19.2.8. %We denote the set of Siegel sets of all elements of $\cb$ by $\g S$, and the set of $I$ of cardinality $n$ by $\g I$. Since an element of $\g S$ defines uniquely its $I$, we get a map $\g S\to \g I$. 
\medskip
For the case of the maximal Schubert cell, i.e. for $I=\{1,\dots, n\}$, the Siegel set is a Siegel matrix of size $(r-n)\times n$. 
\medskip
{\bf Proposition 19.2.8.} The set $I$ of the Schubert cell of $\de(l_*)$ (with respect to the basis $e_1, \dots, e_r$) is the complement of $\al_*$, i.e. $\s\circ\de(l_*)=I=\{1,\dots, r\}-\al_*$. 
\medskip
{\bf Proof} is immediate. A basis of $\Ker \de$ consists of $r-n$ elements from (19.2.7):
$$\vk_{\al_j+i}:=l_{\al_j+i}-(s_{ji1}l_1+s_{ji2}l_{\al_2}+...+s_{jij}l_{\al_j}),\eqno{(19.2.9)}$$ where $j=1,\dots, n$ and $i=1, 2, \dots, \al_{j+1}-\al_j-1$. 

Let us define $g$ be from (19.1.2). It is sufficent to define its non-trivial entries $a_{lm}$. For fixed $l, \ m$ we have: 
$$\exists\ i,\ j, \ k \hbox{ such that } l=\al_i, \ m=\al_j+k$$ where $k=1,\dots, \al_{j+1}-\al_j-1$. 
We let: 
$$a_{lm}:=-s_{jki}$$
This formula defines $g$. 

Example for $n=3$: the lines containing asteriscs have numbers $\al_1=1, \ \al_2, \ \al_3$; the columns whose off-diagonal entries are 0, also have numbers $\al_1=1, \ \al_2, \ \al_3$; all elements below the diagonal are 0: 
$$g=\left(\matrix 1&*&*&\dots&*&*&0&*&*&\dots&*&*&0&*&*&\dots&*&*\\ 
&1&0&\dots&0&0&0&0&0&\dots&0&0&0&0&0&\dots&0&0\\ \\
&&&&\dots&&&&&\dots&&&&&&\dots&&\\ \\
&&&&&1&0&0&0&\dots&0&0&0&0&0&\dots&0&0\\
&&&&&&1&*&*&\dots&*&*&0&*&*&\dots&*&*\\ 
&&&&&&&1&0&\dots&0&0&0&0&0&\dots&0&0\\ \\
&&&&&&&&&&\dots&&&&&\dots&&\\ \\
&&&&&&&&&&&1&0&0&0&\dots&0&0\\
&&&&&&&&&&&&1&*&*&\dots&*&*\\
&&&&&&&&&&&&&1&0&\dots&0&0\\ \\
&&&&&&&&&&&&&&&\dots&&\\ \\
&&&&&&&&&&&&&&&&1&0\\
&&&&&&&&&&&&&&&&&1\endmatrix \right)$$
It is clear that for $k=1,\dots, \al_{j+1}-\al_j-1$ we have $\vk_{\al_j+i}=g(e_{\al_j+i})$. $\square$
\medskip
Example of $g$ for $n=3, \ r=10, \ \al_1=1, \ \al_2=4, \al_3=7$:
$$g=\left(\matrix 1&-s_{111}&-s_{121}&0&-s_{211}&-s_{221}&0&-s_{311}&-s_{321}&-s_{331}\\ 
0&1&0&0&0&0&0&0&0&0\\ 
0&0&1&0&0&0&0&0&0&0\\ 
0&0&0&1&-s_{212}&-s_{222}&0&-s_{312}&-s_{322}&-s_{332}\\ 
0&0&0&0&1&0&0&0&0&0\\ 
0&0&0&0&0&1&0&0&0&0\\ 
0&0&0&0&0&0&1&-s_{313}&-s_{323}&-s_{333}\\ 
0&0&0&0&0&0&0&1&0&0\\ 
0&0&0&0&0&0&0&0&1&0\\ 
0&0&0&0&0&0&0&0&0&1\\ 
\endmatrix \right)$$
\medskip
{\bf 19.3. Action of $GL_r(\n F_q[\th])$ on $\cb$ and on $Gr(r-n,r)$.}
\medskip
There are the natural actions of $GL_r(\n F_q[\th])$ on $\cb$ and of $GL_r(\p)$ on $Gr(r-n,r)$. They are compatible up to the power of $-1$. 
\medskip
{\bf Proof} is immediate: let $\ga=\{\ga_{ij}\}\in GL_r(\n F_q[\th])$; we denote the $(i,j)$-th entry of $\ga^{-1}$ by $\bar \ga_{ij}$. Let $l'_*:=\ga(l_*)$ and $\sigma'$ be its map $\sigma$. 
\medskip
Let $x=\sum \la_ie_i\in \Ker \sigma$. This means that $\sum \la_il_i=0$ in $\p^n$. We must prove that $\ga^{-1}(x)\in \Ker \sigma'$. 
\medskip
We have: $\ga^{-1}(x)=\sum \la_i \sum \bar \ga_{ij}e_j$. Since $\sigma'(e_j)=\sum \ga_{jk}l_k$, we have $$\sigma'(\ga^{-1}(x))=\sum \la_i \sum \bar \ga_{ij}\sum \ga_{jk}l_k=\sum \la_il_i=0.\eqno{\square}$$
\medskip
{\bf 19.4. t-modular forms.}
\medskip
The classical notion of Siegel modular forms can be word-to-word generalized to the above case. Let us give definitions. We shall call them t-modular forms (because of t-motives). 
\medskip
Let $\rho: GL(\p^n)\to GL(\g W)$ be a representation of $GL(\p^n)$ on a finite-dimensional $\p$-vector space $\g W$. Let $\cl$ be the set of all lattices of rank $r$ in $\p^n$. Let $\vf\in GL(\p^n)$. 
\medskip
A function $f:\cl \to \g W$ is called a t-modular form of weight $\rho$, if $\forall \ L\in \cl,\ \vf\in GL(\p^n)$ we have $$f(\vf(L))=\rho(\vf)^{-1}(f(L))$$ + conditions of meromorphy and on the boundary.
\medskip
This $f$ gives us a function (we denote it by $\bar f$) from $\cb$ to $\g W$. Namely, let $l_*$ be a L-basis --- a representative of an element of $\cb$, and $\al_1=1, \al_2, \dots, \al_n$ be from (19.2.3). Recall that the standard basis of $\p^n$ is denoted by $e_1,\dots, e_n$. There exists the only $\vf\in GL(\p^n)$ such that $\forall \ i$ we have $\vf(l_{\al_i})=e_i$. Let $L$ be the lattice generated by $\vf(l_*)$. By definition, 
$$\bar f(l_*):=f(L).$$
Clearly $\bar f$ is well-defined. 
\medskip
Let us recall the classical formula for Siegel modular forms. Let $\Cal H_g$ be the Siegel upper half plane, $f: \Cal H_g \to \g W$ a Siegel modular form of weight $\rho$, and $\ga=\left(\matrix A&B\\ C&D\endmatrix \right)\in GSp_{2g}(\n Z)$, $S\in \Cal H_g$. For all $\ga, \ S$ we have $CS+D\in GL_g(\n C)$. The automorphy condition is 
$$f(\ga(S))=\rho(CS+D)(f(S)),\eqno{(19.4.1)}$$
the matrix $CS+D$ is called the automorphy factor for $\ga$ and $S$. 
\medskip
{\bf 19.4.1a.} Analogously, we can define modular forms on $\Cal H^3_{n,r-n}$ (see (1.8.1) and the below lines). In this case $\ga\in GU(n,r-n)(\n Z)=G^3_{n,r-n}(\n Z)$ (see [Sh], 2.7, page 163), sizes of blocks of $\ga$ are $n, \ r-n$, and for $S\in \Cal H^3_{n,r-n}$ we have always $\ga(S)\in \Cal H^3_{n,r-n}, \ CS+D\in GL_{r-n}(\n C)$. The automorphy condition is the same (19.4.1). 
\medskip
The analog of (19.4.1) for the function field case is the following. We consider $\cb$ as a subset of $Gr(r-n,r)$ via $\de$. Let $x\in \cb$ belong to the maximal Schubert cell of $Gr(r-n,r)$ and has the Siegel matrix $S$ of size $(r-n,\ n)$, and let $\ga=\left(\matrix A&B\\ C&D\endmatrix \right)\in GL_r(\n F_q[\th])$ (sizes of blocks are $r-n, \ n$) be such that $$CS+D\in GL_n(\p).\eqno{(19.4.2)}$$
Condition (19.4.2) is equivalent to the condition that $\ga(x)$ belongs to the maximal Schubert cell of $Gr(r-n,r)$. In this case the formula (19.4.1) is the same:\footnotemark \footnotetext{Maybe for the present situation, in order to form the automorphy factor, we should use not the blocks $C$ and $D$ of $\ga$, but other blocks. This depends on notations and is not important.}$$\bar f(\ga(S))=\rho(CS+D)(\bar f(S)).\eqno{(19.4.3)}$$
\medskip
We can get the analog of (19.4.3) for the case of non-maximal Schubert cells of $Gr(r-n,r)$. Let $x\in \cb$, and the Schubert cell of $x$ corresponds to a subset $\s\circ\de(x)=I\subset \{1,\dots,r\}$. We describe $x$ by its Siegel set (more exactly, by the Siegel set of $\de(x)$) from (19.1.2); we denote it by $S(x)$. Let $\ga\in GL_r(\n F_q[\th])$. Objects $I, \ S(x), \ \ga$ define: 
\medskip
First, a Schubert cell $\s\circ\de(\ga(x))$ (corresponding to a subset of $\{1,\dots,r\}$ denoted by $I'$) of $\ga(x)$;
\medskip
Second, the Siegel set of $\ga(x)$;
\medskip
Third, the automorphy factor for $\ga$ and $S(x)$ (an analog of $CS+D$ for the case when both $I$ and $I'$ are $\{1,\dots,n\}$). 
\medskip
All entries of these objects are rational functions of $S(x)$ and $\ga$, and can be explicitly written for all possible $I$, $I'$. For all $x$ and $\ga$ we have: $\bar f(\ga(x))$ must satisfy an analog of (19.4.3). 
\medskip
{\bf Remark 19.4.4.} I am not sure that t-modular forms really exist (except the zero form). Does there exist, or not, a good quotient space $\cb/GL_r(\n F_q[\th])$? Recall two results from [GL17]. We need the following definitions. Let $\om\in \n F_{q^2}-\n F_q$ (it is an analog of $i=\sqrt{-1}\in \n C$). Let $r=2n$. The maximal Schubert cell for this case is $M_{n\times n}(\p)$, and $GL_{2n}(\n F_q[\th])$ "almost" acts on it ("almost" means that $\ga(S)$ not always belongs to maximal Schubert cell). Let $G_0\subset GL_{2n}(\n F_q[\th])$ be the stabilizer of $\om I_n$. 
\medskip
We have ([GL17], Proposition 1.7.1):
\medskip
For any neighborhood $U$ of $\omega I_n$ in $M_{n\times n}(\p)$ there exist $S_1$, $S_2\in U$, $\gamma\in GL_{2n}(\Bbb F_q[\theta])$ such that $S_2=\gamma(S_1)$ and $\gamma\not\in G_0$.
\medskip
So, the situation for the function field case is worse than for the number field case. Nevertheless, we have ([GL17], Proposition 1.7.2):
\medskip
There exists a neighborhood $U$ of $\omega I_n$ in $M_{n\times n}(\p)$ such that if $\ga\in GL_{2n}(\Bbb F_q[\theta])$ satisfies $\ga(\omega I_n)\in U$ then $\ga\in G_0$.
\medskip
Finally, to develop the theory of t-modular forms, it is important to know that there is a 1 -- 1 correspondence between pure uniformizable t-motives and lattices. We know only that the lattice map for them is injective ([HJ]) and locally surjective ([GL17]). 
\medskip
{\bf 19.5. Application to the duality.}
\medskip
{\bf Proposition 19.5.1.} Let $l_*$ be a L-basis having the subset $\s\circ\de(x)=I\subset \{1,\dots,r\}$. Then the dual basis $\hat l_*$ has the subset $\hat I:=r+1 - (\{1,\dots,r\}-I)$. The Siegel set of $\hat l_*$ is the minus transposed of the Siegel set of $l_*$. 
\medskip
{\bf Proof.} To write (it is immediate). 
\medskip
This is a generalization of Theorem 12.10. 
\medskip
{\bf 19.6. Case of any $N$.}
\medskip
{\bf 19.6.1.} First, we recall the definition and properties of the "classical" flag varieties, see for example [B], Section 1, or [GL18], Section 3b. Let $r$ be, as above, the rank of a t-motive $M$, and $$k_1+k_2+...+k_{m+1}=r$$ an ordered partition of $r$. Here $k_i\ge0$, they come from some invariants of the nilpotent operator $N$ attached to $M$, see below, and $m$ is defined by the condition $N^m=0$. Particularly, if $N=0$, i.e. $m=1$, then $k_1=n, \ k_2=r-n$.
\medskip
The "classical" flag variety $X(k_*)=X(k_1, \dots,k_{m+1})$ is the set of flags of type $k_1, \dots,k_{m+1}$ in $\p^r$, where a flag of type $k_1, \dots,k_{m+1}$ in $\p^r$ is the set of vector subspaces
$$0=V_0\subset V_1 \subset V_2 \subset ...  \subset V_m \subset V_{m+1}=\p^r\eqno{(19.6.1.1)}$$where $\dim V_i/V_{i-1}=k_i$. It is well-known that $$\dim X(k_1, \dots,k_{m+1})=\sum_{1\le u<y\le m+1}k_uk_y\eqno{(19.6.1.2)}$$ 

An equivalent definition: let $G=GL_r$, and $P_{k_*}$ be the parabolic subgroup consisting of upper block diagonal matrices, where sizes of diagonal blocks are $k_1, \dots,k_{m+1}$. We have: $$X(k_1, \dots,k_{m+1})=G(\p)/P_{k_*}(\p)\eqno{(19.6.1.3)}$$
Also, let $W$ be an abstract variable and $I_{k_*}\subset G(\p[[W]])$ be the Iwahori subgroup of type $k_*$ defined as follows: let $\pi: G(\p[[W]]) \to G(\p)$ be the natural projection ($W\mapsto0$), then $I_{k_*}:=\pi^{-1}(P_{k_*}(\p))$. We have:
$$X(k_1, \dots,k_{m+1})=G(\p[[W]])/I_{k_*}\eqno{(19.6.1.4)}$$
Let us mention that it is possible to define the affine flag variety of type $k_*$:
$$Fl_{aff}(k_*):=G(\ \p((W)) \ )/I_{k_*}\eqno{(19.6.1.5)}$$
although apparently it has no applications in the theory of t-motives. 
\medskip
Typical objects of study are the complete (affine) flag varieties that correspond to $k_*=\{1,\dots, 1\}$. 
\medskip
The set of Schubert cells of $X(k_*)$ is the set of increasing sequences of subsets of $\{1,\dots,r\}$:
$$\emptyset=I_0\subset I_1\subset I_2 \subset ... \subset I_m\subset I_{m+1}=\{1,\dots,r\}$$
such that $\#(I_\al)=k_1+...+k_\al$. 
\medskip
{\bf 19.6.1.6.} Equivalently, we can use notations $J_\al:=I_\al-I_{\al-1}$. We have 
\medskip
$\#(J_\al)=k_\al$; \ $J_\al \cap J_\be=\emptyset$, \ $\underset{\al}\to\bigcup \ J_\al=\{1,\dots,r\}$.  What notation is better?
\medskip
We denote a sequence of this type by $\ci=\{I_1\subset I_2 \subset ... \subset I_m\subset I_{m+1}\}$ and the set of these $\ci$ by $\g I$, like in 19.1.3. Also, we have $\s: X(k_*)\to \g I$. 
\medskip
For $\ci\in \g I$ we denote by $E_\ci$ the flag of type (19.6.1.1) such that $\forall \ \al=1,\dots, m$ we have: $V_\al$ is spanned by the elements $e_\be$, where $\be\in I_\al$. The Schubert cell $C_\ci$ corresponding to $\ci$ is $UE_\ci$. 
\medskip
The minimal Schubert cell (it is a point) corresponds to $\ci$ denoted by $\ci_{min}$ having $I_\al=\{1,\dots, k_1+...+k_\al\}$, and the maximal Schubert cell (it is open) corresponds to $\ci$ denoted by $\ci_{max}$ having $I_\al=\{r+1-(k_1+...+k_\al),\dots, r\}$. 
\medskip
Let $U_\ci\subset U$ be the stabilizer of $E_\ci$. It is described as follows. For $j\in \{1,\dots,r\}$ we denote by $\vk(j)$ the minimal number such that $j\in I_{\vk(j)}$. Then 
$$g=\{a_{ij}\}\in U_\ci \iff \forall \ i, \ j \hbox{ such that $i<j$ we have: } \vk(i)>\vk(j)\implies a_{ij}=0$$

We have $U_{\ci_{min}}=U$, $$U_{\ci_{max}}=U\ \cap \hbox { \{the set of block diagonal matrices\} }$$(sizes of blocks are $k_1, \dots,k_{m+1}$). 
\medskip
Inversely, let $U^\ci\subset U$ be a subgroup of $U$ defined as follows: 
$$g=\{a_{ij}\}\in U^\ci \iff \forall \ i, \ j \hbox{ such that $i<j$ we have: } \vk(i)\le \vk(j)\implies a_{ij}=0.$$

Like for the case of the Grassmannians, the map $U^\ci\to C_\ci$ defined by $g\mapsto g\cdot E_\ci$, is a 1 -- 1 map. 
\medskip
We have $U^{\ci_{min}}=1$, $U^{\ci_{max}}$ is the group of block unitriangular matrices, i.e. if $g\in U^{\ci_{max}}$ then the $(\al,\be)$-th block of $g$ is 0, if $\al>\be$, and the identity matrix $I_{k_\al}$ if $\al=\be$. 
\medskip
{\bf 19.6.1.7.} Particularly, if $x\in X(k_*)$ belongs to the maximal Schubert cell $C_{\ci_{max}}$ then the analog of its Siegel matrix is the set of matrices $S_{\al,\be}$, $1\le \al <\be\le m+1$, which are $(\al,\be)$-th blocks of $g\in U^{\ci_{max}}$ corresponding to $x$; $S_{\al,\be}$ is a matrix of size $k_\al\times k_\be$. 
\medskip
{\bf 19.6.2.} The analog of $Gr(r-n,r)$ for the case $N\ne0$ is not the "classical" flag variety $X(k_1, \dots,k_{m+1})$ but its version that is called a t-flag variety (tFV; t from t-motives), see Definition 19.6.4.6 (to prove this fact!). Let us recall that the definition of a $N$-lattice is given in 10.5.1, and the definition of an isomorphism of $N$-lattices is given in 10.5.5 (it is necessary to emphasize that $\ga: V\to V$ such that $\ga(L_1)=L_2$ {\bf must commute with $N$}). 
\medskip
Let us give analogs of the initial definitions of 19.2. Let $l_*=\{l_1,\dots, l_r\}$ be a set of elements of $\p^n$ which are a basis of a $N$-lattice, i.e. elements of (10.5.2) form a $N$-lattice. Elements $l_*$ are called a $NL$-basis. Two $NL$-bases $l_{1*}=\{l_{11}, \dots, l_{1r}\}$ and $l_{2*}=\{l_{21}, \dots, l_{2r}\}$ are called equivalent if $\exists \ \vf\in GL(\p^n)$ {\bf commuting with $N$} such that $\vf(l_{1i})=l_{2i}$. The set of $NL$-bases up to equivalence is denoted by $\Cal B=\Cal B(N)$. 
\medskip
{\bf 19.6.3. Definition of $k_i$.}  Recall that $N: \p^n\to\p^n$ is a nilpotent operator. Let $d_1,\dots, d_\al$ be sizes of its 0-Jordan blocks. Hence, we have a partition $$n=d_1+...+d_\al\eqno{(19.6.3.1)}$$ denoted by $d_*$, where $m=d_1\ge d_2\ge...\ge d_\al>0$. 
Let $$n=c_1+...+c_m$$ be the conjugate to the partition (19.6.3.1). We have 
\medskip
$\al=c_1\ge c_2\ge...\ge c_{m}>0$. 
\medskip
Now, we let $c_0:=r$, $c_{m+1}:=0$ and $$\hbox{for }i=1,\dots, m+1 \ \ \ \ k_i:=c_{i-1}-c_i.\eqno{(19.6.3.2)}$$
\medskip
We have $$k_i\ge0, \ \ \ r=\sum_{i=1}^{m+1}k_{i}, \ \ \ \ \ n= k_2+2k_3+3k_4+... + m k_{m+1}.\eqno{(19.6.3.3)}$$ 

Let us mention other useful formulas: $$\dim \Ker N^{i}=c_1+...+c_i,$$ $$\dim N^{m-1}V=k_{m+1},\eqno{(19.6.3.4)}$$
$$\dim N^{m-2}V=k_m+2k_{m+1},\eqno{(19.6.3.5)}$$ $$\dim N^{m-3}V=k_{m-1}+2k_m+3k_{m+1}\hbox{ etc., }\eqno{(19.6.3.6)}$$
$$\hbox{For }i\ge2:\ k_{i}=\dim (\Ker N^{i-1}/ \Ker N^{i-2})- \dim (\Ker N^{i}/ \Ker N^{i-1}) \eqno{(19.6.3.7)}$$
$$=\dim (\Im N^{i-2} / \Im N^{i-1}) - \dim (\Im N^{i-1} / \Im N^{i})\eqno{(19.6.3.8)}$$

{\bf 19.6.3.9.} It is easy to check that the dimension of the set of $\ga: V\to V$ commuting with $N$ is $$rn-\sum_{1\le i < j\le m+1}(j-i)k_ik_j\eqno{(19.6.3.10)}$$ (it does not depend on $r$, but only on $k_2,\dots,k_{m+1}$). See the below formula 19.6.4.7 for its application. 
\medskip
Does exist a better form of this formula? 
\medskip
{\bf Example 19.6.3.11.} For $m=1$ we have: the pair $(k_1,k_2)$ is $(r-n,n)$.
\medskip
{\bf Remark 19.6.3.12.} U. Hartl and A.-K. Juschka consider in [HJ] another invariants describing $N$, namely, invariants $\la_i$ from Remark 19.6.4.10 instead of $k_*$. 
\medskip
{\bf 19.6.4. t-flag variety.}
\medskip
{\bf Definition 19.6.4.1.} The t-Iwahori subgroup of type $k_*$ (notation: $tI_{k_*}$) is a subgroup of $G(\p[[W]])$ consisting of block matrices such that their diagonal blocks are square blocks of sizes $k_1,\dots,k_{m+1}$, and for any $$1\le u<y\le m+1\eqno{(19.6.4.2)}$$ all entries of the $(y, \ u)$-th block of this matrix belong to $W^{y-u}\cdot \p[[W]]$. 
\medskip
Example for $m=2$:
$$\left(\matrix &&&|&&&&|\\ &*&&|&&*&&|&*&\\ &&&|&&&&|\\ -&-&-&-&-&-&-&-&-&- \\ &&&|&&&&|\\ W\cdot & \p&[[W]]&|&&*&&|&*&\\ &&&|&&&&|\\ -&-&-&-&-&-&-&-&-&- \\ &&&|&&&&|\\ W^2\cdot & \p&[[W]]&|&W\cdot & \p&[[W]]&|&*&\\ &&&|&&&&|\endmatrix \right)\eqno{(19.6.4.3)}$$

{\bf Remark 19.6.4.4.} Some $k_i$ can be 0. For example, if $m=4$ and $k_2=0$, then $tI_{k_*}$ is the set of matrices
$$\left(\matrix &&&|&&&&|\\ &*&&|&&*&&|&*&\\ &&&|&&&&|\\ -&-&-&-&-&-&-&-&-&- \\ &&&|&&&&|\\ W^2\cdot & \p&[[W]]&|&&*&&|&*&\\ &&&|&&&&|\\ -&-&-&-&-&-&-&-&-&- \\ &&&|&&&&|\\ W^3\cdot & \p&[[W]]&|&W\cdot & \p&[[W]]&|&*&\\ &&&|&&&&|\endmatrix \right)\eqno{(19.6.4.5)}$$
(sizes of diagonal blocks are $k_1, \ k_3, \ k_4$). 
\medskip
{\bf Definition 19.6.4.6.} The variety tFV$(k_*)$ is $G(\p[[W]])/tI_{k_*}$. 
\medskip
Clearly $$\dim_{\p}(tFV(k_*))=\sum_{1\le u<y\le m+1}(y-u)k_uk_y\eqno{(19.6.4.7)}$$   i.e. it differs from the above formula (19.6.1.2) by coefficients $y-u$. 
\medskip

{\bf Theorem  19.6.4.8.} The representatives of the maximal Schubert cell of $G(\p[[W]])/tI_{k_*}$ are block lower unitriangular matrices such that $\forall \ y, \ u$ (as earlier $1\le u<y\le m+1$) all entries of their $(y, \ u)$-th block are polynomials in $W$ of degree $\le y-u$. 
\medskip
{\bf Proof} is straightforward. This is an analog of (19.6.1.7). 
\medskip
Example for $m=2$: these representatives are
\medskip
$$\left(\matrix &&&|&&&&|\\ &I_{k_1}&&|&&0&&|&0&\\ &&&|&&&&|\\ -&-&-&-&-&-&-&-&-&- \\ &&&|&&&&|\\ &A_{210}&&|&&I_{k_2}&&|&0&\\ &&&|&&&&|\\ -&-&-&-&-&-&-&-&-&- \\ &&&|&&&&|\\ A_{310}&+&A_{311}W&|&&A_{320}&&|&I_{k_3}&\\ &&&|&&&&|\endmatrix \right)\eqno{(19.6.4.9)}$$
\medskip
where $A_{yu\al}$ is a $(k_y\times k_u)$-matrix with entries in $\p$. (to check: the transposed matrix!)
\medskip
Obviously the quantity of these elements is the same as in the right hand side of (19.6.4.7). 
\medskip

{\bf Remark 19.6.4.10.} A similar, but not exactly the same flag variety, was considered in [R16], Section 2.1. Let $T\subset GL_r$ be the diagonal matrices and $\chi\in X_*(T)^+$. We can identify $\chi$ with a sequence of integers $\la_1\le \la_2 \le ... \le \la_r$. If we fix $k_*$ then a corresponding $\chi$ is defined as follows: $\la_1$ is arbitrary (i.e. $\la_1, \la_2, \dots,  \la_r$ are defined mod $\n Z$), 
\medskip
$\la_1= \la_2 = ... = \la_{k_1}$, 
\medskip
$\la_{k_1+1}= \la_{k_1+2} = ... = \la_{k_1+k_2}=\la_1+1$, 
\medskip
$\la_{k_1+k_2+1}= \la_{k_1+k_2+2} = ... = \la_{k_1+k_2+k_3}=\la_1+2$, 
\medskip
$\dots$
\medskip
$\la_{k_1+k_2+... +k_{m}+1}= \la_{k_1+k_2+... +k_{m}+2} = ... = \la_{r}=\la_1+m$. 
\medskip
This $\chi$, and hence $k_*$, defines a vertex in the Bruhat-Tits building of $G(\p[[W]])$ belonging to the standard apartment. Its stabilizer $G_\chi$ is a subgroup of $G(\p[[W]])$ consisting of block matrices such that their diagonal blocks are square blocks of sizes $k_1,\dots,k_{m+1}$, and for any $u,\ y \in [1, m+1]$ all entries of the $(y, \ u)$-th block of this matrix belong to $W^{y-u}\cdot \p[[W]]$. 
\medskip
{\bf 19.6.4.11.} We see that the difference from the Definition 19.6.4.1 is the absence of the condition (19.6.4.2): $u<y$. Example for $m=2$ (compare with (19.6.4.3)): 
$$\left(\matrix &&&|&&&&|\\ &\p&[[W]]&|&W^{-1}\cdot & \p&[[W]]&|&W^{-2}\cdot & \p&[[W]]&\\ &&&|&&&&|\\ -&-&-&-&-&-&-&-&-&-&-&- \\ &&&|&&&&|\\ W\cdot & \p&[[W]]&|&&\p&[[W]]&|&W^{-1}\cdot & \p&[[W]]&\\ &&&|&&&&|\\ -&-&-&-&-&-&-&-&-&- &-&-\\ &&&|&&&&|\\ W^2\cdot & \p&[[W]]&|&W\cdot & \p&[[W]]&|&&\p&[[W]]&\\ &&&|&&&&|\endmatrix \right)\eqno{(19.6.4.12)}$$

Really, Richarz considers in [R16], Section 2.1 a more general object than 
\medskip
$\chi$ = a vertex in the Bruhat-Tits building of $G(\p[[W]])$. 
\medskip
Namely, he considers a facet $\g a$. It is clear that its stabilizer $G_\g a$ is a subgroup of the above $G\chi$, but it is never $tI_{k_*}$ from Definition 19.6.4.1. 
\medskip
\medskip
{\bf 19.6.5. $tFV(N)$ is an analog of $Gr(r-n,r)$ for the case $N\ne0$.} 
\medskip
{\bf Conjecture 19.6.5.1.} $tFV(N)$ is an analog of $Gr(r-n,r)$ for the case $N\ne0$. Namely, 
let $\Cal B=\Cal B(N)$ be from 19.6.2, and $tFV=tFV(N)$ from 19.6.4.6. There exists a canonically defined map $\de: \cb\to tFV$ which is an analog of the map of 19.2.1, case $N=0$. 
\medskip
At the moment we have no definition of $\de$ for the case $N\ne0$ which is an analog of Definitions 19.2.1, 19.2.1a for the case $N=0$, because we have no interpretation of elements of $tFV(N)$ as subspaces of $\p^r$. We have only 
a coordinate description of $\de$ on the level of representatives of both $\cb$ and $tFV(N)$ for the elements of the maximal Schubert cell of $tFV(N)$. 
\medskip
Hence, in principle it can happen that an analog of $Gr(r-n,r)$ for the case $N\ne0$ is not exactly $tFV(N)$ but some other variety which is "very similar to it".
\medskip
Let us give the above mentioned coordinate description of $\de$ (we follow [GL18], Sections 3a.4, 3a.5). See also 19.6.7 for the case of non-maximal Schubert cell. 
\medskip
Let $l_*=(l_1, \dots, l_r)$ be an element of $\Cal B$. We have: $\de(l_*)$ belongs to the maximal Schubert cell of tFV iff the below conditions (19.6.5.3), (19.6.5.6), (19.6.5.8) etc. hold. 
\medskip
First step: Elements $N^{m-1}l_j$, $j=1,\dots,r$, generate $N^{m-1}V$ as a $\p$-vector space. (19.6.3.4) shows that the dimension of $N^{m-1}V$ is $k_{m+1}$. 
\medskip
{\bf 19.6.5.2.} Condition that $\de(l_*)$ belongs to the maximal Schubert cell of tFV implies that 
$$\hbox{the last $k_{m+1}$ elements from $l_1,\dots,l_r$ form a $\p$-basis of $N^{m-1}V$.}\eqno{(19.6.5.3)}$$ 
We denote these elements by $l_{m+1,1}, \dots, l_{m+1,k_{m+1}}$, and their set by $\hat l_{m+1}$. 
\medskip
Further on (second step), elements $N^{m-2}l_j$, $N^{m-1}l_j$, $j=1,\dots,r$, generate $N^{m-2}V$ as a $\p$-vector space. 
\medskip
{\bf 19.6.5.4.} Elements $N^{m-2}(l_{m+1,i})$, $N^{m-1}(l_{m+1,i})$, $i=1, \dots, k_{m+1}$, are linearly independent over $\p$, see [GL18], 3a.4.1 for a proof. 
\medskip
We have (see (19.6.3.5)) $\dim N^{m-2}V=2k_{m+1}+k_{m}$. Hence, we get:
\medskip
{\bf 19.6.5.5.} Condition that $\de(l_*)$ belongs to the maximal Schubert cell of tFV implies that:
\medskip
The last $k_{m}$ elements from the set $l_1,\dots,l_{r-k_{m+1}}$
(we denote them by $l_{m,1}, \dots, l_{m,k_{m}}$ respectively, and their set by $\hat l_{m}$) have the property:
\medskip
{\bf 19.6.5.6.} $N^{m-2}(l_{m,i})$, $i=1, \dots, k_{m}$, 
\medskip
$N^{m-2}(l_{m+1,i})$, $i=1, \dots, k_{m+1}$, 
\medskip
$N^{m-1}(l_{m+1,i})$, $i=1, \dots, k_{m+1}$, form a $\p$-basis of $N^{m-2}V$.
\medskip
{\bf 19.6.5.7.} The third step of the process: elements 
\medskip
$N^{m-3}l_j$, $N^{m-2}l_j$, $N^{m-1}l_j$, $j=1,\dots,r$, 
\medskip
generate $N^{m-3}V$ as a $\p$-vector space. 
\medskip
Elements $N^{m-3}(l_{m,i})$, $N^{m-2}(l_{m,i})$, $i=1, \dots, k_{m}$, 
\medskip
$N^{m-3}(l_{m+1,i})$, $N^{m-2}(l_{m+1,i})$, $N^{m-1}(l_{m+1,i})$, $i=1, \dots, k_{m+1}$, 
\medskip
are linearly independent over $\p$ (the proof is exactly the same as the proof of (19.6.5.4) ( = the proof of [GL18], 3a.4.1). 
\medskip
The quantity of these elements is $2k_m+3k_{m+1}$. We have (see (19.6.3.6)) $\dim N^{m-3}V=k_{m-1}+2k_m+3k_{m+1}$. 
\medskip
Hence, we get:
\medskip
{\bf 19.6.5.8.} Condition that $\de(l_*)$ belongs to the maximal Schubert cell of tFV implies that: 
\medskip
The last $k_{m-1}$ elements from the set $l_1,\dots,l_{r-k_{m+1}-k_m}$
(we denote them by $l_{m-1,1}, \dots,\ l_{m-1,k_{m-1}}$ respectively, and their set by $\hat l_{m-1}$) have the property:
\medskip
{\bf 19.6.5.9.} $N^{m-3}(l_{m-1,i})$, $i=1, \dots, k_{m-1}$, 
\medskip
$N^{m-3}(l_{m,i})$, $N^{m-2}(l_{m,i})$, $i=1, \dots, k_{m}$, 
\medskip
and $N^{m-3}(l_{m+1,i})$, $N^{m-2}(l_{m+1,i})$, $N^{m-1}(l_{m+1,i})$, $i=1, \dots, k_{m+1}$, 
\medskip
form a $\p$-basis of $N^{m-3}V$.
\medskip
{\bf 19.6.5.10.} Continuing this process, we represent the ordered set $\{l_1,\dots,l_{r}\}$ as a disjoint ordered union of segments: $$\{l_1,\dots,l_{r}\}=\hat l_1\cup \hat l_2\cup \dots\cup \hat l_{m+1}$$ where the length of $\hat l_{i}$ is $k_i$, i.e. 
$$\hat l_{1}=(l_1, \dots, l_{k_1}), \ \  \hat l_{2}=(l_{k_1+1}, \dots, l_{k_1+k_2}), \ \  \hat l_{3}=(l_{k_1+k_2+1}, \dots, l_{k_1+k_2+k_3})\hbox{ etc.}$$ (we use also a notation 
$$(l_1, \dots, l_{k_1})=(l_{1,1}, \dots, l_{1,k_1}); \ (l_{k_1+1}, \dots, l_{k_1+k_2})=(l_{2,1}, \dots, l_{2,k_2}), \hbox{ etc.)}$$ such that $\forall \ u=0,\dots,m-1$ we have:
\medskip
{\bf 19.6.5.11.} A $\p$-basis of $N^uV$ is formed by elements from $N^\al(\hat l_{\be})$, where $\al\in[u,\dots, m-1]$, $\be\in[\al+2,\dots, m+1]$.
\medskip
Explicitly, (19.6.5.11) can be described in a form of the diagram 
$$\matrix V&|& \hat l_{2}& \hat l_{3}& \hat l_{4}& \dots&\hat l_{m-1}& \hat l_{m}& \hat l_{m+1}\\ \\
         N(V)  &|&    & N(\hat l_{3})& N(\hat l_{4})& \dots&N(\hat l_{m-1})& N(\hat l_{m})& N(\hat l_{m+1})\\  \\ N^2(V)  &|&
           && N^2(\hat l_{4})& \dots&N^2(\hat l_{m-1})& N^2(\hat l_{m})& N^2(\hat l_{m+1})\\  \\
                \dots                  &|&    \dots       &\dots&\dots&\dots&\dots&\dots&\dots \ \ \  (19.6.5.12)\\  \\
N^{m-3}(V)&|&   &&&&N^{m-3}(\hat l_{m-1})& N^{m-3}(\hat l_{m})& N^{m-3}(\hat l_{m+1})\\  \\
N^{m-2}(V)&|&  &&&&& N^{m-2}(\hat l_{m})& N^{m-2}(\hat l_{m+1})\\  \\ N^{m-1}(V)&|& &&&&& & N^{m-1}(\hat l_{m+1})\endmatrix$$

Its meaning is the following. Let us consider the line corresponding to $N^u(V)$. Then elements $N^\al(\hat l_{\be})$ situated on this line and below it, form a basis of $N^u(V)$. 
We see that for any $u$ the set of these $(\al, \ \be)$ is a triangle. 
\medskip
Some $k_i$ can be 0. The corresponding sets $\hat l_i$ are empty. 
\medskip
{\bf 19.6.6. Description of $\de(l_*)$ via its Siegel set.}
\medskip
The idea of description of $\de(l_*)$ for the case when $\de(l_*)$ belongs to the maximal Schubert cell of $tFV$, is the following. We define a representative of $\de(l_*)$ in the form of Theorem 19.6.4.8. We can consider this reperesentative as a set of Siegel matrices parametrized by points of a tetrahedron. These matrices come from the fact that the sets $$\hat l_1, \ N(\hat l_{2}), \ N^2(\hat l_{3}), \dots, N^{u-1}(\hat l_u),\dots, N^{m-2}(\hat l_{m-1}), \ N^{m-1}(\hat l_{m})\eqno{(19.6.6.1)}$$ (they are under the bottom entries of the columns of (19.6.5.12) ) 
are linear combinations of all entries of (19.6.5.12) which are in the east --- southeast sector from them. 
\medskip
More exactly, parameters of the points of this tetrahedron are integers $u, \ z, \ y$ where $u$ is from (19.6.6.1), $y$ is the number of the column of (19.6.5.12), and $z$ is the number of the line of (19.6.5.12) (we neglect the fact that the numeration starts from 2 for columns, and from 0 for lines). The condition that $N^z\hat l_{y}$ is in the east --- southeast sector from $N^{u-1}(\hat l_u)$, is the following: 
$$u\in [1,m], \ \ z\in [u-1, m-1], \ \  y\in [z+2,m+1]\eqno{(19.6.6.2)}$$
It turns out that it is convenient to introduce one more parameter $v$ uniquely defined by $u$:
$$v=u-1\eqno{(19.6.6.3)}$$i.e. $N^{u-1}(\hat l_u)=N^{v}(\hat l_u)$. This notation will be used in duality theory, see [GL18], Secions 3, 4. 
\medskip
{\bf Definition 19.6.6.4.} The Siegel object of $\de(l_*)$ (case of the maximal Schubert cell of tFV) is the set of Siegel matrices $S_{uvyz}$, where $u, \ v, \ y, \ z$ satisfy (19.6.6.2), (19.6.6.3), of size $k_u\times k_y$ with entries in $\p$ such that the following holds:
$$N^{v}\hat l_{u}=-\sum_{z=u-1}^{m-1}\sum_{y=z+2}^{m+1}S_{uvyz}\cdot N^z\hat l_{y}\eqno{(19.6.6.5)}$$
Here $\hat l_{u}, \ \hat l_{v}$ are matrix columns;  if some $k_*$ are 0 then the corresponding $S_{****}$ are empty. The sign minus is not of principle, it comes from applications to the theory of Anderson t-motives. 
\medskip
{\bf 19.6.6.6.} Example for $m=3$, $u=1, \ v=0$:
$$\matrix \hat l_{1}=-(S_{1020}\cdot\hat l_{2} &+& S_{1030}\cdot\hat l_{3} &+& S_{1040}\cdot\hat l_{4} \\ \\
&+& S_{1031}\cdot N(\hat l_{3}) &+& S_{1041}\cdot N(\hat l_{4}) \\ \\ &&&+& S_{1042}\cdot N^2(\hat l_{4}) \ )\endmatrix \eqno{(19.6.6.7)}$$
(terms of a fixed column of this formula correspond to a fixed $y$ and different $z$ of (3a.5.4.1), and terms of a fixed row of this formula correspond to a fixed $z$ and different $y$ of (3a.5.4.1) ).
\medskip
{\bf 19.6.6.8.} For $m=2$, $u=1, \ v=0$ we have the same diagram (19.6.6.7) without the rightmost column.
\medskip
{\bf Remark 19.6.6.9.} The above matrices $S_{uvyz}$ ($v=u-1$) are not in 1 -- 1 correspondence with $A_{yu\al}$ from (19.6.4.9), but they are their (simple) linear combinations. (To write more details). 
\medskip
{\bf 19.6.6.10.} The notion of t-modular forms for the case of $N\ne0$ is similar to the case $N=0$. The only differences are the following. Let $GL(N; \p^n)$ be the subgroup of $GL(\p^n)$ consisting of elements that commute with $N$. Let $\cl$ be the set of $N$-lattices in $\p^n$, and let $\rho: GL(N; \p^n)\to GL(\g W)$ be a representation. A modular form $f: \cl \to \g W$ is defined like in 19.4. 
\medskip
{\bf 19.6.7. Case of any Schubert cell.} 
\medskip
The contents of the present subsection is a subject in development; apparently the situation is more complicated than the one described here. We think (maybe this is wrong!) that the set of Schubert cells of $tFV$ is the same as for the "classical" flag variety. Let a NL-basis  $l_*$ be fixed. We define sets $J_\al$ from 19.6.1.6 like it was made in 19.2.3 -- 19.2.5 (the order is inverted, sorry). Practically, we are going to describe $\s\circ \de(l_*)$ in absence of the analog of Definition 19.2.1 of $\de$ --- we have only a coordinate definition of $\de$. 
\medskip
First, there exists a unique sequence of numbers $\al_{1*}:=\{\al_{11}, \dots, \al_{1,k_{m+1}}\}$ such that $r=\al_{11}>\al_{12}> ... > \al_{1,k_{m+1}}\ge1$ defined as follows: 
\medskip
(19.6.7.1) $\forall \  \vk\in [1,\dots, k_{m+1}]$ we have: 
The dimension of the $\p$-linear envelope of 
$N^{m-1}(l_{\al_{1,\vk}}), N^{m-1}(l_{\al_{1,\vk}+1}), \dots, N^{m-1}(l_{r})$ is $\vk$;
\medskip
(19.6.7.2) $\forall \  \be\in [\al_{1,\vk+1}+1,\dots, \al_{1,\vk}]$ we have: the dimension of the $\p$-linear envelope of $N^{m-1}(l_{\be}), N^{m-1}(l_{\be+1}), \dots, N^{m-1}(l_{r})$ is $\vk$ (i.e. $\al_{1,\vk+1}$ is the maximal number among the numbers $\nu$ such that the dimension of the $\p$-linear envelope of $N^{m-1}(l_{\nu}), N^{m-1}(l_{\nu+1}), \dots, N^{m-1}(l_r)$ is $\vk+1$). 
\medskip
We denote the set of numbers $r=\al_{11}, \ \al_{12},\dots, \ \al_{1,k_{m+1}}$ by $J_1$, and the set of elements $l_{\al_{11}}, \dots, l_{\al_{1,k_{m+1}}}$ by $\hat l_{m+1}$ (like in 19.6.5.2). 
\medskip
The set of elements $N^{m-2}(\hat l_{m+1}), \ N^{m-1}(\hat l_{m+1})$ is linearly independent over $\p$ (the proof is the same as in [GL18], 3a.4.1). 
\medskip
Second, we define the set $J_2$ as follows. It is a decreasing sequence of numbers $J_2=\al_{2*}=\{\al_{21}>\al_{22}>...> \al_{2,k_{m}}\}$ such that $\al_{2*}\subset \{1,\dots,r\} - J_1$, and the following analogs of (19.6.7.1), (19.6.7.2) hold: 
\medskip
(19.6.7.3) $\forall \  \vk\in [1,\dots, k_{m}]$ we have: 
\medskip
The dimension of the $\p$-linear envelope of the set of elements $N^{m-2}(l_{\nu})$ for all $\nu\in\{\al_{2,\vk}, \al_{2,\vk}+1, \dots, r\}-J_1$ and of elements 

$N^{m-2}(\hat l_{m+1}), \ N^{m-1}(\hat l_{m+1})$ is $\vk+2k_{m+1}$;
\medskip
(19.6.7.4) $\forall \  \be\in [\al_{2,\vk+1}+1,\dots, \al_{2,\vk}], \be\not\in J_1$, we have: the dimension of the $\p$-linear envelope of the following elements: 
\medskip
$N^{m-2}(l_{\ga})$ for all $\ga \in [\be, \dots, r], \ \ga\not\in J_1$, and of the elements $N^{m-2}(\hat l_{m+1})$, $N^{m-1}(\hat l_{m+1})$, 
\medskip
is $\vk+2k_{m+1}$ (i.e. $\al_{2,\vk+1}$ is the maximal number among the numbers $\nu\in \{1,\dots,r\} - J_1$ such that the dimension of the $\p$-linear envelope of the following elements: 
\medskip
$N^{m-2}(l_{\ga})$ for all $\ga \in [\nu, \dots, r], \ \ga\not\in J_1$, and of the elements $N^{m-2}(\hat l_{m+1})$, $N^{m-1}(\hat l_{m+1})$, 
\medskip
is $\vk+1+2k_{m+1}$). 
\medskip
We denote the set of elements $l_{\al_{21}}, \dots, l_{\al_{2,k_{m}}}$ by $\hat l_{m}$ (like in 19.6.5.5). 
\medskip
The set of elements 
\medskip
$N^{m-3}(\hat l_{m+1}); \ N^{m-2}(\hat l_{m+1}); \ N^{m-1}(\hat l_{m+1})$;  
\medskip
$N^{m-3}(\hat l_{m}); \ N^{m-2}(\hat l_{m})$ 
\medskip
is linearly independent over $\p$ (the proof is the same as in [GL18], 3a.4.1). 
\medskip
Third, we define the set $J_3$ as follows. It is a decreasing sequence of numbers $J_3=\al_{3*}=\{\al_{31}>\al_{32}>...> \al_{3,k_{m-1}}\}$ such that $\al_{3*}\subset \{1,\dots,r\} -\{ J_1\cap J_2\}$, and the following analogs of (19.6.7.3), (19.6.7.4) hold: 
\medskip
(19.6.7.5) $\forall \  \vk\in [1,\dots, k_{m-1}]$ we have: 
\medskip
The dimension of the $\p$-linear envelope of the set of elements $N^{m-3}(l_{\nu})$ for all $\nu\in\{\al_{3,\vk}, \al_{3,\vk}+1, \dots, r\}-\{J_1\cup J_2\}$ and of elements $N^{m-3}(\hat l_{m+1}); \ N^{m-2}(\hat l_{m+1}); \ N^{m-1}(\hat l_{m+1})$;  
\medskip
$N^{m-3}(\hat l_{m}); \ N^{m-2}(\hat l_{m})$ 
\medskip
is $\vk+2k_{m}+3k_{m+1}$;
\medskip
(19.6.7.6) $\forall \  \be\in [\al_{3,\vk+1}+1,\dots, \al_{3,\vk}], \be\not\in (J_1\cup J_2)$, we have: the dimension of the $\p$-linear envelope of the following elements: 
\medskip
$N^{m-3}(l_{\ga})$ for all $\ga \in [\be, \dots, r], \ \ga\not\in (J_1\cup J_2)$, and of the elements 
\medskip
$N^{m-3}(\hat l_{m+1}); \ N^{m-2}(\hat l_{m+1}); \ N^{m-1}(\hat l_{m+1})$;  
\medskip
$N^{m-3}(\hat l_{m}); \ N^{m-2}(\hat l_{m})$ 
\medskip
is $\vk+2k_{m}+3k_{m+1}$ (i.e. $\al_{3,\vk+1}$ is the maximal number among the numbers $\nu\in \{1,\dots,r\} - (J_1\cup J_2)$ such that the dimension of the $\p$-linear envelope of the following elements: 
\medskip
$N^{m-3}(l_{\ga})$ for all $\ga \in [\nu, \dots, r], \ \ga\not\in (J_1\cup J_2)$, and of the elements 
\medskip
$N^{m-3}(\hat l_{m+1}); \ N^{m-2}(\hat l_{m+1}); \ N^{m-1}(\hat l_{m+1})$;  
\medskip
$N^{m-3}(\hat l_{m}); \ N^{m-2}(\hat l_{m})$ 
\medskip
is $\vk+1+2k_{m}+3k_{m+1}$). 
\medskip
We denote the set of elements $l_{\al_{31}}, \dots, l_{\al_{3,k_{m-1}}}$ by $\hat l_{m-1}$ (like in 19.6.5.8). 
\medskip
{\bf 19.6.7.7.} Continuing this process, we define sets $J_1,\dots, J_m$ in $\{1,\dots,r\}$, and the sets $\hat l_{m+1}, \dots, \hat l_{2}$. They are disjoint; we define $J_{m+1}:=\{1,\dots,r\}-(J_1\cup ... \cup J_m)$, $\hat l_{1}:=$ the complement to $\hat l_{m+1}\cup ... \cup \hat l_{2}$. This gives us the Schubert cell of $\de(l_*)$. 
\medskip
Definition of the Siegel set of $l_*$, i.e. analogs of the formula 19.2.7: to write. 
\medskip
{\bf 19.6.8. Duality and tensor product.} If $l_{1*}\subset \p^{n_1}$, $l_{2*}\subset \p^{n_2}$ are two NL-bases of cardinality $r_1$, $r_2$ respectively, then we can consider their tensor product $l_{1*}\otimes l_{2*}$. Let $\s\circ\de(l_{1*})=\{J_{11},\dots, J_{1,m_1+1}\}, \ \s\circ\de(l_{2*})=\{J_{21},\dots, J_{2,m_2+1}\}$. How to describe $\s\circ\de(l_{1*}\otimes l_{2*})$? Clearly it will be related with the sets $J_{1\al}\times J_{2\be}\subset \{1,\dots,r_1\}\times \{1,\dots,r_2\}$. 
\medskip
Let $\hat l_{*}$ the NL-basis dual to $l_{*}$. We denote $\s\circ\de(\hat l_{*})$ by $\{\hat J_{1},\dots, \hat J_{m+1}\}$. We have: $\hat J_{\al}=r+1-\{J_\al\}$ (maybe $r+1-\{J_{m+1-\al}\}$; maybe this depends on notations). 
\medskip
{\bf Remark.} Since the Schubert cells depend on a basis of the ambient space, we should first choose an arrangement (ordering) of the set $ \{1,\dots,r_1\}\times \{1,\dots,r_2\}$. 
\medskip
{\bf 19.6.8.1. Tensor product of t-modular forms.} Let $i=1, \ 2.$ Our objects are from 19.6.6.10, namely: $\cl_i$ are the sets of $N$-lattices of ranks $r_i$ in $\p^{n_i}$, $\rho_i: GL(N_i, \p^{n_i})\to \g W_i$ are representations, and $f_i: \cl_i\to \g W_i$ are t-modular forms. Recall that for $L_i\in \cl_i$ their tensor product $L_1\otimes L_2$ is a lattice in $\p^n$ where $n=n_1r_2+n_2r_1$. There exists the operator $N$ on $\p^n$ (recall that $N\ne0$ even if $N_1=N_2=0$). I think that we can define canonically $\g W$, $\rho: GL(N, \p^{n})\to \g W$ and a t-modular form $f_1\otimes f_2$. First, we should have $f_1\otimes f_2(L_1\otimes L_2)=f_1(L_1)\otimes f_2(L_2)$ (to define $f_1(L_1)\otimes f_2(L_2)$ ! ), and second, we should define $f_1\otimes f_2(L)$ for $L$ which are not of the form $L_1\otimes L_2$.
\medskip
{\bf 19.6.8.2. Application to the number field case.} Recall (see 5.18) that if $A$ is an abelian variety with MIQF of dimension $r$ and signature $(1,r-1)$, then there exists an abelian variety with MIQF (denoted by $\la^k(A)$) of dimension $\binom{r}{k}$ and signature $(\ \ \binom{r-1}{k-1},\binom{r-1}{k}\ \ )$ such that the lattice of $\la^k(A)$ is $\la^k(L(A))$. This means that if $f$ is a modular form on $\Cal H^3_{1,r-1}$ then (most likely) it will be possible to define $\la^k(f)$. It will be a modular form on $\Cal H^3_{\binom{r-1}{k-1},\binom{r-1}{k}}$.
\medskip
{\bf 99. Generalizations of Anderson t-motives.}
\medskip
Most generalizations use a description of an Anderson t-motive as a $\p[T]$-module with a $\tau$-action. Usually it is not convenient to consider $\tau$ as a skew map (i.e. $\tau(am)=a^q\tau(m)$\ ), hence most definitions use a formalism of a tensor product by a Frobenius map. Namely, let $A$ be a ring (it can be $\p[T]$ or a similar ring) and $\vf: A\to A$ an isomorphism called a Frobenius isomorphism. For example, for $A=\p[T]$ and $C=\sum_{i=0}^\infty c_iT^i$ we define $\vf(C)=\sum_{i=0}^\infty c_i^qT^i$.
\medskip
For a left $A$-module $M$ we can consider an $A$-module $\s^*(M):=M\underset{A}\to{\otimes}A$ where the tensor product is taken with respect to the map $\vf^{-1}: A \to A$.
\medskip
For $m\in M$ we denote the element $m\otimes 1\in \s^*(M)$ by $\bar m$, and multiplication of elements of $\s^*(M)$ by elements of $A$ will be denoted by asterisk, in order do not confuse it with multiplication of elements of $M$ by elements of $A$. All elements of $\s^*(M)$ are of the form $\bar m$ for some $m\in M$. There is a formula $a*\bar m=\overline{\vf(a)\cdot m}$.
\medskip
This means that a skew map $\tau: M\to M$ that enters in the definition of an Anderson t-motive can be considered as an $A$-module map (denoted by $\tau$ as well) $\tau: \s^*(M)\to M$, defined by the formula $\tau(\bar m)= \tau(m)$.
\medskip
If $M$ is free over $A$ of dimension $r$ then an analog of Condition 5.2.1 holds automatically. Analogs of Conditions 5.2.2, 5.2.3 are either omitted, of formulated by some other manner.
\medskip
{\bf Example 99.1.} Let $A=\overline{\n F_q}(T)$, $\vf(a)=a^q$ for $a\in \overline{\n F_q}$, $\vf(T)=T$. Let $M=V$ a finite dimensional vector space over $A$. By definition, a $\vf$-space is a bijective map $\tau: \s^*(V)\to V$.
\medskip
{\bf Example 99.2.} Let $A=\overline{\n F_q}((T))$, and let $\vf$, $V$ be the same. A bijective map $\tau: \s^*(V)\to V$ is called a Dieudonn\'e module.\footnotemark \footnotetext{In order to define localizations of a $\vf$-space at places of $\n F_q(T)$ it is necessary to use a version of the present definition of Dieudonn\'e module.}

We see that for these objects analogs of Conditions 5.2.2, 5.2.3 are omitted.
\medskip
{\bf 99.3.} A definition of an object generalizing Anderson t-motives is given in [HJ], Definition 3.1 (it is called an A-motive in [HJ]). Namely (we give a slightly simplified version), an A-motive is a pair $(M,\tau)$ where $M$ is a free $\p[T]$-module of dimension $r$ and $\tau: \s^*(M)[N^{-1}]\to M[N^{-1}]$ is an isomorphism of $\p[T][N^{-1}]$-modules (here $N=T-\th\in \p[T]$).
\medskip
We see that here it is required that $\tau$ is an isomorphism. An Anderson t-motive in the meaning of the present paper is an effective A-motive, see [HJ], Definition 3.1c. They are defined by the condition that $\tau$ comes from a $\p[T]$-homomorphism $\s^*(M)\to M$.

This definition has an advantage that Hom of A-motives is always defined, and is an A-motive.
\medskip
A natural generalization of a module over a ring is a sheaf over a scheme. Hence, instead of the ring $\p[T]$
(affine line) we consider any (projective) curve $X$ over $\p$ (or its subfields if we consider Anderson t-motives
over fields), and instead of $M$ as a free $\p[T]$-module we consider a locally free sheaf $F$
on $X$. A typical example that should be kept in mind is $X=P^1(\p)$.
\medskip
An analog of the above map $\vf: A \to A$ is a scheme automorphism of $X$ denoted by $\vf$ as well. For example, if $X=P^1(\p)$ and $a\in \p\cup\infty=P^1(\p)$ then the action of $\vf$ on $a$ is the Frobenius action: $\vf(a)=a^q$. On functions, we have $\vf(T)=T$.
Analogously to the case of modules, we define $\s^*(F):=F\underset{X}\to{\otimes}X$ where the tensor product is taken with respect to the map $\vf^{-1}: X \to X$.
\medskip
Let us consider analogs of a map $\tau: \s^*(F) \to F$. We need to prescribe its behavior not only at the point $\th\in X$ (i.e. an analog of 5.2.3), but also its behavior at a point on $X$ which is an analog of $\infty\in P^1(\p)$ (this point is denoted by $\infty$ as well). It turns out that we should consider not one map $\tau$, but a diagram of maps.
\medskip
There are various versions of the definitions. The simplest of them is
\medskip
{\bf Definition} ([D87]; notations of [G96], p. 191, (*)). A right $F$-sheaf is a diagram
$$\s^*(F_0) \overset{\al}\to{\to} F_1, \ \ \ F_0 \overset{\be}\to{\to} F_1$$
where $F_0, \ F_1$ are locally free sheaves on $X$, $\al$ and $\be$ are inclusions such that the supports of their cokernels consist of one point.
\medskip
The support of Coker $(\al)$ is called the zero of $F$-sheaf (it is an analog of $\th$ of 5.2.3), and the support of Coker $(\be)$ is called the pole of $F$-sheaf (it is an analog of $\infty\in P^1(\p)$).
\medskip
Really, Drinfeld considers the case dim Coker $(\al)$ = dim Coker $(\be)=1$, i.e. analog of the case $n=1$, because he did not need an analog of Anderson t-motives. Also, he considers a relative case, i.e. a case when the whole diagram is over a base scheme $S$.
\medskip
The notions of $F$-sheaf and its generalizations were introduced in order to prove Langlands conjecture for $GL_2$ (Drinfeld), $GL_r$ (L. Lafforgue) et al. Nevertheless, they can be used for solutions of much more elementary problems.
\medskip
As an elementary example of application of the notion of $F$-sheaf, we can mention a theorem of existence of an Anderson t-motive having a complete multiplication\footnotemark \footnotetext{Complete multiplication of Anderson t-motives is an analog of the complex multiplication of abelian varieties. See [A87] for an analog of the main theorem of complex multiplication for this case.} over a field $\g K$ with a given complete multiplication type. The idea is to consider a curve $X$ corresponding to $\g K$ and its Picard variety $\Pic_0(X)$. The map $fr - Id$ is an algebraic isogeny on $\Pic_0(X)$, hence it is surjective. A complete multiplication type defines a divisor on $X$ and hence an element of $\Pic_0(X)$. Its preimage with respect to the map $fr - Id$ defines a sheaf corresponding to a desirable t-motive of complete multiplication. See [GL07], Theorem 12.6 for details.\footnotemark \footnotetext{This proof was indicated to the authors by V. Drinfeld.}
\medskip
%\newpage
{\bf References}
\medskip
[A86] Anderson, Greg W. t-motives, Duke Math. J. 53 (2) (1986) 457 -- 502.
\medskip
[A87] Anderson, Greg W. On a question arising from complex multiplication theory. Galois representations and arithmetic algebraic geometry (Kyoto, 1985/Tokyo, 1986), 221 -- 234, Adv. Stud. Pure Math., 12, Amsterdam, 1987.
\medskip
[AT90] Anderson, G.W., Thakur, D.S. Tensor powers of Carlitz module and zeta values. 
Annals of Mathematics, 132 (1990), p. 159 -- 191.
\medskip
[A00] Anderson, Greg W. An elementary approach to $L$-functions mod $p$. J. Number theory 80 (2000), no. 2, 291 -- 303.
\medskip
[B02]  B\"ockle, Gebhard. Global $L$-functions over function fields. Math. Ann. 323 (2002), no. 4, 737 -- 795.
\medskip
[B05] B\"ockle, Gebhard. Arithmetic over function fields: a cohomological approach. Number fields and function fields -- two parallel worlds, 1 -- 38, Progr. Math., 239, Birkhäuser Boston, Boston, MA, 2005.
\medskip
[B12] B\"ockle, Gebhard. Cohomological theory of crystals over function fields and applications.
In "Arithmetic Geometry in Positive Characteristic". Advanced Courses in Mathematics. CRM Barcelona. Birkh\"auser Verlag, Basel (2012).
\medskip
[BH] Bornhofen, Matthias; Hartl, Urs. Pure Anderson motives over finite fields. J. Number Theory 129 (2009), no. 2, 247--283.
\medskip
[B] Brion, Michel. Lectures on the geometry of flag varieties. Topics in cohomological studies of algebraic varieties, 33–85,
Trends Math., Birkhäuser, Basel, 2005.
\medskip
[De] Deligne, Pierre. Travaux de Shimura. S\'eminaire Bourbaki, 23\`eme ann\'ee (1970/1971), Exp. No. 389, pp. 123--165
Lecture Notes in Math., Vol. 244 Springer-Verlag, Berlin-New York, 1971
\medskip
[D76] V.G. Drinfeld, Elliptic modules, Math. USSR Sb. 4 (1976) 561 -- 592.
\medskip
[D87] V.G. Drinfeld, Functional Anal. Appl. 21 (1987), no. 2, 107 -- 122.
\medskip
[EGL] Ehbauer S.; Grishkov, A.; Logachev, D. Calculation of $h^1$ of some Anderson t-motives. J. Algebra Appl. 2022, vol. 21, no. 1, Paper No. 2250017, 31 pp. 
https://arxiv.org/pdf/2006.00316.pdf 
\medskip
[Ga] Gardeyn, F. A Galois criterion for good reduction of $\tau$-sheaves. J. Number Theory 97 (2002), no. 2, 447--471.
\medskip
[GR] E.-U. Gekeler, M. Reversat. Jacobians of Drinfeld modular curves. J. Reine Angew. Math. 476 (1996), 27 -- 93.
\medskip
[G92] Goss, David. $L$-series of $t$-motives and Drinfeld modules. The arithmetic of function fields (Columbus, OH, 1991), 313 -- 402, Ohio State Univ. Math. Res. Inst. Publ., 2, de Gruyter, Berlin, 1992.
\medskip
[G95] Goss, David.
The adjoint of the Carlitz module and Fermat's last theorem.
Finite Fields Appl. 1 (1995), no. 2, 165 -- 188.
\medskip
[G96] Goss, David. Basic Structures of Function Field Arithmetic, Springer-Verlag, Berlin, 1996.
\medskip
[GL07] Grishkov, A.; Logachev, D. Duality of Anderson t-motives. https://arxiv.org/pdf/0711.1928.pdf
\medskip
[GL09] Grishkov, A.; Logachev, D. Anderson t-motives and abelian varieties with MIQF: results coming from an analogy. J. Algebra Appl. 21 (2022), no. 9, Paper No. 2250171. 
https://arxiv.org/pdf/0907.4712.pdf
\medskip
[GL16] Grishkov, A.; Logachev, D. Resultantal varieties related to zeroes of $L$-functions of Carlitz modules. Finite Fields Appl. 38 (2016), 116 -- 176.

arxiv.org/pdf/1205.2900.pdf
\medskip
[GL17] Grishkov, A.; Logachev, D. Lattice map for Anderson t-motives: first approach.
J. Number Theory, 180 (2017), 373 -- 402. https://arxiv.org/pdf/1109.0679
\medskip
[GL18] Grishkov, A.; Logachev, D.  Lattice of the dual of an Anderson $t$-motive in terms of Siegel matrices of its flag variety.

https://arxiv.org/pdf/1812.11576.pdf
\medskip
[GL21] Grishkov, A.; Logachev, D. $h^1 \ne h_1$ for Anderson t-motives. J. of Number Theory. 2021, vol. 225, p. 59 -- 89.
https://arxiv.org/pdf/1807.08675.pdf
\medskip
[GLZ22] Grishkov, A.; Logachev, D.;  Zobnin, A. L-functions of Carlitz modules, resultantal varieties and rooted binary trees - I. J. of Number Theory,  2022, vol. 238, p. 269 -- 312.
https://arxiv.org/pdf/1607.06147.pdf.
\medskip
[GL23] Grishkov A., Logachev D. $h^1,\ h_1$ of Anderson t-motives, systems of affine equations and non-commutative determinants. 

https://arxiv.org/pdf/2302.13480.pdf
\medskip
[GL24] Grishkov A., Logachev D. Non-injectivity of the lattice map for non-mixed Anderson t-motives, and a result towards its surjectivity. J. of Algebra and its applications, 2026 (to appear). 
https://arxiv.org/pdf/2405.17162.pdf
\medskip
[H] Urs Hartl, Uniformizing the Stacks of Abelian Sheaves.

http://arxiv.org/abs/math.NT/0409341
\medskip
[HJ] Hartl U.; Juschka A.-K. Pink's theory of Hodge structures and the Hodge conjecture over function fields.  "t-motives: Hodge structures, transcendence and other motivic aspects", Editors G. B\"ockle, D. Goss, U. Hartl, M. Papanikolas, European Mathematical Society Congress Reports 2020, and https://arxiv.org/pdf/1607.01412.pdf
\medskip
[K] Kauers, M. The Holonomic Toolkit. 

https://www3.risc.jku.at/publications/download etc.
\medskip
[KP] Kauers, Manuel; Paule, Peter. The concrete tetrahedron. Texts Monogr. Symbol. Comput.
Springer Wien, New York, Vienna, 2011, x+203 pp.
\medskip
[L] Lafforgue, Vincent. Valeurs sp\'eciales des fonctions $L$ en caract\'eristique $p$.  J. Number Theory  129  (2009),  no. 10, 2600 -- 2634.
\medskip
[Mu78] Mumford, David. An algebro-geometric construction of commuting operators and of solutions to the Toda lattice equation, 
Korteweg deVries equation and related nonlinear equation.
Proceedings of the International Symposium on Algebraic Geometry (Kyoto Univ., Kyoto, 1977), pp. 115 -- 153
Kinokuniya Book Store Co., Ltd., Tokyo, 1978
\medskip
[NP] Namoijam Ch., Papanikolas M. Hyperderivatives of periods and quasi-periods for Anderson t-modules. https://arxiv.org/pdf/2103.05836.pdf
\medskip
[Pa23] Papikian, Mihran. Drinfeld modules. Grad. Texts in Math., 296 Springer, Cham, 2023. xxi+526 pp.
\medskip
[P] Richard Pink, Hodge structures over function fields. Universit\"at Mannheim.
Preprint. September 17, 1997.
\medskip
[R16] Richarz, Timo. Affine Grassmannians and Geometric Satake Equivalences. 
Int. Math. Res. Not. IMRN 2016, no. 12, 3717 -- 3767
\medskip
[Sh] Shimura, Goro. On analytic families of polarized abelian varieties and automorphic functions.
Annals of Math., 1 (1963), vol. 78, p. 149 -- 192
\medskip
[Tae] Taelman, L. Special $L$-values of Drinfeld modules. Ann. of Math. (2) 175 (2012), no. 1, 369 -- 391.
\medskip
[Tag95] Taguchi, Yuichiro. A duality for finite $t$-modules.
J. Math. Sci. Univ. Tokyo 2 (1995), no. 3, 563 -- 588.
\medskip
[Tag96] Taguchi, Yuichiro. On $\vf$-modules. Journal of number theory 60, 124 -- 141 (1996).
\medskip
[TW] Taguchi, Yuichiro; Wan, D. $L$-functions of $\phi$-sheaves and Drinfeld modules. 
J. Amer. Math. Soc. 9 (1996), no. 3, 755--781.
\enddocument